\documentclass[12pt,letterpaper,reqno]{amsart}

\usepackage{amssymb}
\usepackage{amsmath}
\usepackage{amsthm}
\usepackage{color}
\usepackage{eucal}
\usepackage{graphicx}
\graphicspath{ {pics/} }

\usepackage{caption}
\usepackage{subcaption}

\usepackage{wrapfig}
\usepackage{floatflt}

\usepackage{wrapfig}

\newcommand{\R}{\mathbb R}

\newcommand{\sech}{\textmd{sech}}

\newcommand{\cR}{\mathbb{R}}

\newtheorem{theorem}{Theorem}[section]

\theoremstyle{remark}
\newtheorem{remark}[theorem]{Remark}

\title[1d bi-NLS with mixed dispersion]{Dynamics of solutions in the 1d bi-harmonic \\
nonlinear Schr\"odinger equation}

\author[C. Klein]{Christian Klein}
\address{Institut de Math\'ematiques de Bourgogne, UMR 5584\\
Institut Universitaire de France,\\
Universit\'e Bourgogne Europe, 
9 avenue Alain Savary, 21078 Dijon
Cedex, France} 
\email{Christian.Klein@u-bourgogne.fr}

\author[I. Petrenko]{Iryna Petrenko}
\address{Department of Mathematics \& Statistics\\
Florida International University,  Miami, FL 33199, USA}
\curraddr{}
\email{ipetr006@fiu.edu}

\author[S. Roudenko]{Svetlana Roudenko}
\address{Department of Mathematics \& Statistics\\
Florida International University,  Miami, FL 33199, USA}
\curraddr{}
\email{sroudenko@fiu.edu}

\author[N. Stoilov]{Nikola Stoilov}
\address{Institut de Math\'ematiques de Bourgogne, UMR 5584\\
Universit\'e Bourgogne Europe, 
9 avenue Alain Savary, 21078 Dijon
Cedex, France} 
\email{Nikola.Stoilov@u-bourgogne.fr}

\date{}


\begin{document}

\begin{abstract}
We consider the one dimensional 4th order, or bi-harmonic, nonlinear Schr\"odinger (NLS) equation,  namely, $i 
u_t - \Delta^2 u - 2a \Delta u + |u|^{\alpha} u = 0, ~ x,a \in \R$, $\alpha>0$, 
and investigate the dynamics of its solutions for various powers of 
$\alpha$, including the ground state solutions and their 
perturbations, leading to scattering or blow-up dichotomy when $a \leq 0$, or to a trichotomy when $a>0$. Ground 
state solutions are numerically constructed, and their stability is 
studied, finding that the ground state solutions may form two branches, 
stable and unstable, which dictates the long-term behavior of 
solutions.  Perturbations of the ground states on the unstable branch 
either lead to dispersion or the jump to a stable ground state. 
In the critical and supercritical cases, blow-up in 
finite time is also investigated, and it is conjectured that the 
blow-up happens with a scale-invariant profile (when $a=0$) 
regardless of the value of $a$ of the lower dispersion. 
The blow-up rate is also explored.  
\end{abstract}

\maketitle

\tableofcontents

\section{Introduction}

We consider the 4th order nonlinear Schr\"odinger (NLS) equation, often referred to as the bi-harmonic NLS equation (especially when no lower dispersion present, i.e., when $a=0$): 
\begin{equation}\label{biNLS}
\text{(bi-NLS)} 
\qquad \qquad 
i\, u_t - \Delta^2 u - 2a \Delta u + |u|^{\alpha} u = 0, \qquad x \in \R^d, \quad t \in \cR, \qquad
\end{equation}
where $u(t,x)$ is a complex-valued function, $a \in 
\mathbb R$, $\Delta$ is the standard Laplacian, and the nonlinearity power $\alpha>0$. In this work we examine the one dimensional case, $d=1$, hence, the Laplacian is equivalent to $\partial^2_{x}$, consequently, $\Delta^2 \equiv \partial^4_x$. This equation is a higher dispersion generalization of the well-known nonlinear Schr\"odinger equation:
\begin{equation}
\label{NLS}
\text{(NLS)} \qquad \qquad
i\, u_t + \Delta u + |u|^{\alpha} u = 0, \qquad x \in \cR^d, \qquad t \in \cR. \qquad
\end{equation}

\subsection{Background}

A first mentioning of the bi-harmonic NLS model was in the early 90s by Karpman \cite{K1991}, and Karpman and Shagalov \cite{KS1991}, where the influence of higher order dispersion was included into the NLS equation to model intense laser beam propagation. Lately, the bi-harmonic NLS model has been increasingly attracting attention as the quartic solitons have certain favorable properties such as flattening or stabilization in applications and experiments. A recent experimental work in silicon photonic crystal waveguides for the first time produced pure quartic solitons on a chip  \cite{Nature}, where the leading order dispersion in the model was quartic (instead of a typical quadratic dispersion as in the NLS model \eqref{NLS}). Evolution from Gaussian data into pure quartic solitons and other features were further followed up in \cite{Tam2019}, and a more general model with a mixed dispersion such as \eqref{biNLS} is discussed in \cite{Tam2020}. Even more recently, considering that the quartic temporal solitons have been experimentally achieved, 
the spatio-temporal solitons (SKS), or light bullets, were described by the quartic dispersion \cite{LightBullets2024}. For a review and related work, see \cite{Aceves2023}, \cite{PA2024}, \cite{HA2024}, \cite{Mal2024}. Therefore, long-term behavior of solutions and stability of solitary waves from the mathematical point of view are timely questions to investigate.

During their lifespan, solutions $u(t)$ to \eqref{biNLS} conserve mass and energy (or Hamiltonian):
\begin{equation}\label{MC}
M[u(t)]=\int_{\cR^d} |u(t)|^2 \, dx \equiv M[u(0)]
\end{equation}
and
\begin{equation}\label{EC}
E[u(t)]=\dfrac{1}{2}\int_{\cR^d} |\Delta u(t)|^2 \; dx -a \int_{\cR^d} |\nabla u(t)|^2 \; dx - \dfrac1{\alpha+2} \int_{\cR^d} |u(t)|^{\alpha+2} \; dx \equiv E[u(0)].
\end{equation}

Similar to the NLS equation, the bi-harmonic NLS has time, space and phase invariances; the one, which is especially useful in the evolution equations, is {\it the scaling invariance}, which states that an appropriately rescaled version of the original solution is also a solution of the equation. For the equation \eqref{biNLS} due to the different dispersion terms, there is no simple suitable symmetry like that, however, if one considers the lower dispersion absent ($a=0$, a pure bi-harmonic NLS), 
then the scaling is
\begin{equation}
	u_\lambda(t, x)=\lambda^{\frac{4}{\alpha}} u(\lambda^4 t, \lambda x).
	\label{scaling}
\end{equation}
This symmetry makes a specific Sobolev norm $\dot{H}^s$ invariant, i.e.,
\begin{equation*}
\|u(0,\cdot,\cdot) \|_{\dot{H}^s}=\lambda^{\frac{4}{\alpha}+s-\frac{d}{2}} \|u_0\|_{\dot{H}^s},
\end{equation*}
and the index $s$ gives rise to the critical-type classification of equations.
For the bi-harmonic NLS equation \eqref{biNLS} (with $a=0$) the critical index is 
$$
s=\frac{d}2-\frac4{\alpha},
$$
and when $s=0$ (or $\alpha=8/d$) the equation \eqref{biNLS} is $L^2$-critical, when $s<0$ (or $\alpha<8/d$) it is subcritical and $s>0$ (or $\alpha>8/d$) is supercritical. 
While the general equation \eqref{biNLS} does not have the scaling invariance, we nevertheless use the values of the above scaling index $s$ for its critical-type classification. 

The local well-posedness of the initial value problem \eqref{biNLS} with $u(0,x)=u_0$ in the energy space $H^2(\mathbb{R}^d)$ was established by Ben-Artzi, Koch \& Saut \cite{BAKS2000}. 
Papanicolaou, Fibich \& Ilan obtained sufficient conditions for global $H^2$ solutions in \cite{Fibich2002} for some cases of cubic and quintic power and did asymptotic analysis (and numerical simulations). Pausader obtained global solutions in the energy-critical and some subcritical cases and investigated scattering or ill-posedness, see \cite{Pausader2007}, \cite{Pausader2009_1}, \cite{Pausader2009_2}. 
Improvements and clarifications about the well-posedness was done by Dinh in \cite{Dinh2018} via the Strichartz esimates corresponding to the quartic flow (developed in \cite{BAKS2000}). For $\alpha \geq 1$ the local well-posedness is also known in $H^1$ in dimension one, and for any $\alpha >0$ in weighted subspaces of $H^s(\mathbb R^d)$ with certain conditions on the initial data for some $s>s_0>0$ via arguments not involving Strichartz estimates, see the work of the second and third authors in  \cite{PRR}.

Provided there exists a suitable local well-posedness, one can obtain {\it global} well-posedness in the energy space $H^{2}(\mathbb{R}^d)$ in the subcritical case ($s<0$ or $\alpha <8/d$) via the corresponding Gagliardo-Nirenberg inequality and the energy conservation \eqref{EC}. The same argument will show the global existence in the critical ($s=0, \alpha = 8/d$) case, provided a bounded condition on the mass holds (i.e., mass less than that of a ground state, which we define below), 
see, for example, Fibich, Ilan \& Papanicolaou \cite{Fibich2002}.
In the supercritical case ($s>0, \alpha >8/d$) an argument as in 
Holmer \& Roudenko \cite{HR2008, HR2007} using the invariant quantities (expressed via mass, energy and such) gives  the dichotomy for global existence and scattering vs.\ blow-up for radial functions in 2d and higher, see \cite{BL2017}, also \cite{Dinh2018}, \cite{Dinh2021}, for a more general case, see \cite{PRR}.

\subsection{Solitary waves}

The 4th order NLS equation \eqref{biNLS} has a family of (standing) solitary waves, called waveguide solutions in nonlinear optics, 
\begin{equation}\label{Eq:SW}
u(t,x) = e^{ibt} \, Q(x), 
\end{equation}
with $Q(x) \to 0$ as $|x| \to + \infty$, and we take $b>0$ in this paper. 
Here, $Q$ is a 
{\it ground state} solution in $H^2(\cR^d)$, the energy space, of the nonlinear elliptic equation 
\begin{equation}\label{E:groundstate}
\Delta^2 Q + 2a\,\Delta Q +  b\, Q  - |Q|^{\alpha} Q = 0.
\end{equation}

In the pure quartic case, $a=0$, a simple way to define a ground state is as an optimizer of the Gagliardo-Nirenberg inequality (or equivalently, of the corresponding Weinstein functional): for $u \in H^{2}(\mathbb R^d)$ and $\frac8{d} < \alpha<\frac8{d-4}$ (or $0<s<2$), 
$$
\|u\|_{L^{\alpha+2}}^{\alpha+2} \leq C_{GN} \, \|\Delta u\|_{L^2}^{\frac{\alpha d}4} \|u\|_{L^2}^{2-\frac{\alpha}4 (d-4)}, 
$$ 
where $C_{GN}$, an optimal constant, depends on the power $\alpha$ and the dimension $d$, e.g., see \cite{BL2017}. 
Uniqueness of ground states in general is not known, however, since we consider only even powers of $\alpha$ (in one dimension), the ground state can be chosen to be radially symmetric, real-valued, and continuous (actually, $Q \in H^\infty(\mathbb R)$), with $Q(0)>|Q(r)|$, $r = |x| >0$ (e.g., see appendix A in \cite{BL2017} or Prop. 3.6 in \cite{BCSN2018}). We emphasize that the ground states in this case are non-monotonic, non-positive, and oscillate around the $x$-axis, see for instance, \cite{Fibich2002} and an example on the right of Fig. \ref{QA}.    
We note that in this pure quartic case, the following scaling is useful in our simulations:
\begin{equation}\label{Qscaling}
Q_b(x) = b^{1/\alpha} Q(b^{1/4} x),
\end{equation}
which produces a family of solutions ${Q_b}$, provided $Q$ is a solution of \eqref{E:groundstate} with $a=0$.

Furthermore, it is sufficient to fix one of the parameters (for example, say, $a$) and consider dependence on the other (in this case would be $b$), then to understand behavior for another value of $a$, one would just rescale the ground state and consider the equation \eqref{E:groundstate} replacing $b$ with $b/a^2$, see Remark \ref{R:1}.

In the general case ($a \in \mathbb R$), a ground state is defined as the least energy solution of 
some action functional, which it minimizes (and typically constrained either under the mass, the $L^2$-norm, or the potential energy, the $L^{\alpha+2}$-norm). Thus, define a quadratic form
$$
q_{a,b}(u) = \|\Delta u\|^2_{L^2} - 2a \|\nabla u\|^2_{L^2} +b\|u\|^2_{L^2}, 
$$
with the energy functional corresponding to the stationary equation \eqref{E:groundstate}:
$$
E_{a,b}(u) = \frac12 q_{a,b}(u) - \frac1{\alpha+2} \|u\|^{\alpha+2}_{L^{\alpha+2}}.
$$
Note that 
$$
q_{a,b}(u) = \int g_{a,b}(|\xi|) |\hat u(\xi)|^2 \, d\xi,
$$ 
where 
$$
g_{a,b}(|\xi|) = |\xi|^4-2a|\xi|^2+b  = (|\xi|^2-a)^2 +b-a^2.
$$  
When $a<0$, the multiplier $g_{a,b}(x)$ is an upward parabola with the minimum at $x=0$, $g_{a,b}(0)=b$; it is decreasing on $(-\infty,0)$ and increasing $(0,\infty)$; thus, the multiplier is monotone and positive, resembling the pure quartic, scaling-invariant case $a=0$. Therefore, due to the negative coefficient $a$,  the lower dispersion does not interfere with the higher order dispersion, and in a sense `helps' solutions to behave similar to the scale-invariant case. This makes it easier to use or identify some thresholds in global behavior via the conserved quantities of the ground state (such as energy or mass, e.g., see \cite{BL2017}, \cite{BCSN2018}, \cite{BN2015}). 
Furthermore, it was shown in \cite[Thm 1.1]{BN2015} that if $a \leq -\sqrt{b}$, then any least energy solution (i.e., a ground state), does not change sign, is radially symmetric (around some point) and is strictly radially decreasing (see our numerical confirmation in Figures \ref{QA} and \ref{QC}), some positive explicit solutions are discussed in Section \ref{S:GSexact}.   

For $a>0$ 
it is easy to notice that $q_{a,b}$ is positive-definite if and only 
if $a^2<b$ (also, observe that the parabola $g_{a,b}(x)$ has 3 local 
minima (at $x = -a, 0, a$, and is no longer monotonic on each side of 
$x=0$).  Under this condition, following \cite{LW2021}, we minimize $q_{a,b}$ under the potential energy $L^{\alpha+2}$ constraint:  
\begin{equation}\label{E:PE-constrain}
R_{a,b}(\alpha):=\inf \{ q_{a,b}(u): u \in H^2(\mathbb R^d) \setminus \{0\}, \|u\|_{L^{\alpha+2}}=1 \} 
\equiv \inf \frac{ q_{a,b}(u)}{\|u\|^2_{L^{\alpha+2}}}.
\end{equation}
Then a ground state is the function $Q$, on which $\inf R_{a,b}(\alpha)$ is attained. As noted in \cite{LW2021}
the value of the least energy among all non-trivial solutions of \eqref{biNLS} is characterized as 
$$
\inf_{u} \sup_{t \geq 0} E_{a,b} (tu) \equiv (R_{a,b}(\alpha))^{\frac{\alpha+2}{\alpha}}.
$$
These minimizers correspond (up to multiplication by a positive factor) to non-trivial solutions of \eqref{biNLS}, 
where the {\it least energy} value is attained. 
While the energy is minimal on such solutions $Q$ (which are often 
referred to as ``minimum action solutions"), the constraint in \eqref{E:PE-constrain} is not on the mass (or $L^2$ norm) but rather on the potential energy. 
Theorem 1.3 in \cite{FJMM} 
shows that these minimizers correspond to the minimizers of the 
energy under a fixed mass constraint. 
Thus, regardless of the constraint, the set of minimizers in this case are the same, and therefore, we compute ground states as critical points of the energy. 

Summarizing, for the purpose of this work in 1d, ground states can be chosen radially symmetric and positive for $a\leq -\sqrt b$, and real-valued but sign-changing for $-\sqrt b < a < \sqrt b$, with oscillatory behavior as $|x| \to \infty$, see e.g., \cite{BF2011}, \cite{LW2021} and  Section \S \ref{S:GS} for examples. 
Stability of the set of minimizers (and its connection with ground state solutions) and other properties have been investigated starting from the work of Albert \cite{Albert} (who also found an explicit positive ground state solution), for some recent progress refer to \cite{BCSN2018}, \cite{BCGJ2019}, \cite{FJMM}, \cite{NataliPastor2015}, \cite{dA2023} and references therein.  
Unlike the 1d NLS equation, which has explicit ground state solutions 
for any $\alpha$ (in terms of the {\it sech} function), explicit 
solutions of the elliptic problem \eqref{E:groundstate} are known 
only in a few specific cases. We mention some of them in Section \ref{S:GS}. 
\smallskip

Numerical simulations of solutions to the bi-NLS equation 
\eqref{biNLS}, including solitary wave solutions to \eqref{E:groundstate}, go back to work of Karpman and Shagalov work in \cite{KS1991}, and then more thorough investigations by Fibich, Ilan \& Papanicolaou in  \cite{Fibich2002} and their follow-up work, especially, the collapsing or blow-up solutions in critical and supercritical cases \cite{BFM2010,BFM2010b,BF2011}. 
Unlike the NLS equation, the bi-harmonic NLS equation does not have a convenient or rather simple virial identity, which in a standard NLS typically gives a straightforward proof of existence of collapse or blow-up solutions. Numerical investigations of finite time blow-up in the bi-harmonic NLS equation was initially done 
by Karpman and Shagalov in \cite{KS1991}, then by Fibich et al \cite{Fibich2002} and their follow-up work
in \cite{Fibich2002}, \cite{BFM2010}, \cite{BF2011}. 
The breakthrough for proving analytically the existence of blow-up in the bi-harmonic NLS equation ($d \geq 2$, $\alpha \geq 8/d$) was done by Boulenger and Lenzmann in \cite{BL2017}, see further progress in \cite{BCGJ2019, LW2021}.
The question of existence of blow-up is entirely open in one dimension, as is the finite time blow-up in a pure quartic case, $a=0$, in dimension two and higher, or if there are any blow-up solutions when $a>0$. Investigating this in the {\it one-dimensional} case as well as global behavior of solutions and dynamics of solitary waves in the $1d$ are the goals of this paper. 

\subsection{Main results}

In this paper we consider the 1d bi-harmonic NLS with several even 
nonlinear powers $\alpha$ that correspond to the $L^2$-subcritical, 
critical and supercritical cases (if $a=0$), and address the question 
of solutions behavior globally, 
specifically, the dynamics of solitary waves and their perturbations. We are especially interested in the behavior of sign-changing ground state solutions, their stability (or not) when $\alpha \leq 8$ and stable blow-up when $\alpha \geq 8$. 
\smallskip

Our first, and most surprising, observation is that the typical dichotomy in solutions behavior (scattering vs. finite time blow-up) does not necessarily hold in the mixed dispersion equation \eqref{biNLS}. 
More precisely, in 1d there exist positive $\alpha_\ast$ and $\alpha^\ast$ with $2 \leq \alpha_\ast < \alpha^\ast <10$ such that for any power $\alpha \in (\alpha_\ast,\alpha^\ast)$ there are two branches of ground state solutions, with one of them being a {\it stable} branch and another one {\it unstable}, see Fig. \ref{F:ME-subcritical}, \ref{F:ME-crit+supercrit}. Perturbations of solitary waves with mass slightly larger than that of the unstable ground state will jump to the stable branch (the one with the lower energy), exhibiting an oscillatory asymptotic approach to the stable ground state solution (rescaled and shifted), this can also be thought as `scattering' to the stable branch; perturbations with slightly lower mass of the unstable ground state will disperse away (for examples, see Fig. \ref{F:superG}, \ref{F:sech}). 
Perturbations of the stable branch with small deviations show a stable asymptotic oscillatory behavior below or above the mass of the ground state (e.g., Fig. \ref{F:AQ0-alpha6-mass2.9} bottom row or Fig. \ref{F:alpha8branch} left column). 

{\it Remark.} This behavior and branching has resemblance to the combined nonlinearity case recently shown in \cite{CKS}, see also    
\cite{BJM2009}, \cite{CJS2009} for the definition and description of ground states as minimizers in that  case. 

We next recall the blow-up alternative in the energy-subcritical cases, $s<2$ (e.g., \cite{BL2017}): either the solution $u\in C^0([0,T], H^2(\mathbb R^d))$ of \eqref{biNLS} extends to all times $t \geq 0$ or $\displaystyle \lim\limits_{t \to t*} \|\Delta u(t)\|_{L^2} = + \infty$.   
For $\alpha \geq 8$ (in 1d) the local theory gives a lower bound on 
the blow-up rate \eqref{E:rate-lb}, see Sections 
\ref{S:rate-critical} and \ref{S:rate-supercritical}. Furthermore, 
when $\alpha = 8$ (the critical case of \eqref{biNLS}), similar to 
the critical NLS case, the convergence to the self-similar blow-up to 
a (rescaled) blow-up profile is slow (compared to the supercritical 
case $\alpha >8$, where the convergence is exponentially fast), which 
affects the blow-up rate. We discuss that in \S \ref{S:blowup} with 
numerical confirmations of the blow-up profile and the rate, which we 
compare to  the conjectured rate in \cite{BFM2010}. 
\\

The main results of our studies, including numerical simulations, confirm the following three conjectures about solutions to the equation \eqref{biNLS} in 1d: 

\textbf{Conjecture I:} (1d subcritical case, $\alpha<8$)

\begin{figure}[!htb]
\includegraphics[width=0.54\hsize,height=0.33\hsize]{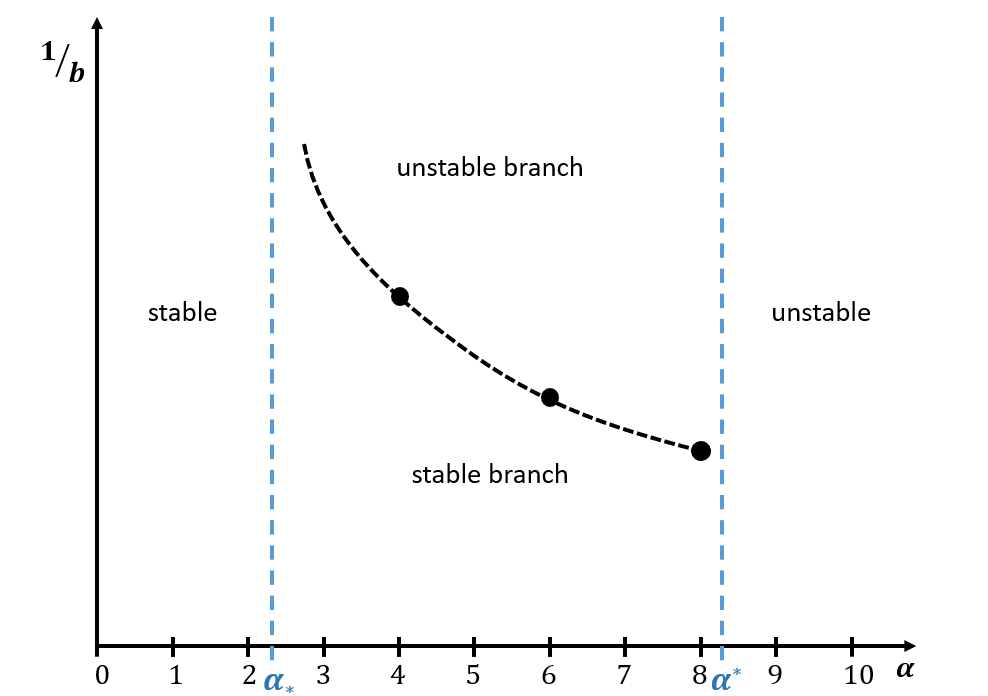}
\includegraphics[width=0.45\hsize,height=0.33\hsize]{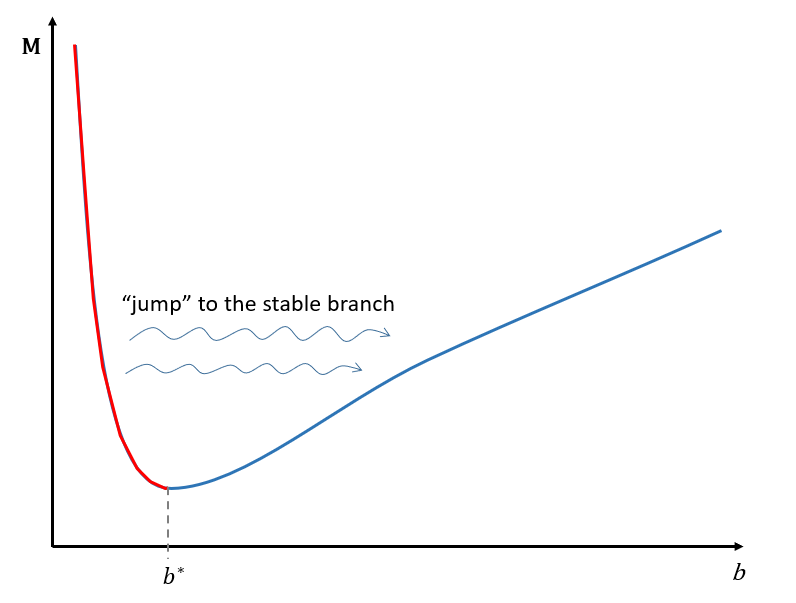}
\caption{\footnotesize  Schematic representation of stability of ground states (left) for different powers of nonlinearity $\alpha$ as stated in Conjectures. Below some $\alpha_\ast \leq 2$ all ground states are stable (perturbations asymptotically approach a rescaled ground state). Above some $\alpha^\ast > 8$ all ground states are unstable (perturbations either radiate away or blow up). In between, there are two branches of ground states: for $b>b^\ast$ a stable branch (on the graph $1/b^\ast$ curve is indicated), and for $b<b^\ast$ it is unstable (perturbations `jump' to a stable branch of ground states as shown on the right for the subcritical case, Conjecture I part (2); for critical case see Figure \ref{F:scheme}). }
\label{F:scheme2}
\end{figure}

There exist $2 \leq \alpha_\ast < 4$ 
and $8 < \alpha^\ast < 10$ such that 
\begin{enumerate}
\item 
The ground state solutions are asymptotically stable in the subcritical case for small powers $\alpha \leq \alpha_*$ (no branching of ground states occurs). 

\item
The ground state solutions of \eqref{E:1dGS} form two branches of {\it stable} and {\it unstable} ground states for power $\alpha \in (\alpha_\ast,\alpha^\ast)$ (in particular, in the subcritical range $\alpha \in (\alpha_\ast,8)$). Fixing such power $\alpha <8$, there exists a value $b^\ast$ such that the ground state solutions of \eqref{E:1dGS} with $b < b^\ast$ belong to the unstable branch and with $b>b^\ast$ belong to the stable branch. Furthermore,
\begin{enumerate}
\item
the perturbations of the unstable branch lead to either (i) jumping onto the stable branch or (ii) dispersing away (in other words, scattering to zero); 
\item
the perturbations of the stable branch are asymptotically stable, i.e., approach a rescaled version of a (stable branch) ground state.  
\end{enumerate}

\item 
The long time behavior of solutions to the subcritical 4th order mixed dispersion NLS equation \eqref{biNLS} for initial data in the 
Schwartz class is characterized by the appearance of ground states 
plus radiation (in accordance with the \emph{soliton resolution 
conjecture}).
\end{enumerate}
We show an illustration of Conjecture I in Figure \ref{F:scheme2}. 

In the critical case, in \cite{BFM2010} it was stated that sufficiently localized initial data with a mass larger than 
the mass of the ground state blows up in finite time and disperses if it is below the ground state mass,  which resembles the standard NLS equation.  We investigate this further and find that in the case of mixed dispersion, the blow-up may not happen for slightly supercritical mass (of a ground state), instead it happens for larger values; similarly, initial data with mass just slightly below the mass of a ground state may not disperse (or scatter down) to zero, but instead approach asymptotically a different final state. This happens since the scaling invariance is broken (by having two different dispersions), and that produces {\it a gap} in the typical dichotomy (scattering vs. blow-up) of solutions; thus, forming a {\it trichotomy} (scattering to zero/linear solution or dispersing away, scattering to or asymptotically approaching a stable soliton, and finite time blow-up).  Specifically, 
\smallskip

\textbf{Conjecture II:} (1d critical case, $\alpha=8$)

Let $u_0 \in \mathcal{S}(\mathbb{R})$ be 
the Schwartz class of smooth rapidly decaying functions and let $Q^{(a)}$ denote a ground state solution of \eqref{E:1dGS} for a given $a$ (varying with a positive $b$).  

\begin{enumerate}
\item
Similar to the subcritical case, the ground state solutions of \eqref{E:1dGS} form two branches of {\it stable} and {\it unstable} ground states, i.e., there exists a value $b^\ast>0$ such that the ground state solutions of \eqref{E:1dGS} with $0<b < b^\ast$ belong to the unstable branch and with $b>b^\ast$ belong to the stable branch. Furthermore,
\begin{enumerate}
\item
small perturbations of the unstable branch lead to either (i) jumping onto the stable branch or (ii) dispersing away (in other words, scattering to zero); 
\item
small perturbations of the stable branch are asymptotically stable, i.e., approach a rescaled version of a (stable branch) ground state;
\item
larger amplitude perturbations of the stable branch lead to blow-up solutions.   
\end{enumerate}
\smallskip

\item[(2a)]
If $a \leq 0$ and $\|u_0\|_{L^2} >  \|Q^{(a)}\|_{L^2}$, then 
the solution $u(t)$ of \eqref{biNLS} with initial condition $u_0$ blows up in finite time $t^{*}$ in a self-similar blow-up regime of the form 
\begin{equation}\label{conjbu}
u(x) - 	\frac{f(t)}{\lambda(t)^{4/\alpha}}\mathcal{P} \left(\frac{x-x_0(t)}{\lambda(t)}\right)\to \tilde{u}, \quad \tilde{u}\in L^{2}(\mathbb{R}).
\end{equation}

\item[(2b)]
If $a>0$ and the mass $\|u_0\|_{L^2} > (1+ \mu(a)) 
\|Q^{(a)}\|_{L^2}$ for some $\mu(a) >0$ (i.e., larger than the mass 
of a ground state with some non-trivial gap), then the solution $u(t)$ also blows up in a self-similar manner \eqref{conjbu}.  

Regardless of $a$, the blow-up profile in both cases (2a) and (2b) is given by the scale-invariant case $\mathcal{P}= Q^{(0)}$, where 
$Q^{(0)}$ is the ground state solution of \eqref{E:groundstate} with  $a=0$.  
Furthermore, 
\begin{equation}	\label{Lc}
	\lambda(t)=\frac{c}{(t^{*}-t)^{1/3}},
\end{equation}
the blow-up rate in the $\dot{H}^1$ and $L^\infty$ norms is given in \eqref{E:infty-norm} and \eqref{E:L2-norm}, correspondingly; 
and the correction function $f(t)$ in \eqref{conjbu} is such that $\lim_{t\to t^{*}} 
f(t)(t^{*}-t)^{\kappa}=0$ for all $\kappa>0$, and cannot be 
determined numerically. (Recall that in the standard NLS case the correction 
is given by $f(t)=\ln|\ln (t^{*}-t)|$. ) 
\smallskip

\item[(3)] 
If $a>0$ and $(1 - \nu(a)) \|Q^{(a)}\|_{L^2}  \leq \|u_0\|_{L^2} < (1+ \mu(a)) \|Q^{(a)}\|_{L^2}$ for some $\mu(a) >0$ and $\nu(a)\geq 0$, then the solution $u(t)$ from such initial condition approaches asymptotically (possibly in oscillatory manner) a stable, rescaled ground state solution. 
\smallskip

\item[(4)]
If $a>0$ and $ \|u_0\|_{L^2} < (1 - \nu(a)) \|Q^{(a)}\|_{L^2}$ for some $\nu(a) \geq 0$, or if $a \leq 0$ and $\|u_0\| <  \|Q^{(a)}\|_{L^2}$, then the solution disperses away. 
\end{enumerate}

\begin{figure}[!htb]
\includegraphics[width=0.45\hsize,height=0.33\hsize]{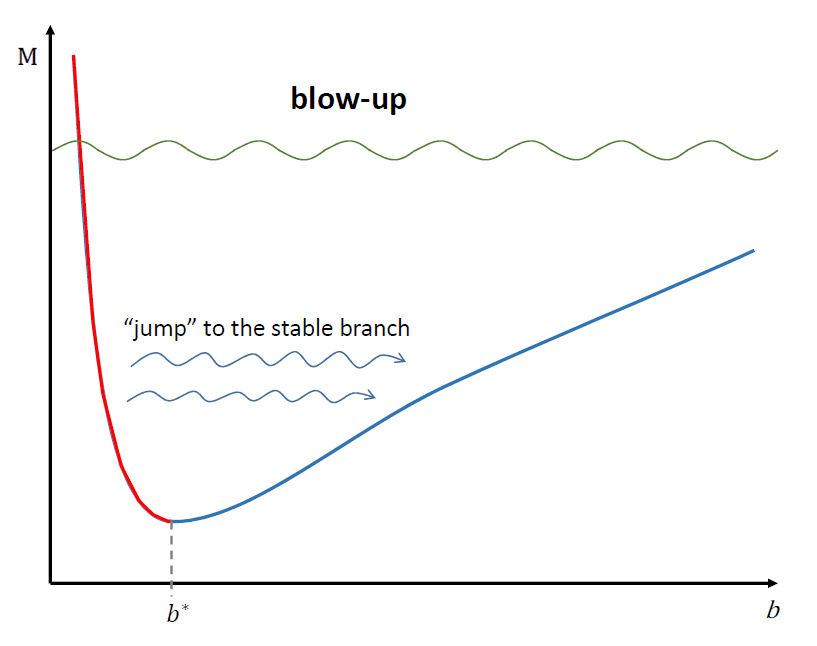}
\includegraphics[width=0.54\hsize,height=0.33\hsize]{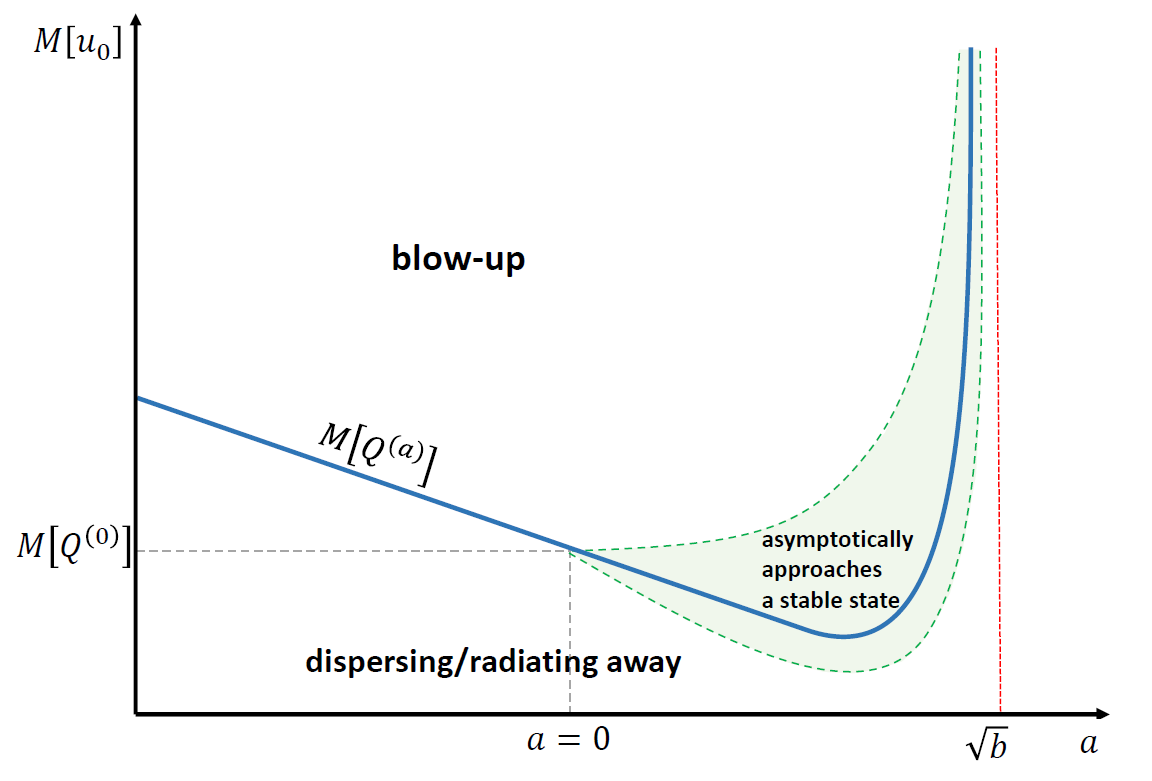}
\caption{\footnotesize Schematic representation of solutions behavior in the critical case: Conjecture II, part 1 (left) and  Conjecture II, parts 2a, 2b, 3, 4 (right). }
\label{F:scheme}
\end{figure}

We show a schematic depiction of Conjecture II in Figure \ref{F:scheme}.

\textbf{Conjecture III:} (1d supercritical case, $\alpha>8$)

In the supercritical case a stable blow-up happens with a self-similar profile as in \eqref{conjbu}, where $\mathcal P$ is a localized smooth solution of the equation \eqref{E:Q1} with a single maximum conjectured to exist and 
\begin{equation}	\label{Lsc}
	\lambda(t)=\frac{f(t)}{(t^{*}-t)^{1/4}},
\end{equation}
with $f(t)$ converging exponentially to a constant (similar to the supercritical blow-up in the standard NLS). 
\smallskip


The structure of this paper is as follows: in Section \ref{S:GS} we 
review the ground state solutions to \eqref{E:groundstate}, which in 
some special cases are known explicitly, and in others are 
constructed numerically. When solving numerically \eqref{E:groundstate} 
for ground states, for some powers of nonlinearities we observe   
two branches in the graphs of the energy vs mass dependence, thus, we investigate 
that bifurcation phenomenon in Section \ref{S:bifurcation}. 
In Section \ref{S:evolv} we describe our numerical approach to track the time evolution of the solution to the 1d bi-harmonic NLS \eqref{biNLS} with a given datum. 
In Section \ref{S:nearQ} we investigate the near soliton dynamics, that is, perturbations of the ground 
states, and finding stable and unstable branches (when these exist) 
of the 1d ground state equation \eqref{E:1dGS}; we then describe various behaviors of solutions  in these branches.    
In Section \ref{S:subcritical-general} we investigate  
solutions to different types of data in the subcritical case of the 1d bi-harmonic NLS, including 
subcritical cases with and without branching of ground states.
In Section \ref{S:critical} we study the critical case ($\alpha=8$) as well as a 
few examples in the supercritical case ($\alpha=10$) and confirm in 1d the 
existence of finite time blow-up solutions, and comment about their rates and profiles. 

{\bf Acknowledgments.} 
The research of this project started during February 2023 visit of 
I.P. and S.R. to the Institut de Math\'ematiques de Bourgogne (IMB), 
they would like to thank the IMB for hospitality. I.P. and S.R. were 
partially supported by the NSF grant DMS-2055130. This material is 
based on research supported by the Swedish Research Council under 
grant no. 2016-06596 while two of the authors (C.K. and S.R.) were in 
residence at Institut Mittag-Leffler (IML) in Djursholm, Sweden, 
during the Fall 2023 semester. Both of them would like to thank for 
the hospitality of the IML and its staff. 
The  work of C.K. and N.K. was partially supported by 
the ANR-17-EURE-0002 EIPHI and the ANR project ISAAC-ANR-23-CE40-0015-0.


\section{Ground States: exact solutions and numerical construction}\label{S:GS}
In this section we discuss solutions to the ground state equation \eqref{E:1dGS}: first, in \S \ref{S:GSexact} 
for some special cases of parameters $a,b$, and $\alpha$, we provide a few explicit solutions of the ground state $Q$; then we write Pokhozhaev identities with several consequences in \S \ref{S:Pokh}; afterwards in \S \ref{S:construction} we construct numerical ground states for any set of parameters, see Fig. \ref{QA}-\ref{QC}. While obtaining numerical ground states (as critical points of energy), we observe bifurcations in the energy vs. mass behavior and investigate that in \S \ref{S:bifurcation}.

\subsection{Exact ground state solutions}
\label{S:GSexact}

In 1d the equation \eqref{biNLS} becomes  
\begin{equation}\label{E:explicit1}
i \partial_{t} u -\partial_{x}^{4}u - 2a \partial_{x}^{2}+|u|^{\alpha}u=0.
\end{equation}
Letting $u(x,t) = e^{i bt} Q(x)$ with $Q$ real (for this work we take $b>0$), one gets that $Q$ satisfies
\begin{equation}\label{E:1dGS}
Q^{(4)} +2a\,Q^{\prime\prime}+b\,Q-Q^{\alpha+1}=0.
\end{equation}

\begin{remark}\label{R:1}
Observe that if we fix one of the parameters, say $a=1$, and consider the equation \eqref{E:1dGS} with this $a$, namely, $Q^{(4)} +2\,Q^{\prime\prime}+b\,Q-Q^{\alpha+1}=0$, then in order to understand the behavior of the solution to \eqref{E:1dGS} with a different value of $a \neq 0$, $a \neq 1$, it suffices to rescale the ground state as $\tilde Q(x) = a^{\frac2{\alpha}} Q (a^{\frac12} x)$. Then $\tilde Q$ solves $\tilde Q^{(4)} +2\,\tilde Q^{\prime\prime}+\frac{b}{a^2}\,\tilde Q- \tilde Q^{\alpha+1}=0$. Thus, the results for $a=1$ could be applied, rescaling them for $\tilde b = \frac{b}{a^2}$. For negative $a<0$, the rescaling is replaced with $|a|$.  
\end{remark}

For the parameters as listed below, the following explicit solutions 
are known, see  \cite{Albert},  \cite{Wazwaz}, \cite{NataliPastor2015}, \cite{PS2020}, 
which we use later to test our numerical solutions:
\begin{itemize}
\item 
\underline{$\alpha=2$} (subcritical): ~for $a < 0$ and $b=\frac{16}{25} a^2$, 

\begin{equation}\label{E:a2}
Q(x) =   \sqrt{\tfrac{6}{5}} \, |a| \,  \sech^2 \Big(\sqrt{\tfrac{|a|}{10}} \, x \Big).
\end{equation}

\item 
\underline{$\alpha=8$} (critical): ~for $a < 0$ and $b=25 (\frac{a}{13})^2$

\begin{equation}\label{E:a8}
Q(x) = \left(\sqrt{105} \,\tfrac{|a|}{13} \right)^{1/4}\, \sech^{1/2} \Big(2 \sqrt{\tfrac{|a|}{13}} \, x \Big).
\end{equation}

\item 
\underline{$\alpha=10$} (supercritical): ~~for $a < 0$ and $b=\big(\tfrac{12\, a}{37}\big)^2 $

\begin{equation}\label{E:a10}
Q(x) = \left( \sqrt{714} \, \tfrac{|a|}{37} \right)^{1/5}\,  \sech^{2/5} \Big(5\,\sqrt{\tfrac{|a|}{74}}x \Big) .
\end{equation}

\end{itemize}

\subsection{Pokhozhaev identities}\label{S:Pokh}

We record the Pokhozhaev identities in the 1d case, which are useful later. For the first one the equation \eqref{E:1dGS} is multiplied by $Q$ and integrated to obtain
\begin{equation}\label{E:P1}
\|\partial_x^2 Q\|^2_{L^2} - 2a \|\partial_x Q\|^2_{L^2} + b\|Q\|^2_{L^2} - \|Q\|_{L^{\alpha+2}}^{\alpha+2} = 0.
\end{equation}
For the second one, the equation \eqref{E:1dGS} is multiplied by $x \partial_xQ$ and integrated to get
\begin{equation}\label{E:P2}
3 \|\partial_x^2 Q\|^2_{L^2} - 2a \|\partial_x Q\|^2_{L^2} - {b} \|Q\|^2_{L^2} + \tfrac{2}{\alpha+2} \|Q\|_{L^{\alpha+2}}^{\alpha+2} = 0.
\end{equation}
Solving for  $\|\partial_x^2 Q\|_{L^2}^2$ and $\|\partial_x Q\|^2_{L^2}$ from \eqref{E:P1} and \eqref{E:P2}, 
we obtain
\begin{equation}\label{E:delta}
\|\partial_x^2 Q\|^2_{L^2} = {b} \, M[Q] - \tfrac{\alpha+4}{2(\alpha+2)} \|Q\|_{L^{\alpha+2}}^{\alpha+2}, 
\end{equation}
\begin{equation}\label{E:grad}
\|\partial_x Q\|^2_{L^2} = \tfrac{b}{a} \, M[Q] - \tfrac{3\alpha+8}{4a(\alpha+2)} \|Q\|_{L^{\alpha+2}}^{\alpha+2}. 
\end{equation}

Recalling the energy and the mass from \eqref{EC} and \eqref{MC}, we obtain the following relation between the mass, energy and the potential term:
\begin{equation}\label{E:EMP}
E[Q] = - \frac{b}2 \, M[Q] + \frac{\alpha}{2(\alpha+2)} \|Q\|_{L^{\alpha+2}}^{\alpha+2},
\end{equation}
or the energy in terms of the mass and first derivative,
\begin{equation}\label{E:EMP2}
E[Q] = \frac{\alpha-8}{3\alpha+8} \, \frac{b}{2} \, M[Q] - \frac{2a \,\alpha}{3\alpha+8} \, \| \partial_x Q\|_{L^{2}}^{2}.
\end{equation}
From the last expression, one can observe that the following holds: 
\begin{itemize}
\item
in the subcritical and critical cases, $\alpha \leq 8$:
if the lower order dispersion coefficient is positive, $a>0$, then the energy of ground states $Q^{(a)}$ is negative, $E[Q^{(a)}] < 0$ (here, the superscript indicates the dependence on $a$);
\item
in the critical case ($\alpha=8$) \eqref{E:EMP2} becomes $E[Q^{(a)}] = -\frac{2a \alpha}{3\alpha +8} \|\partial_x Q\|^2_{L^2}$, and hence,
\begin{itemize}
\item
in the pure quartic case $a=0$:  the energy is zero, $E[Q^{(0)}]=0$,  
\item
when $a<0$: energy is positive, $E[Q^{(a)}]>0$,
\item
when $a>0$: energy is negative, $E[Q^{(a)}]<0$. 
\end{itemize}
\item
in the supercritical case $\alpha>8$: the energy is positive when $a \leq 0$.
\end{itemize}
Furthermore, for the pure quartic case $a=0$, the equations \eqref{E:P1}-\eqref{E:P2} yield:
\begin{itemize}
\item
the energy 
is directly proportional 
to the ground state mass: 
$$
E[Q^{(0)}] = \frac{\alpha-8}{2(3\alpha+8)} b \, M[Q^{(0)}].
$$

\end{itemize}
We confirm some of these observations in our numerical computations in Section \S \ref{S:bifurcation}. 

\subsection{Numerical construction of ground states}\label{S:construction}
In this part, we numerically construct stationary solutions to 
the equation \eqref{E:1dGS} for different parameters. First, we outline the numerical approach and test it for the explicit example in the subcritical case with $\alpha=2$ in \eqref{E:a2}, then we consider examples for various powers $\alpha$ and values of the parameter $a$. 

\subsubsection{Numerical approach} \label{S:GS-num}
We are interested in smooth solutions $Q$ of \eqref{E:1dGS} that are 
critical points of the energy and vanish at infinity. These solutions decay exponentially (see, e.g., \cite{Fibich2002}), thus, Fourier spectral methods are very efficient in this case. Concretely, we apply the same 
approach as in \cite{KS15}, a Fourier spectral approach with a Newton-Krylov iteration. The solution can be chosen to be real (e.g., see \cite{BLS2022}) and having a positive global maximum at the origin. This is enforced 
during the iteration. 

This means we consider equation \eqref{E:groundstate} in the Fourier 
domain. We define the Fourier transform $\hat{u}$ of a function $u\in 
L^{2}(\mathbb{R})$ as 
\begin{align}
	\hat{u}(k) &= \int_{\mathbb{R}}^{}u(x)e^{-ikx}dx, \nonumber\\
	u(x)& = \frac{1}{2\pi}\int_{\mathbb{R}}^{}\hat{u}(k)e^{ikx}dk,
	\label{uhat}
\end{align}
and write \eqref{E:1dGS} in the form 
\begin{equation}
\mathcal{F}(\hat{Q}):=\hat{Q}- \frac{\widehat{Q^{\alpha+1}}}{k^{4}-2ak^{2}+b} = 0.
\label{Qhat}
\end{equation}
Note that the parameter $b$ is always chosen such that $b>a^{2}$ implying 
$k^{4}-2ak^{2}+b 
>0$. For the purpose of this section we take $b=2$. 

To numerically solve the equation \eqref{Qhat}, we approximate the 
Fourier transform by a discrete Fourier transform (DFT), which is 
conveniently computed with a Fast Fourier Transform (FFT). This means 
that the problem is treated as a periodic problem on $L[-\pi,\pi]$, where $L>0$ 
is chosen large enough that the considered functions and their relevant derivatives vanish at the domain boundaries to machine 
precision (we work here with double precision, which is on the order of $10^{-16}$). 
We introduce the standard discretization of the FFT for $N$ Fourier modes, $h = 2\pi L/N$, $x_{j}=-L\pi+hj$, 
$j=1,\ldots,N$. In an abuse of notation, we denote the discrete Fourier transform with the same symbol as the Fourier transform.
Then the equation \eqref{Qhat} becomes a system, call it $\mathcal F$, of $N$ nonlinear equations 
$\mathcal{F}(\hat{Q})=0$  that we solve with a standard Newton iteration
$$
\hat{Q}^{(n+1)} = \hat{Q}^{(n)} -
\left[\mbox{Jac} \big[ \mathcal{F} (\hat{Q})\big]{\big|_{\hat{Q}^{(n)}}} \right]^{-1} \, \mathcal{F}(\hat{Q}^{(n)}),
$$
where $\hat{Q}^{(n)}$ is the $n$th iterate, 
and where $\mbox{Jac} (\mathcal{F})$ denotes the Jacobian of $\mathcal{F}$. The action of the Jacobian on $\mathcal{F}$ is 
computed iteratively with the Krylov subspace method GMRES 
\cite{SS86}.
Note that a further advantage of Fourier spectral methods is that the 
discrete Fourier transform $\hat{u}_k$, $k=(0,\ldots,N)/L$ decreases for analytic functions 
exponentially for  large $|k|$. This allows to control the accuracy of the approximation of a function via the highest terms of 
the DFT, which are of the order of the numerical error, see for instance the discussion in \cite{trefethen}. In this paper we always control the spatial resolution in this way. 

\subsubsection{Test against an exact solution} \label{S:num-exact}
We now test this method for the explicitly known example \eqref{E:a2}, for which we set, for instance, $b=2$. We choose $L=20$, $N=2^{8}$, and the initial iterate $Q^{(0)}(x)=0.5\,e^{-x^{2}}$. The convergence of Newton iterations obviously depends strongly on the choice of the initial iterate, which is here similar in size to the exact solution \eqref{E:a2}, but with a considerably faster decay 
for large $x$ (that is, $e^{-x^2}$ vs. $e^{-a|x|}$). Nevertheless,  the computation converges after 10 iterations (to be precise, it is 
stopped once the residual $\|\mathcal{F}\|_{\infty} \,  < 10^{-10}$). 
\begin{figure}[htb!]
  \includegraphics[width=0.49\textwidth,height=0.35\textwidth]{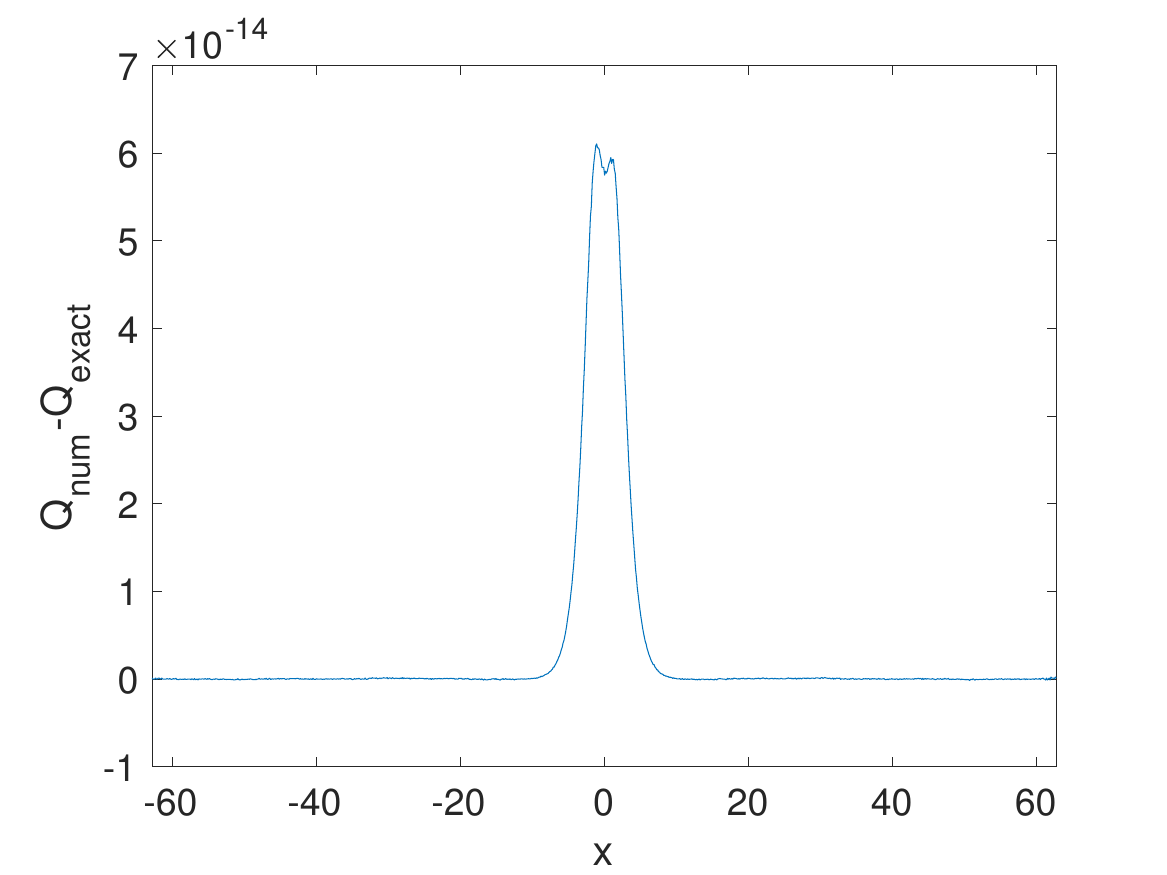}
\caption{\small The difference between the numerically computed $Q$ from \eqref{Qhat} and the exact solution from \eqref{E:a2}.}
 \label{F:num-exact}
\end{figure}
The difference between numerical and exact solution is on the order of 
$10^{-13}$, i.e., roughly on the order of the rounding error, which we illustrate in Fig.~\ref{F:num-exact}. 
We note that the residual $\|\mathcal{F} \|_{\infty}$ for the exact solution 
is on the order of $10^{-15}$.

\subsubsection{Examples}
Below we show several examples of solutions to 
\eqref{E:groundstate} with a fixed value $b=2$, varying the nonlinearity power $\alpha$ and the coefficient of the lower (second order) dispersion $a$. Here, we take $L=10$, $N=2^{10}$ and the initial 
iterate $Q^{(0)}(x) = 1.5e^{-x^{2}}$ unless noted otherwise. \\

\begin{figure}[htb!]
  \includegraphics[width=0.49\textwidth,height=0.36\textwidth]{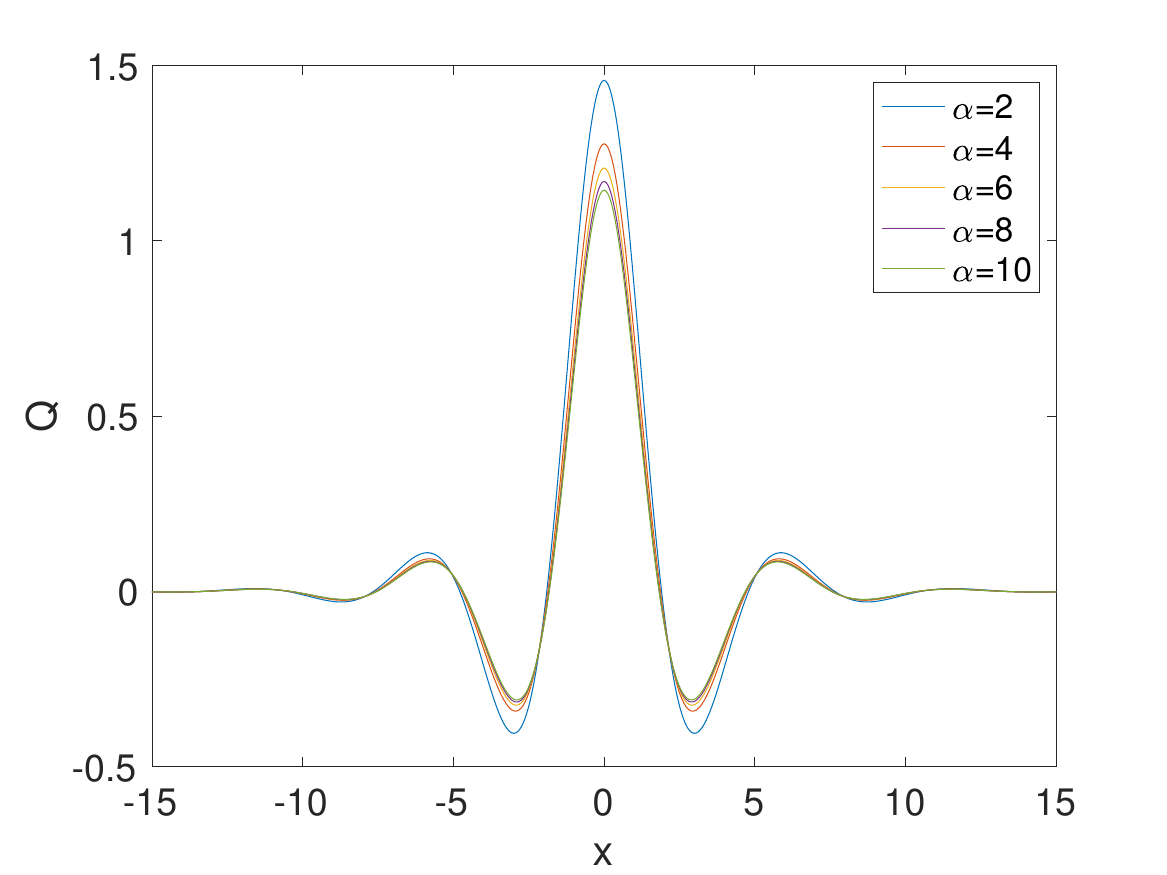}
  \includegraphics[width=0.49\textwidth,height=0.36\textwidth]{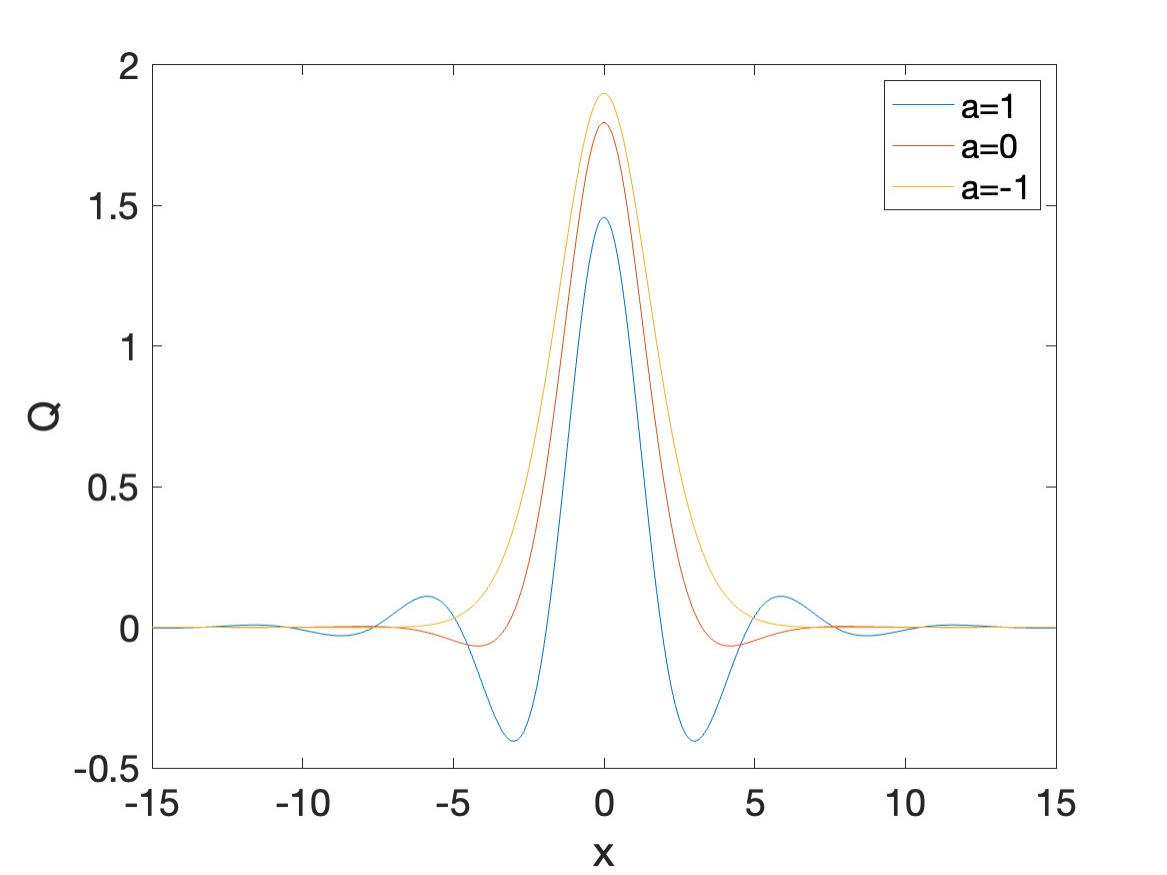}
\caption{\small Profiles of ground state solutions to \eqref{E:1dGS} with $b=2$. Left: $a=1$, $2 \leq \alpha \leq 10$. Right: cubic nonlinearity ($\alpha=2$), coefficient of lower dispersion $a= -1,0,1$.}
 \label{QA}
\end{figure}

On the left of Fig.~\ref{QA} one can see the profiles of the ground states for a fixed $a=1$ with varying nonlinearity power $\alpha$: between $2$ (subcritical case) and $10$ (supercritical case), noting that $\alpha=8$ corresponds to the critical case. All solutions 
are non-monotonic and all have a depression into negative values 
around the central hump and then continuing out with damped 
oscillations; this is due to the higher order dispersion (and mixed 
dispersion), which breaks positivity: in this equation the fourth 
order dispersive term is coupled with the second order term, and unless the second order dispersion is `helping' the higher order dispersion with the very negative coefficient ($a\ll 0$), the profile will have oscillations; similar phenomena are seen in other equations, 
for instance, in the case of the Benjamin equation \cite{TBB2,ABR}. There is no decisive effect of different nonlinearities on the overall shape of the solutions for a fixed $a$, just the overall height is slightly decreasing and the larger values of the nonlinearity $\alpha$ lead to slightly smaller amplitudes. 

On the right of Fig.~\ref{QA}, one can observe changes in the profiles for a (fixed) cubic nonlinearity ($\alpha=2$) with varying $a$, the coefficient of the lower (second order) dispersion. In particular, the more negative $a$ becomes, the less and less oscillations can be seen. We discuss this in further details in the next figure. 

As far as the numerical computations of the profiles, we point out that some relaxation is needed for negative values of $a$ when computing the profile: 
instead of $Q^{(n+1)}$  of the Newton iteration, the value 
$\mu Q^{(n+1)}+(1-\mu)Q^{(n)}$ with $\mu\in(0,1)$ is chosen as the 
new iterate (for instance, with $\mu=0.1$).

\begin{figure}[!htb]
  \includegraphics[width=0.49\textwidth,height=0.36\textwidth]{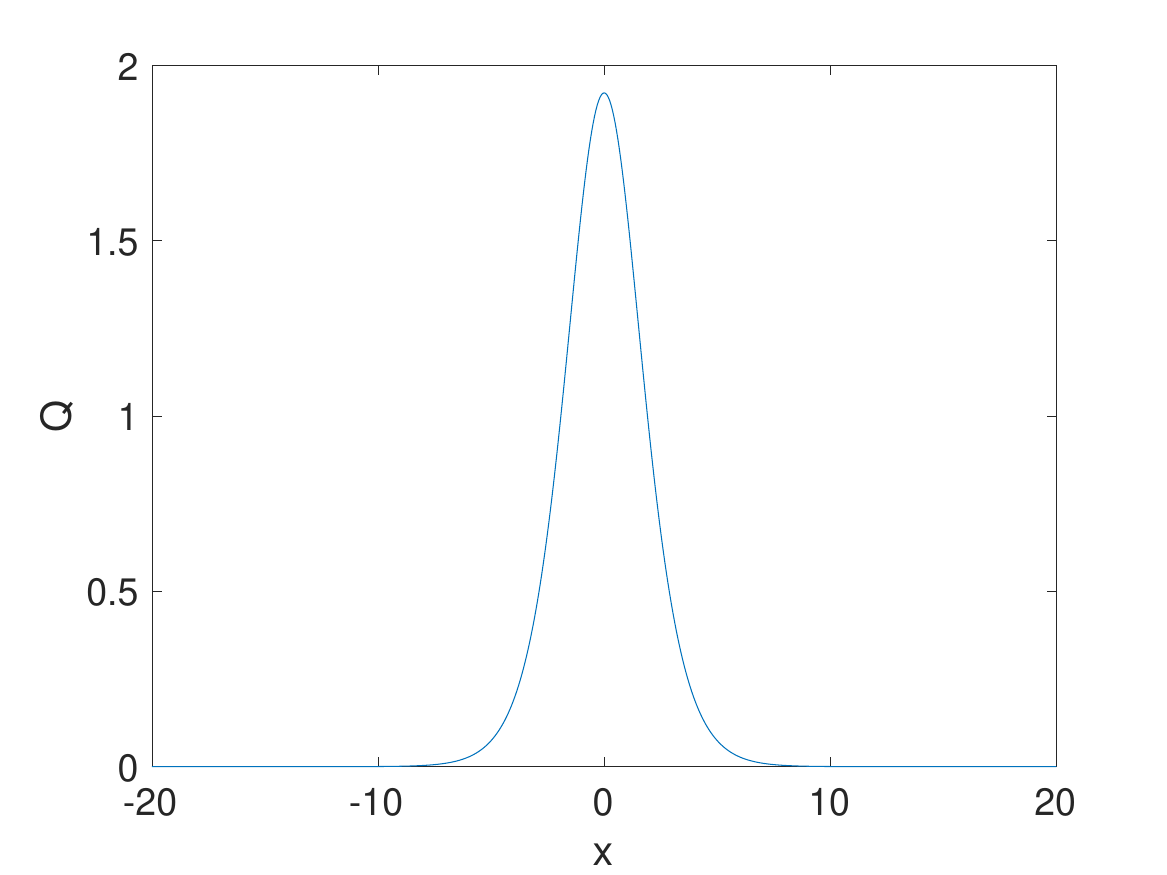}
  \includegraphics[width=0.49\textwidth,height=0.36\textwidth]{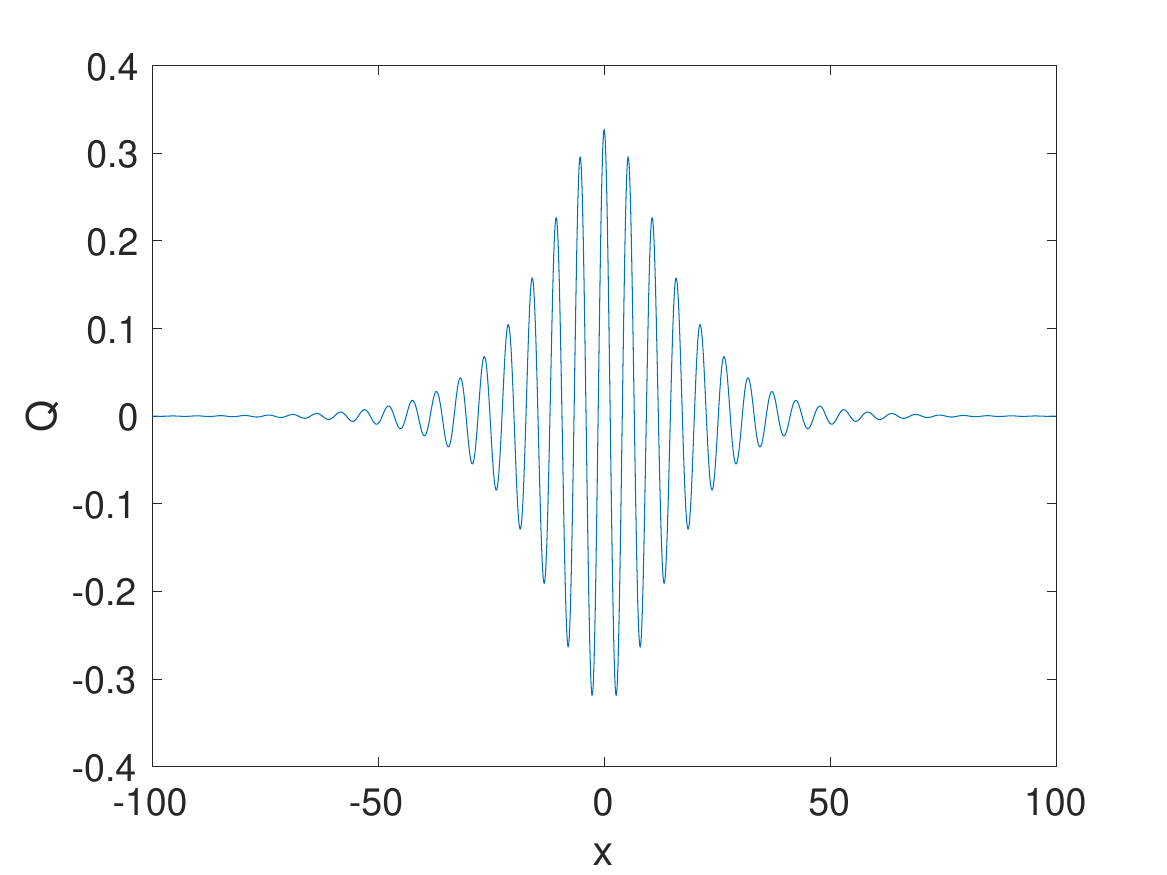}
\caption{\small {Profiles of ground state solutions to \eqref{E:1dGS} with $b=2$ and cubic nonlinearity ($\alpha=2$). Left: $a=-\sqrt{2}$. Right: $a=1.4$. }}
 \label{QC}
\end{figure}
To make further clarification about the oscillatory vs. monotone nature of the profiles, we plot the solution for 
$a=-\sqrt{2}$ on the left of Fig.~\ref{QC} (for this example we chose 
the initial iterate $Q^{(0)}(x)=2e^{-x^{2}}$), noting that it is positive and monotone (as a function of $|x|$). This property of positivity and monotonicity (in $|x|$) will be shared by all profiles with $a \leq -\sqrt 2$, 
which is a confirmation of the results on positivity and radiality of ground states in the regime $a \leq -\sqrt{b}$ in a general case from \cite{BN2015} and \cite{BCSN2018}.  
As we increase $a$ from the value $-\sqrt 2$, we start observing more and more oscillatory behavior of the profiles while decreasing in height as shown in Fig.~\ref{QA}, and as $a\to \sqrt{2}$, solutions become very oscillatory, with height eventually diminishing to zero. 
We show the (almost limiting) case $a=1.4$ on the 
right of Fig.~\ref{QC} (the profiles exist for $a < \sqrt b = \sqrt 2$). For this computational example we chose $L=40$, $N=2^{11}$ 
and the initial iterate $Q^{(0)}(x)=0.5e^{-x^{2}}$ with the relaxation 
discussed above. Note that the height of the solution has decreased to 
0.35 and the number of oscillations has significantly increased 
compared to the cases $a=0, 1$ as in Fig.~\ref{QA}.

\subsection{Bifurcation of ground states}\label{S:bifurcation} 
We investigate the dependence of ground state solutions $Q$ on the parameters $a$ and $b$ more carefully.

\subsubsection{Dependence on $a$.}
\begin{wraptable}{r}{5.2cm}
\centering
{\footnotesize
\begin{tabular}{ |c|l|  }
\hline
$a$ & ~~$M[Q^{(a)}]$ \\ 
\hline
\hline
 1.35 & ~3.41262917  \\ 
 1 & ~2.465972485370718 \\ 
    0 & ~2.986792978326142 \\
    -1 & ~3.604140616082845 \\
    -2 &  ~4.17308102\\
\hline
\end{tabular}
}
\caption{\footnotesize {Mass $M[Q^{(a)}]$ for different $a$ with fixed $b=2$, $\alpha=8$.
}}
\label{T:1}
\end{wraptable}
First, we fix $b=2$ and study how the properties of ground state solutions to \eqref{E:1dGS} change with the parameter $a$, the coefficient of the lower order dispersion in \eqref{biNLS}. For a given value of $a$, we denote the ground states as $Q^{(a)}$ and recall that ground states only exist for $a < \sqrt b \equiv \sqrt 2$. 

We consider the critical case $\alpha=8$ and in Table \ref{T:1} give 
the values of the ground states mass $M[Q^{(a)}]$ for several values 
of $a$. One can observe that the mass is not {\it monotonic} in this dependence.

\begin{figure}[!htb]
\includegraphics[width=0.34\hsize, height=0.3\hsize]{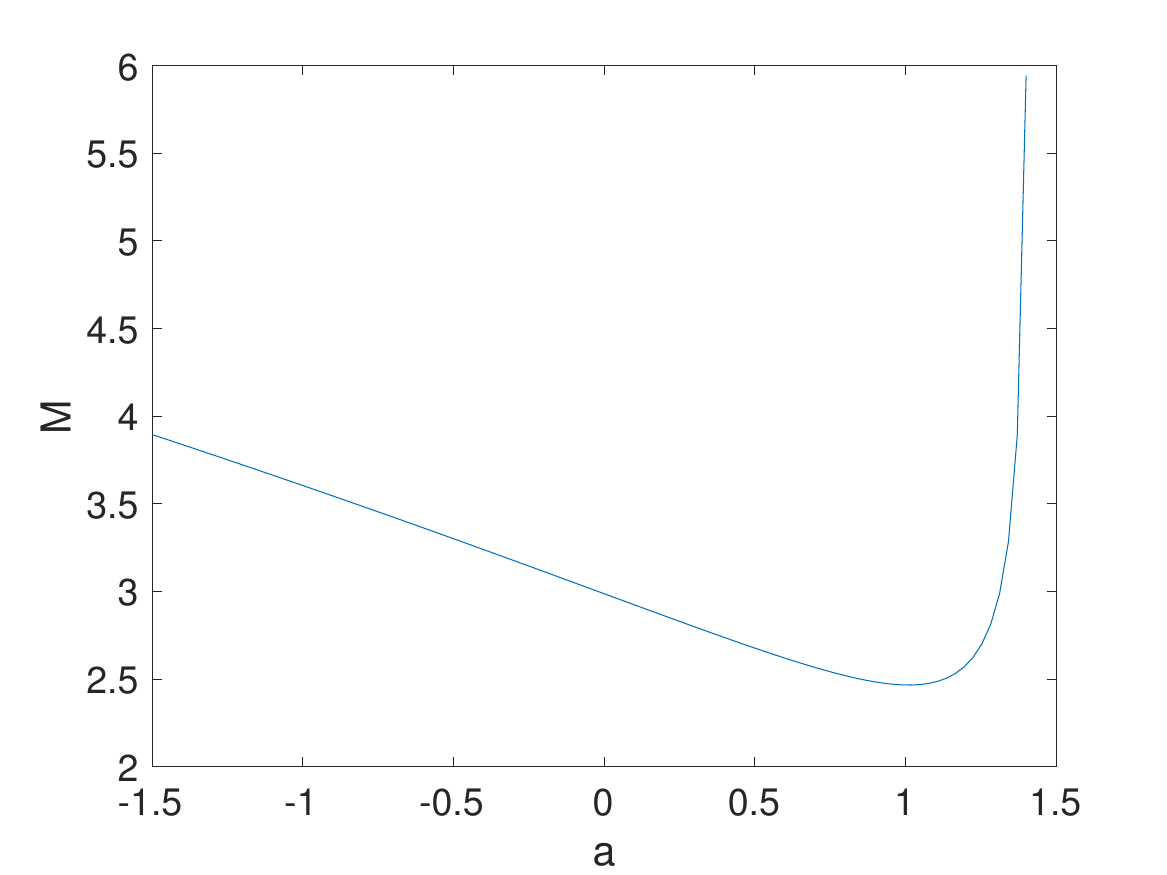}
\includegraphics[width=0.32\hsize, height=0.3\hsize]{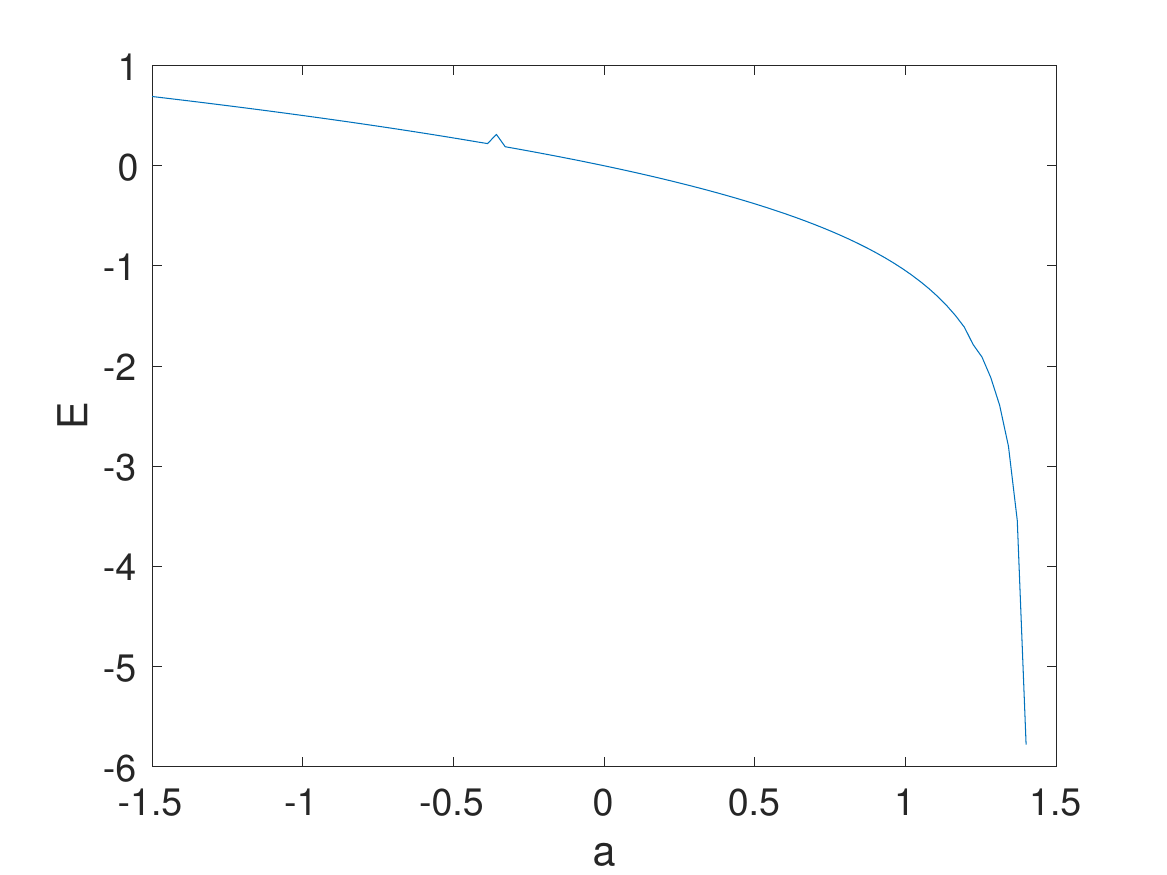}
\includegraphics[width=0.32\hsize, height=0.3\hsize]{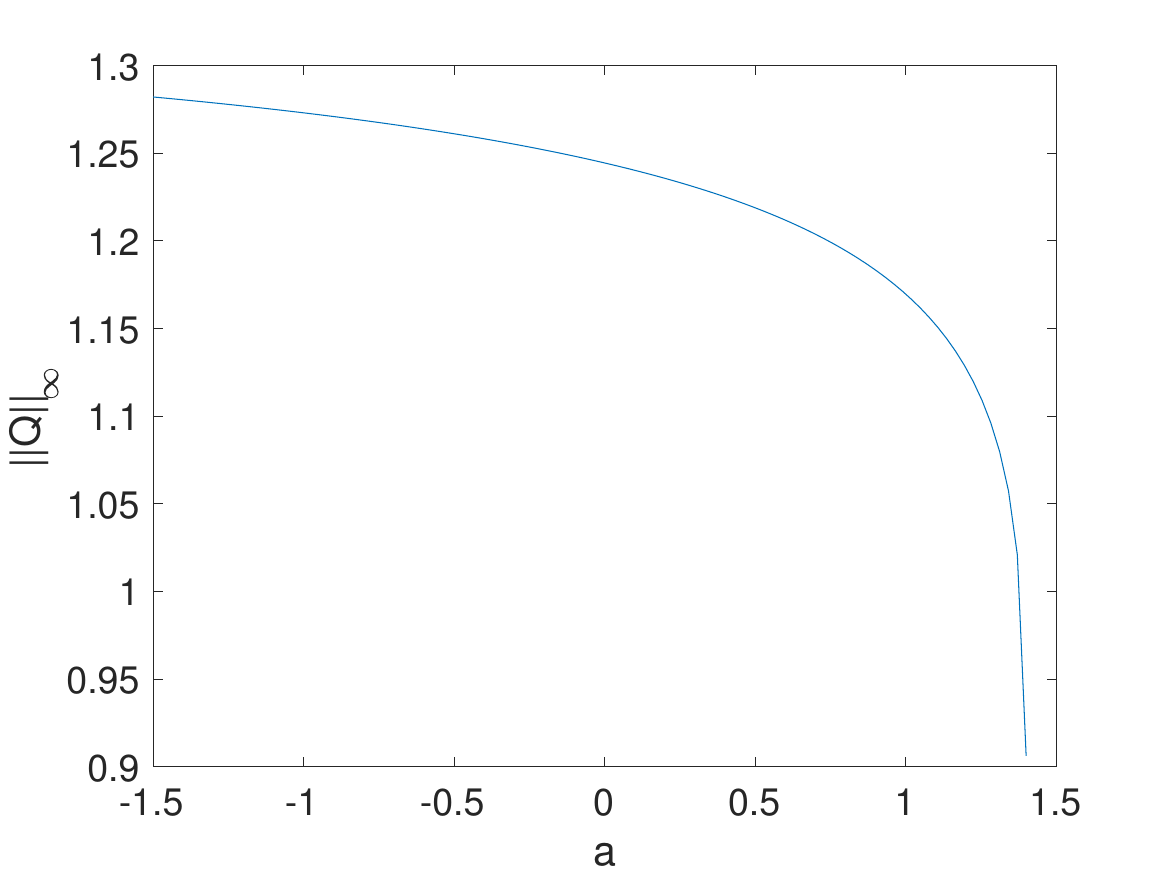}
\caption{\small $\alpha=8$. Dependence of the ground state mass $M[Q^{(a)}]$ on the parameter $a$ for a fixed $b=2$. Mass (left),  energy (middle), $L^{\infty}$ (right).}
\label{F:alpha8}
\end{figure}

To study this further we investigate the dependence of $a$ on several quantities of $Q^{(a)}$, 
such as the mass, energy, and the $L^\infty$ norm, see Fig.~\ref{F:alpha8}. 

One notices that the mass has a minimum around $a \approx 1$, while 
the energy and the sup norm are decreasing as $a \to \sqrt 2$. The plot of the 
$L^\infty$ norm shows that the ground states decrease in their 
height as $a$ increases, being consistent with Fig.~\ref{QC} and the right plot of Fig. 
~\ref{QA}, and since the mass is increasing, they gain more and more oscillatory behavior 
as $a \to \sqrt 2$ as shown in Fig.~\ref{QC}. Note from the middle plot of Fig.~\ref{F:alpha8} that the energy is zero when $a=0$ (the scaling invariant case), $E[Q^{(0)}] = 0$, it is positive for $a<0$ and negative for $a>0$, as it was proved at the end of Section \S \ref{S:Pokh}. This will play an important role in the investigation of blow-up in Section \S \ref{S:critical}.

\begin{figure}[!htb]
\includegraphics[width=0.34\hsize, height=0.3\hsize]{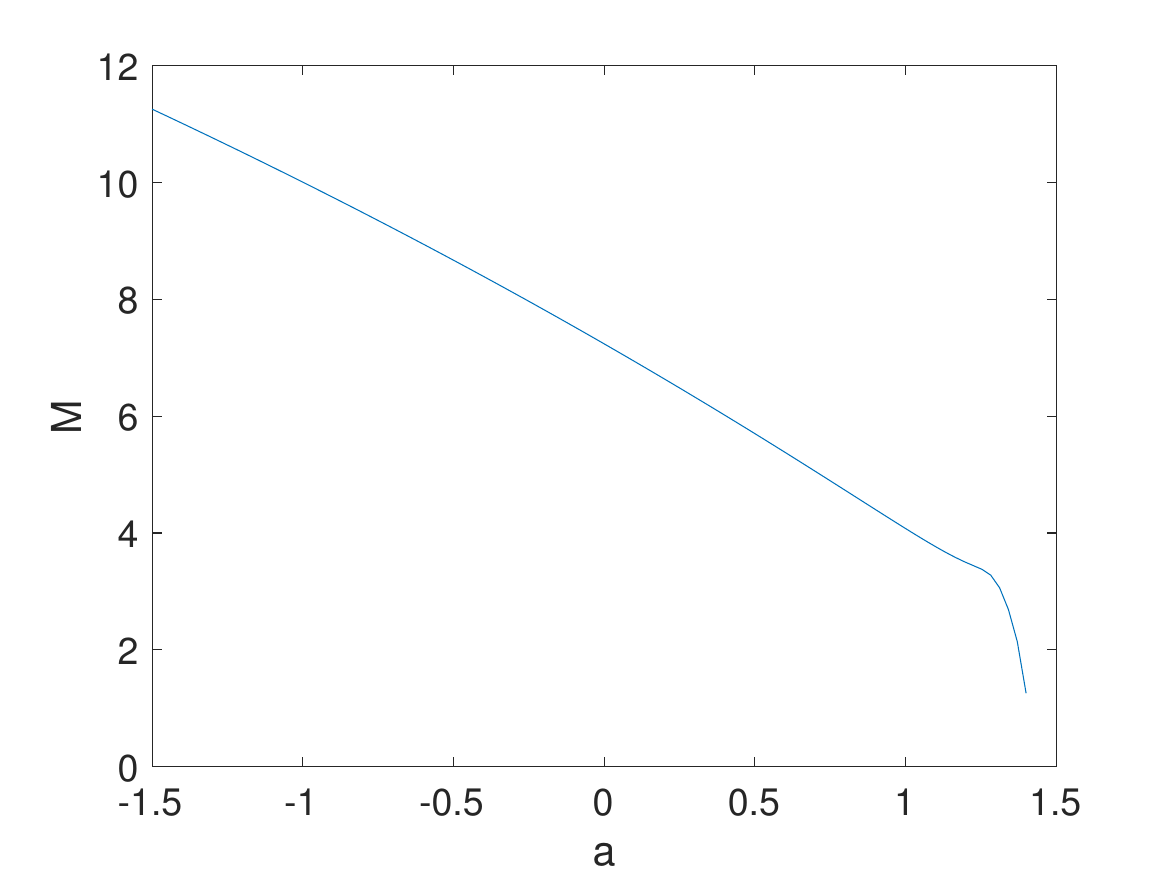}
\includegraphics[width=0.32\hsize, height=0.3\hsize]{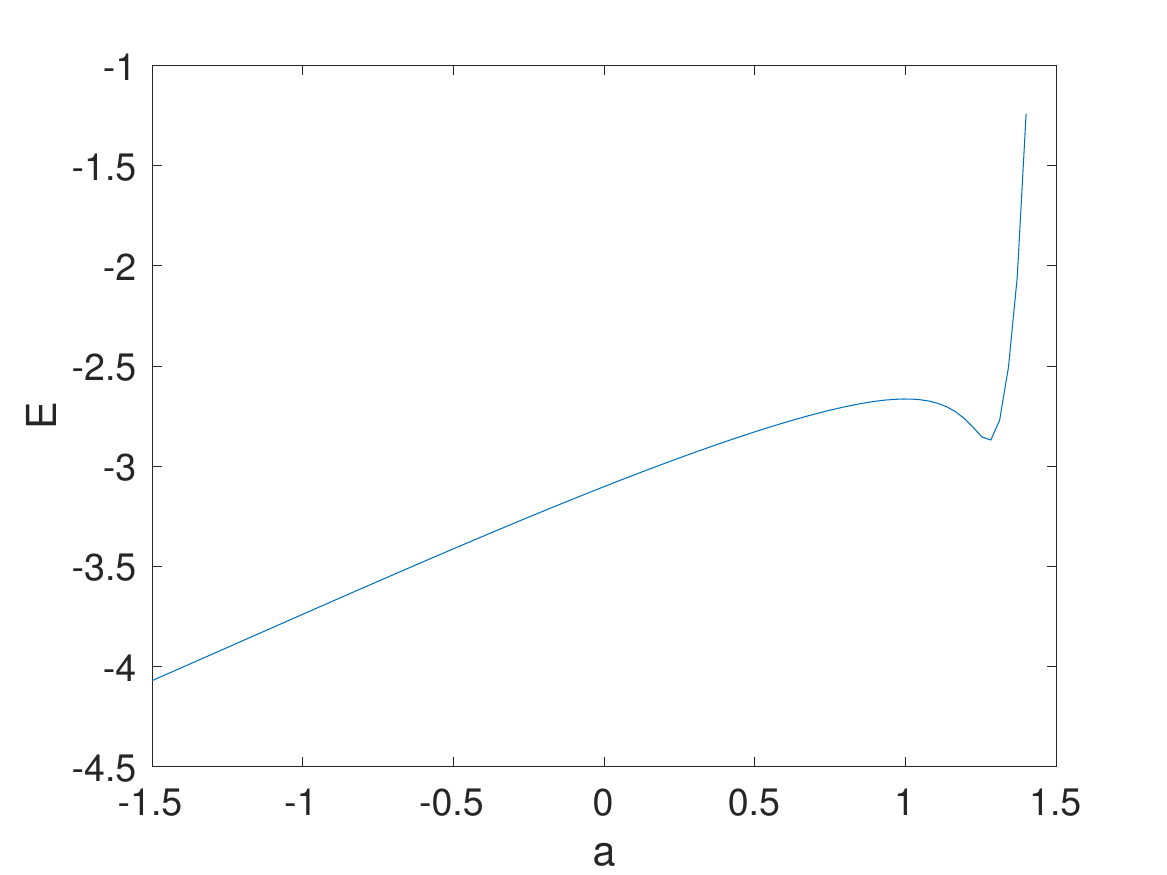}
\includegraphics[width=0.32\hsize, height=0.3\hsize]{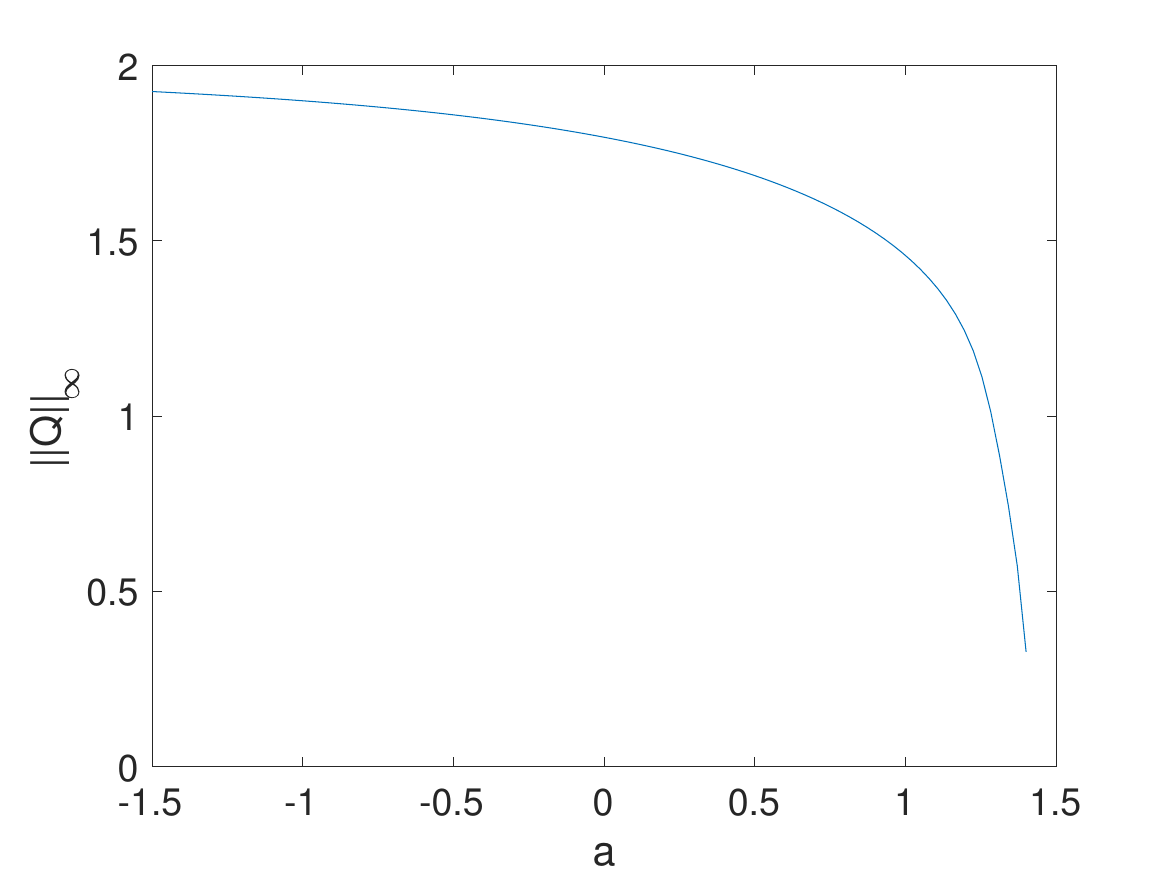}
\caption{\small $\alpha=2$. Dependence of the ground state mass $M[Q^{(a)}]$ on the parameter $a$ for a fixed $b=2$. Mass (left),  energy (middle), $L^{\infty}$ (right).}
\label{F:alpha2}
\end{figure}

We investigate the subcritical case $\alpha = 2$ in Fig. ~\ref{F:alpha2} and observe that, opposite to the critical case, the mass is decreasing as $a$ increases to $\sqrt 2$ value (left plot). The energy is increasing, having a small dip around $a\approx 1.3$ (in the middle plot), and the sup norm decreases as before (right plot), indicating that the oscillatory envelope of the ground state is decreasing in height, and since the mass is decreasing, the decay of the envelope of oscillations is much faster than in the critical case.  

\subsubsection{Dependence on $b$.}\label{S:branching}
We next fix $a=1$, the lower order dispersion coefficient (note that in this case both dispersions work against each other), and track 
the  dependence of the ground state $Q^{(1)}$ quantities on $b$. 
We investigate cases from subcritical to supercritical: $\alpha = 2,4,6,8,10$, to show a new phenomenon about the ground states.

\begin{figure}[!htb]
\begin{subfigure}{.32\textwidth}
\includegraphics[width=1\linewidth,height=0.85\linewidth]{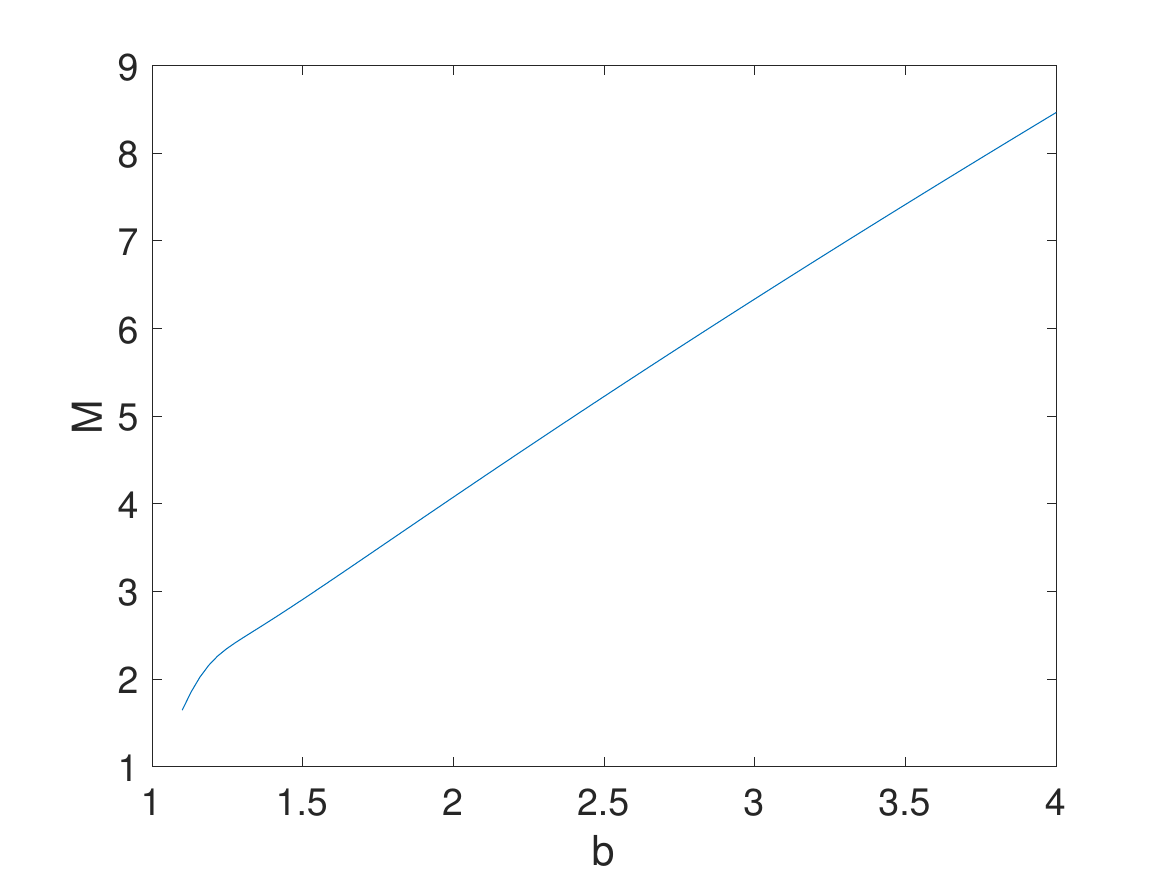}
\subcaption[]{{\footnotesize {$\alpha=2$, $M[Q^{(1)}] = M(b)$.}}}
\end{subfigure}
\begin{subfigure}{.32\textwidth}
\includegraphics[width=1\linewidth,height=0.85\linewidth]{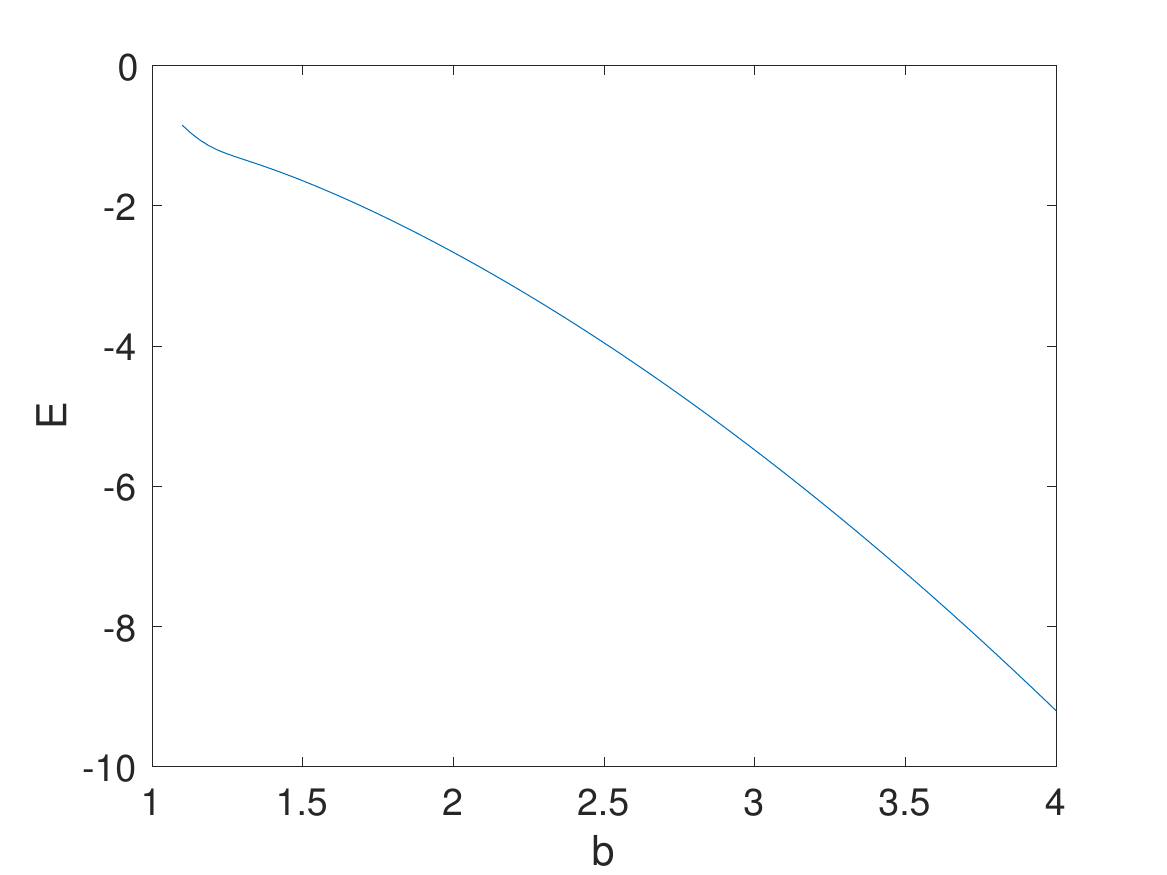}
\subcaption[]{{\footnotesize $E[Q^{(1)}]$ as function of $b$.}}
\end{subfigure} 
\begin{subfigure}{.32\textwidth}
  \includegraphics[width=1\linewidth,height=0.85\linewidth]{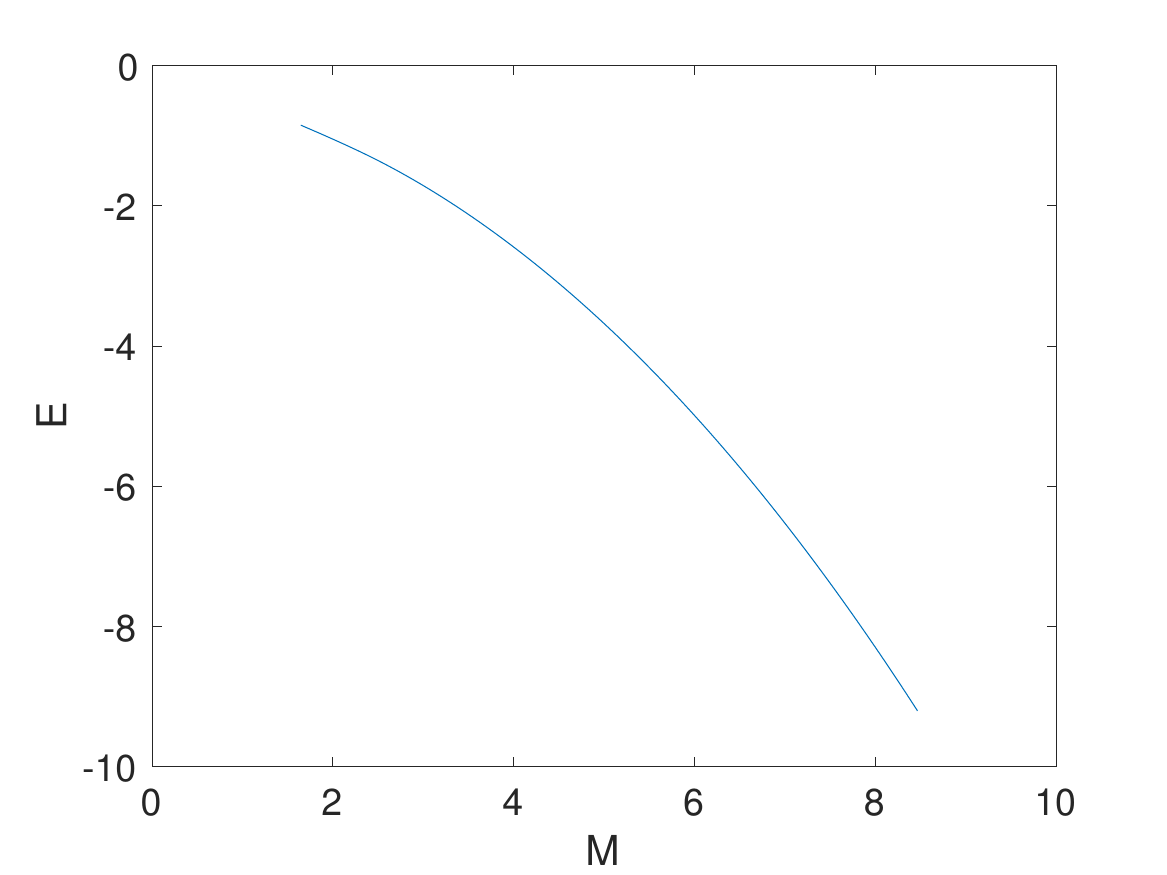}
\subcaption[]{{\footnotesize $E = E(M)$.}}
\end{subfigure}\\
\begin{subfigure}{.32\textwidth}
\includegraphics[width=1\linewidth,height=0.85\linewidth]{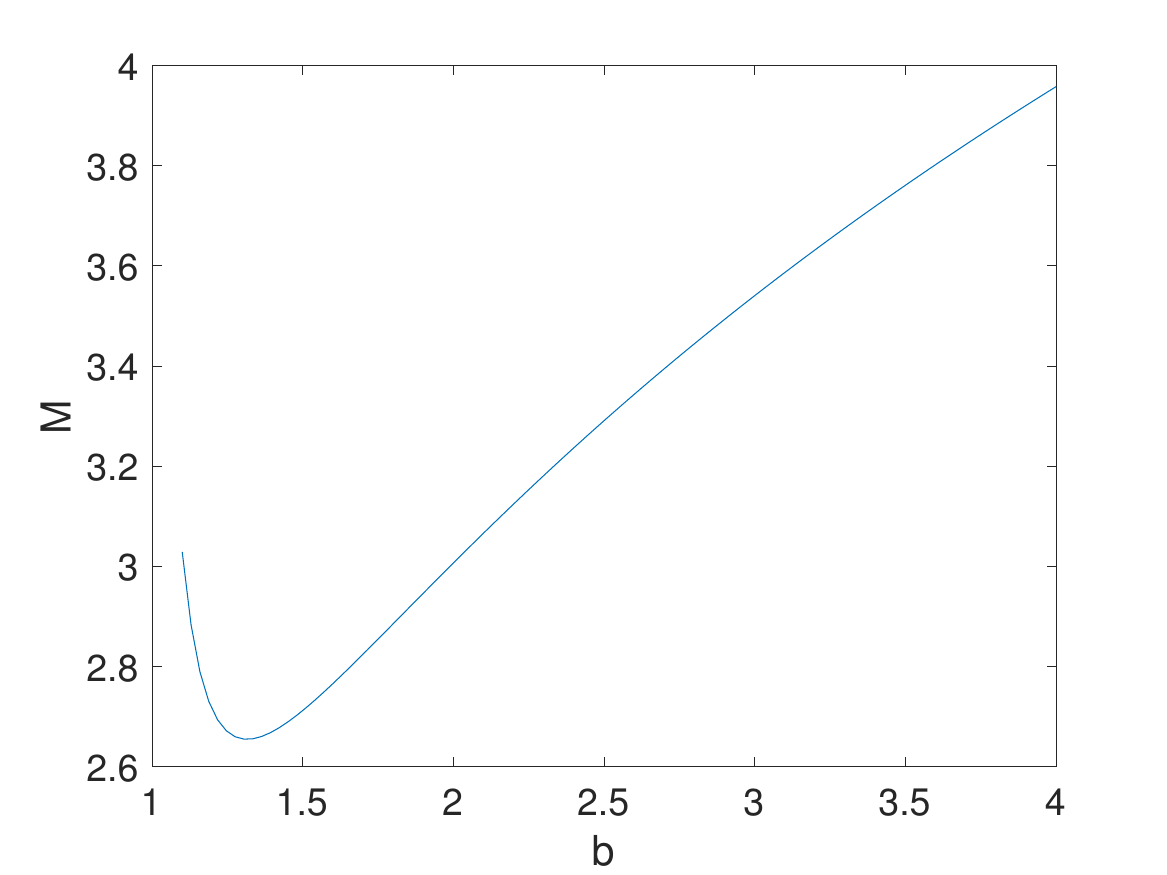}
\subcaption[]{{\footnotesize $\alpha=4$, $M[Q^{(1)}]=M(b)$.}}
\end{subfigure}
\begin{subfigure}{.32\textwidth}
  \includegraphics[width=1\linewidth,height=0.85\linewidth]{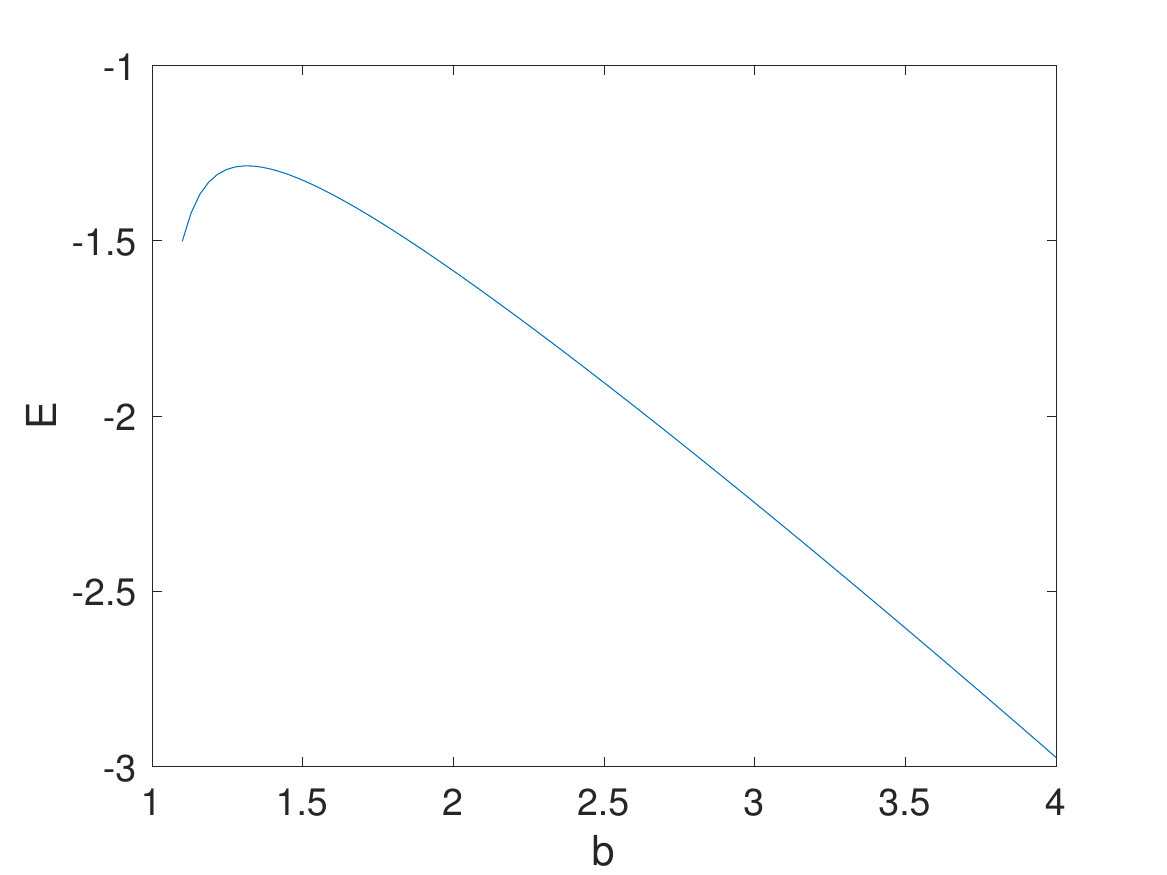}
\subcaption[]{{\footnotesize $E[Q^{(1)}]$ as function of $b$.}}
\end{subfigure}
\begin{subfigure}{.32\textwidth}
\includegraphics[width=1\linewidth,height=0.85\linewidth]{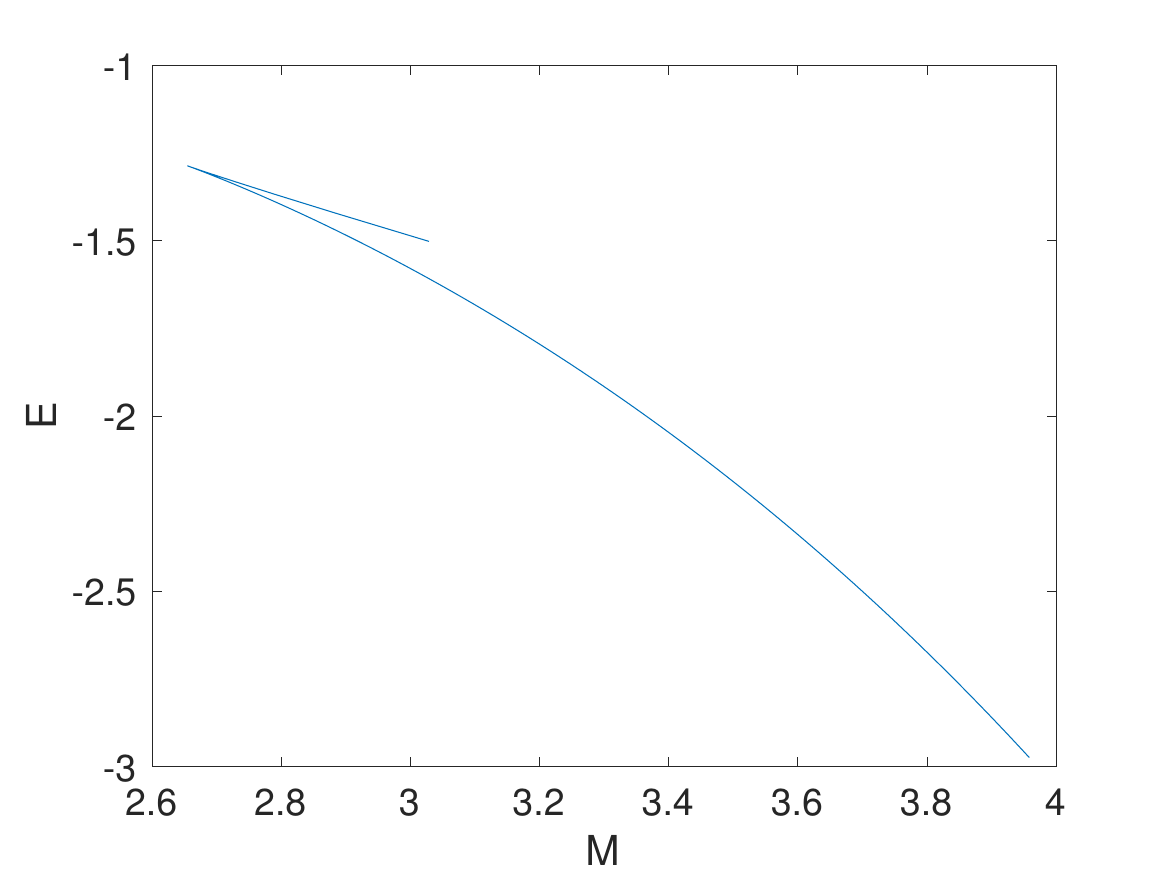}
\subcaption[]{{\footnotesize $E=E(M)$.}}
\end{subfigure}\\
\begin{subfigure}{.32\textwidth}
\includegraphics[width=1\linewidth,height=0.85\linewidth]{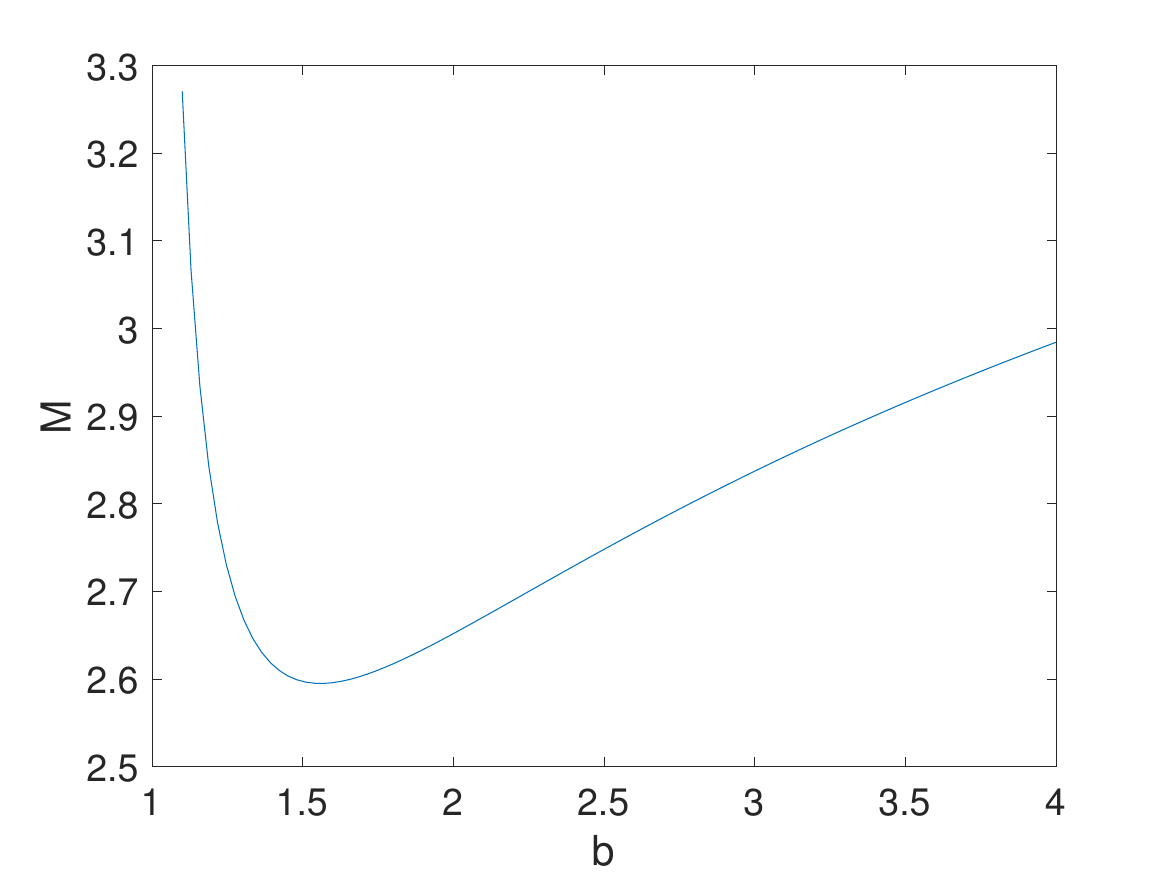}
\subcaption[]{{\footnotesize $\alpha=6$, $M[Q^{(1)}]=M(b)$.}}
\end{subfigure}
\begin{subfigure}{.32\textwidth}
  \includegraphics[width=1\linewidth,height=0.85\linewidth]{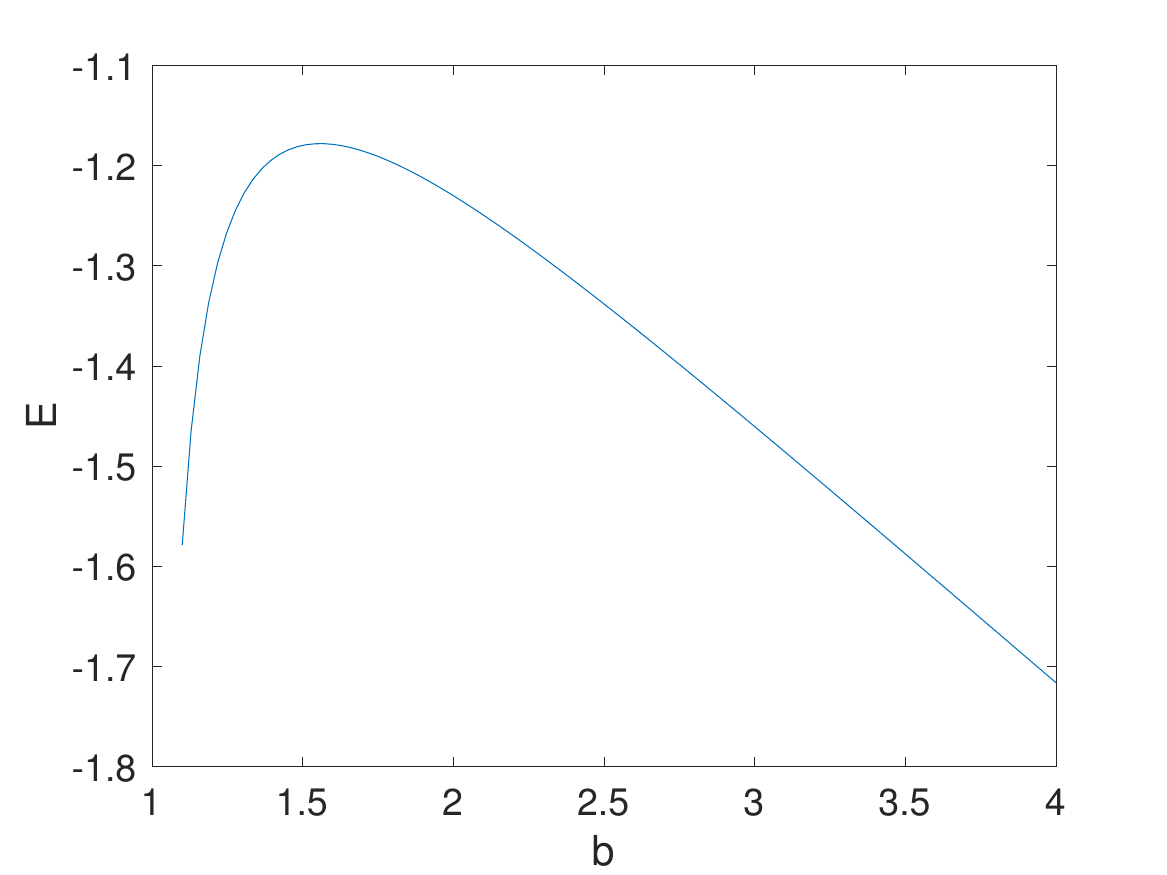}
\subcaption[]{{\footnotesize $E[Q^{(1)}]$ as function of $b$.}}
\end{subfigure}
\begin{subfigure}{.32\textwidth}
\includegraphics[width=1\linewidth,height=0.85\linewidth]{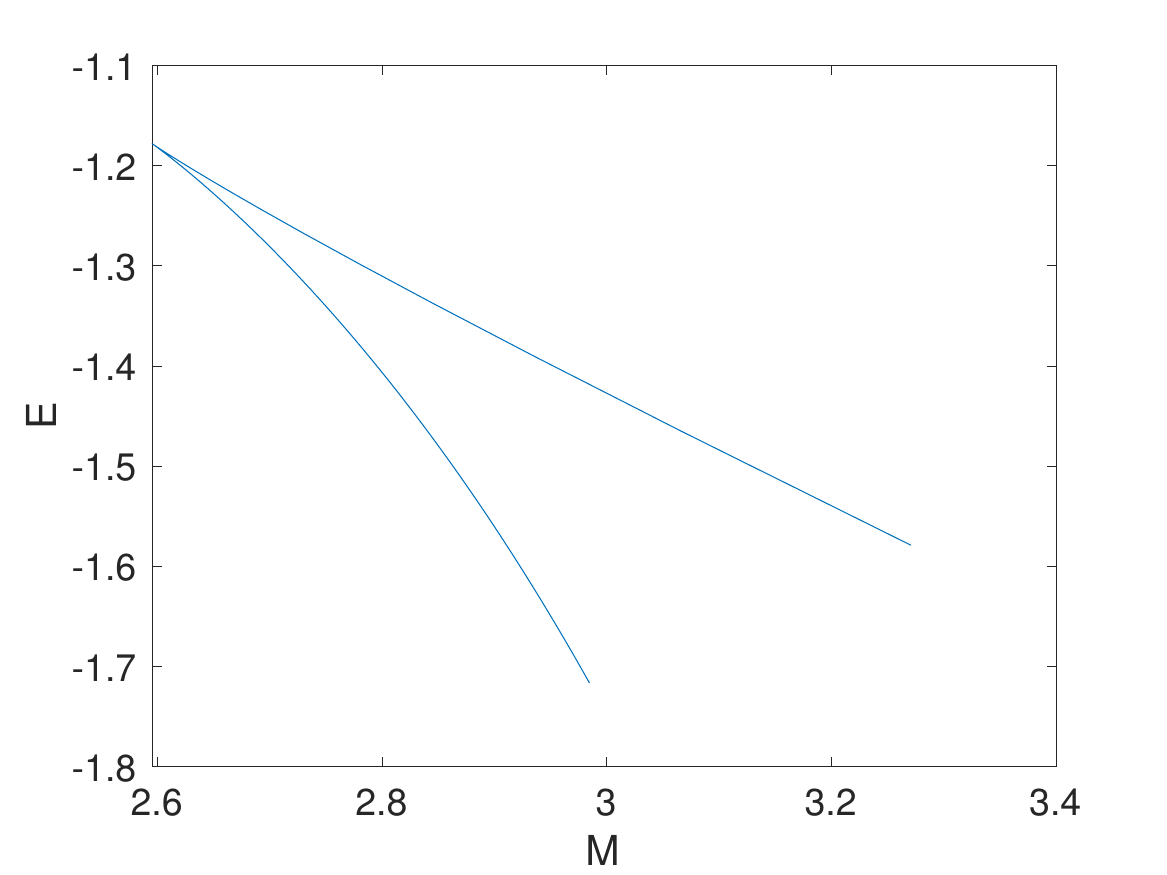}
\subcaption[]{{\footnotesize $E=E(M)$.}}
\end{subfigure}
\caption{\footnotesize {Dependence in sub-critical cases $\alpha=2$ (top), $\alpha=4$ (middle), $\alpha=6$ (bottom), of 
$M[Q^{(1)}]$ and $E[Q^{(1)}]$ on the parameter $b$ for a fixed $a=1$ (left and middle columns). Dependence of energy as a function of mass, $E = E(M)$ (right column).}}
\label{F:ME-subcritical}
\end{figure}

\begin{figure}[htb]
\begin{subfigure}{.32\textwidth}
\includegraphics[width=1\linewidth,height=0.85\linewidth]{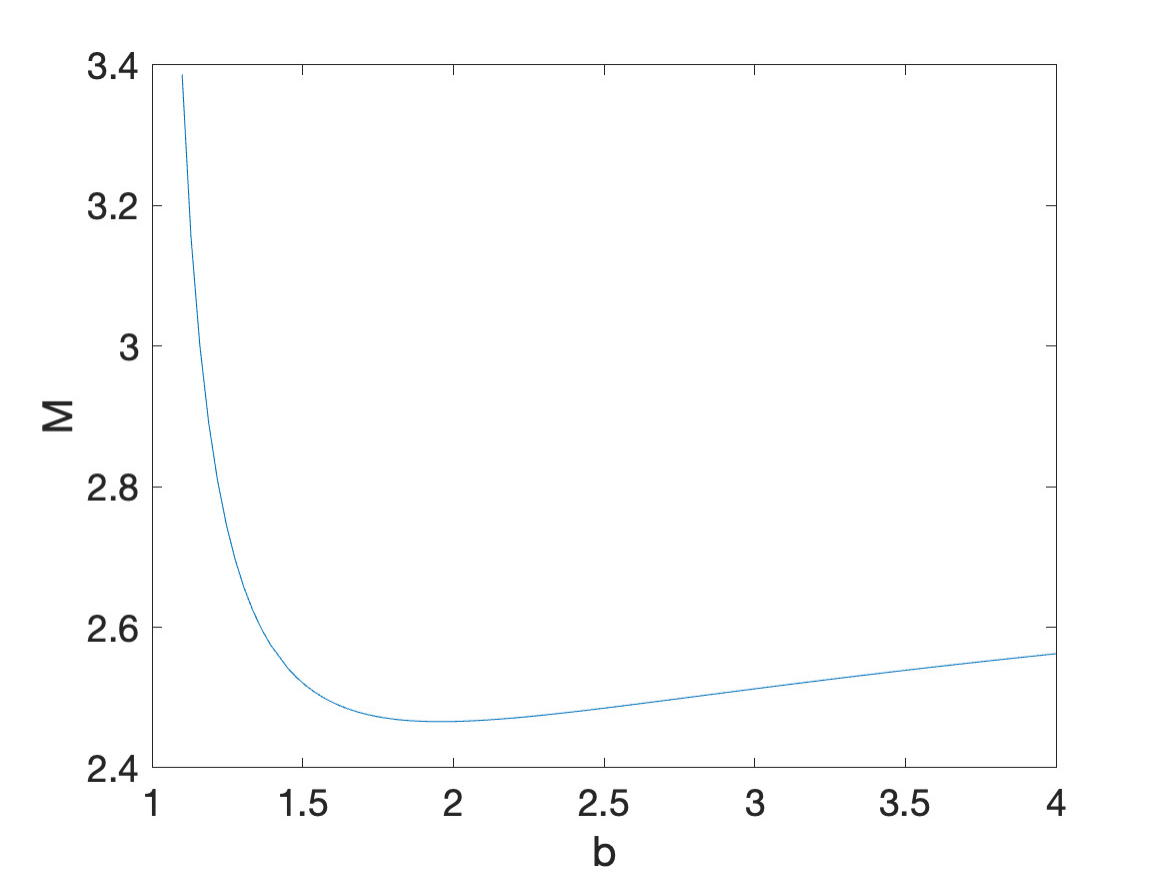} 
\subcaption[]{{\footnotesize $\alpha=8$, $M[Q^{(1)}]=M(b)$.}}
\end{subfigure}
\begin{subfigure}{.32\textwidth}
  \includegraphics[width=1\linewidth,height=0.85\linewidth]{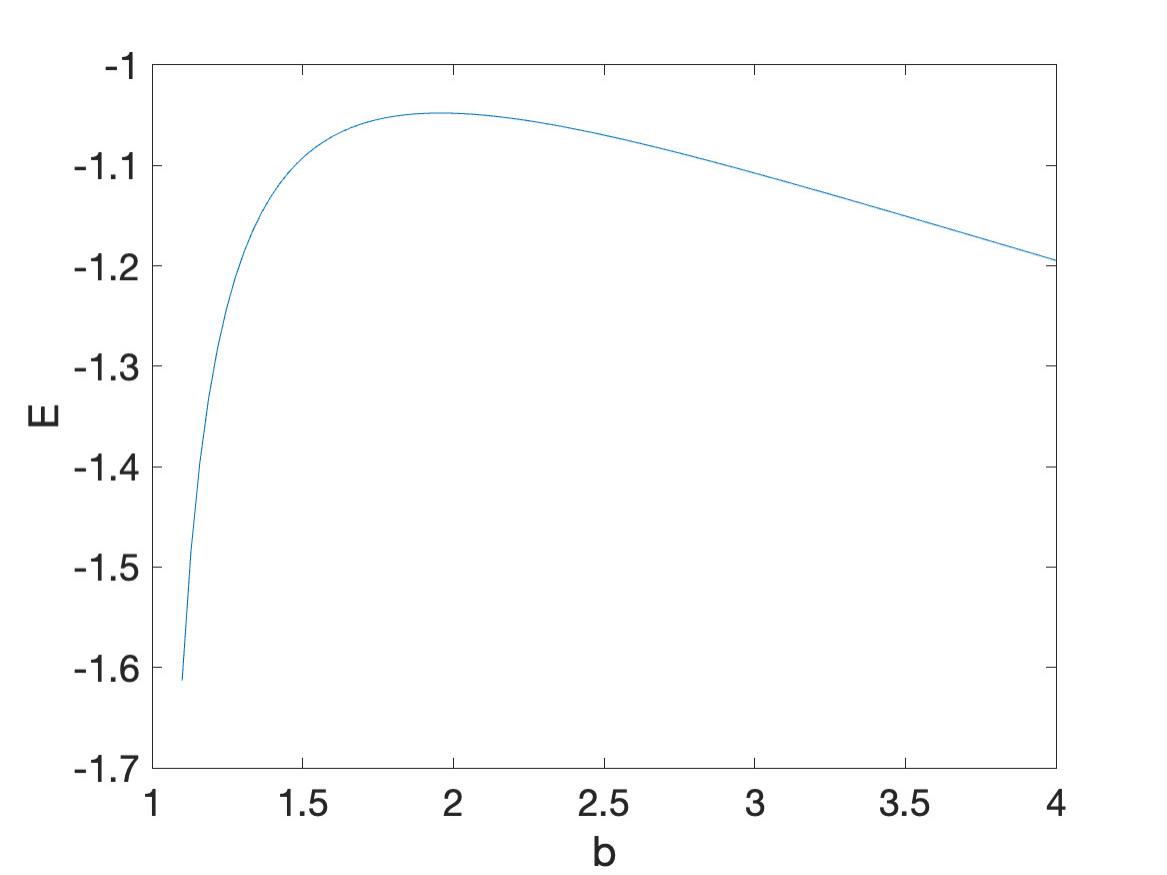}
\subcaption[]{{\footnotesize $E[Q^{(1)}]$ as function of $b$.}}
\end{subfigure}
\begin{subfigure}{.32\textwidth}
\includegraphics[width=1\linewidth,height=0.85\linewidth]{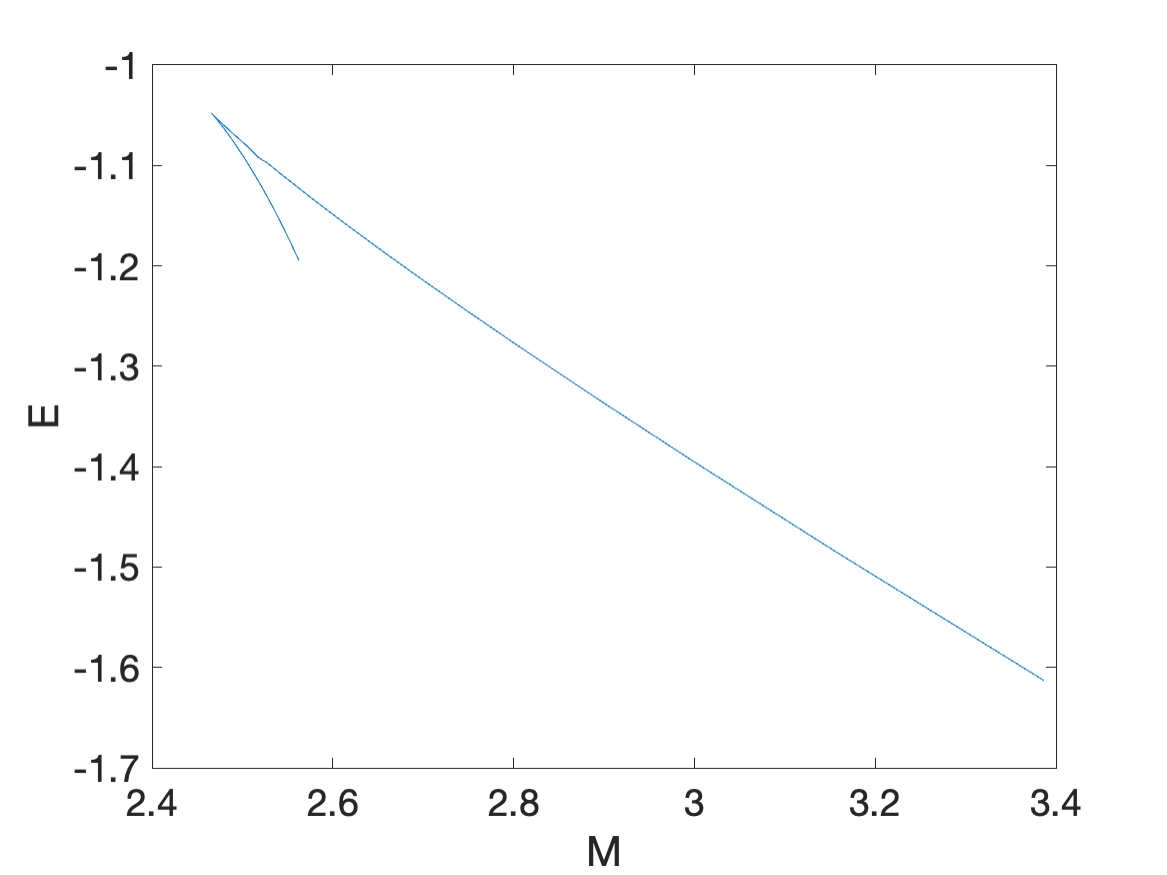}
\subcaption[]{{\footnotesize $E=E(M)$.}}
\end{subfigure}\\
\begin{subfigure}{.32\textwidth}
\includegraphics[width=1\linewidth,height=0.85\linewidth]{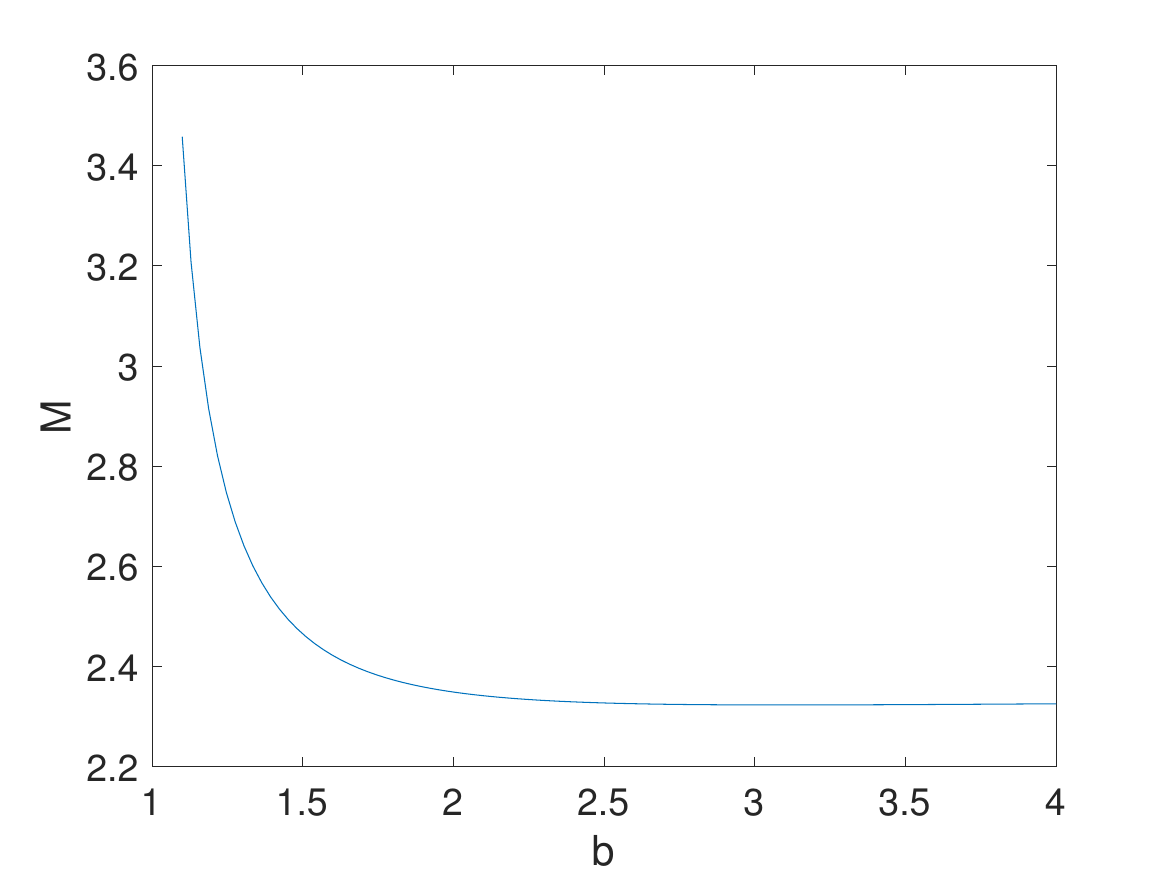}
\subcaption[]{{\footnotesize $\alpha=10$, $M[Q^{(1)}]=M(b)$.}}
\end{subfigure}
\begin{subfigure}{.32\textwidth}
  \includegraphics[width=1\linewidth,height=0.85\linewidth]{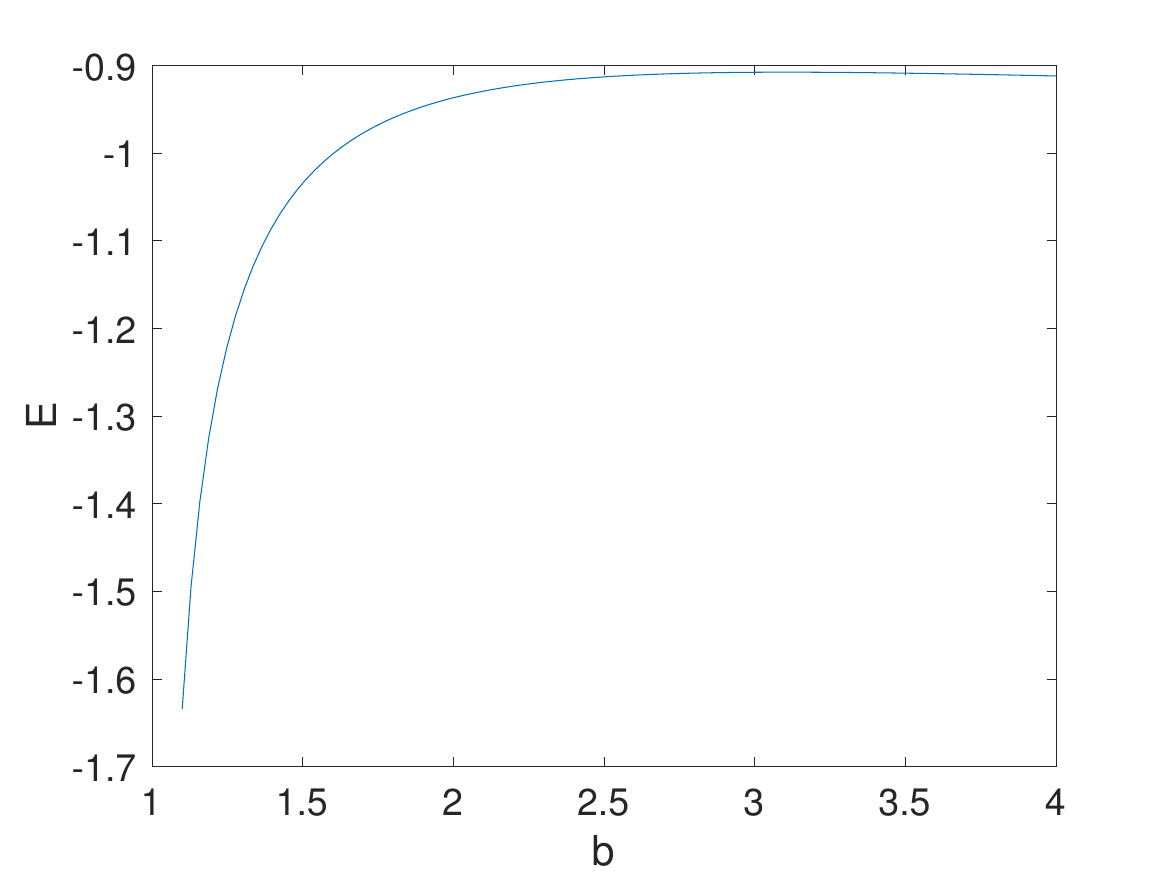}
\subcaption[]{{\footnotesize $E[Q^{(1)}]$ as function of $b$.}}
\end{subfigure}
\begin{subfigure}{.32\textwidth}
\includegraphics[width=1\linewidth,height=0.85\linewidth]{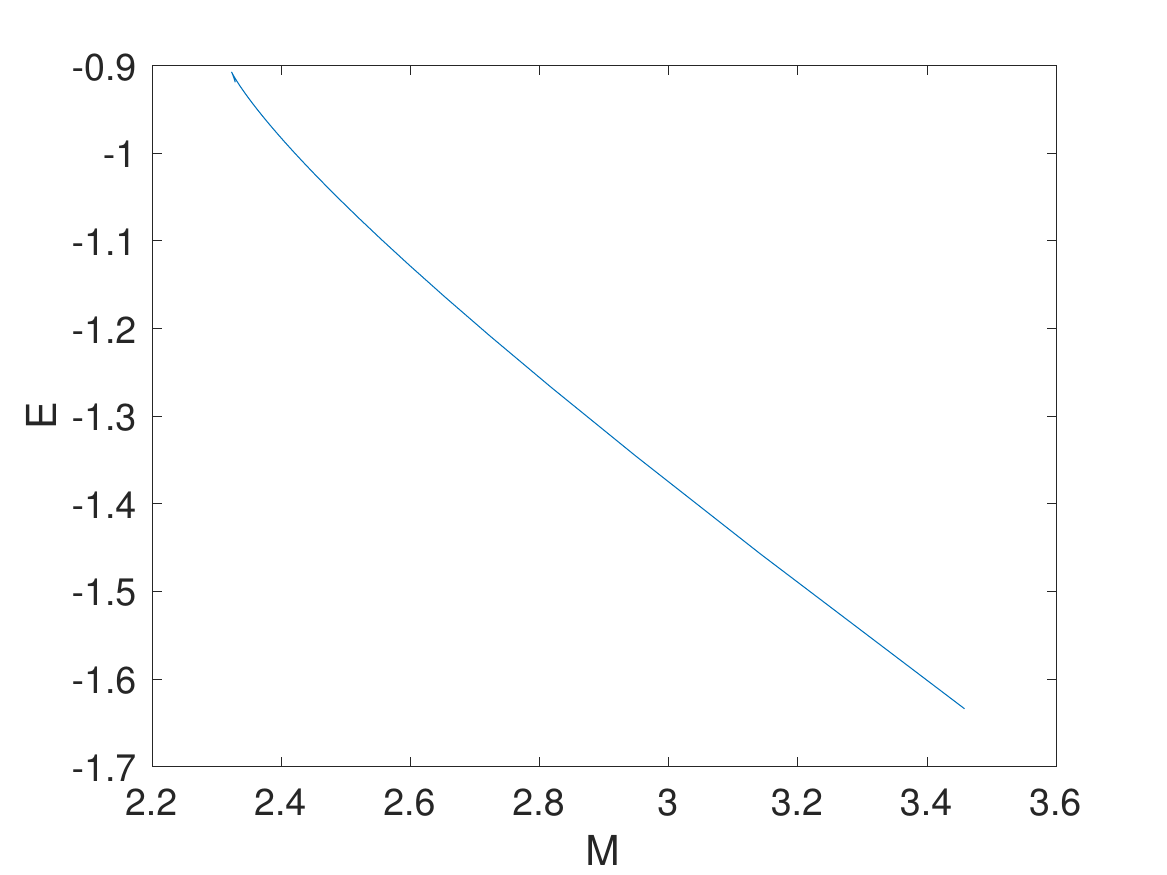}
\subcaption[]{{\footnotesize $E=E(M)$.}}
\end{subfigure}
\caption{\footnotesize {Dependence in the critical $\alpha=8$ (top) and super-critical $\alpha=10$ (bottom) cases, of 
$M[Q^{(1)}]$ and $E[Q^{(1)}]$ on the parameter $b$ for a fixed $a=1$ (left and middle columns). Dependence of energy as a function of mass, $E = E(M)$ (right column).}}
\label{F:ME-crit+supercrit}
\end{figure}

We plot cases $\alpha = 2,4,6$ in Fig. \ref{F:ME-subcritical} and $\alpha=8,10$ in Fig. \ref{F:ME-crit+supercrit}. 
\medskip

First, note that while in the sub-critical case $\alpha=2$ all graphs look to be monotone, in the cases $\alpha=4,6$, there is a different behavior: 
initially, the mass is decreasing as $b$ increases 
and then it starts increasing, while the energy does exactly the 
opposite; the reverse behavior is consistent with the dependence 
shown in \eqref{E:EMP}. This change in monotonicity, when plotted as 
a function, where the energy dependence on the mass (in the computed 
range of $b$ between $1$ and $4$), i.e., $E = E(M)$, shows the appearance of \textit{\textbf {two branches}} in the energy, see the right plots in Fig. \ref{F:ME-subcritical}.  
This means that for the same mass there are 
two solutions for the ground state. This, in some sense, resembles 
the behavior of the NLS equation with a combined nonlinearity or a 
double well potential, for example, see \cite{CKS}. We, therefore, 
call the ground state from the upper branch an \textit{\textbf 
{unstable branch}}, and the lower one  - a \textit{\textbf {stable 
branch}}. As we show later, when the ground state from the upper 
branch is perturbed 
such that it has a larger mass 
than the unperturbed ground state, it will 
jump to the lower branch (with the lower energy) and will try to 
approach asymptotically that ground state. On the other hand if the 
perturbation leads to a situation with less mass than the unperturbed 
state, the initial data are simply dispersed.
This is investigated in more details in Section \ref{S:jumping}.

In Fig.~\ref{F:ME-crit+supercrit}, we plotted the dependencies of the same conserved quantities (mass and energy) on the parameter $b$ in the critical and supercritical cases. 
Observe that in the critical case the behavior of the quantities is similar to the sub-critical cases with $\alpha=4$ and $6$, producing a bifurcation in the energy vs. mass plot, $E = E(M)$, see top row in Fig.~\ref{F:ME-crit+supercrit}. 

In the supercritical case, $\alpha=10$, the dependence of mass and energy becomes monotone (as in the case $\alpha=2$), 
and thus, no more bifurcation is present in this case, see the bottom row of Fig.~\ref{F:ME-crit+supercrit} (at least in the range of  $b$ that we computed). 

\begin{wrapfigure}{l}{7.9cm}
\centering
\includegraphics[width=0.98\hsize, height=0.55\hsize]{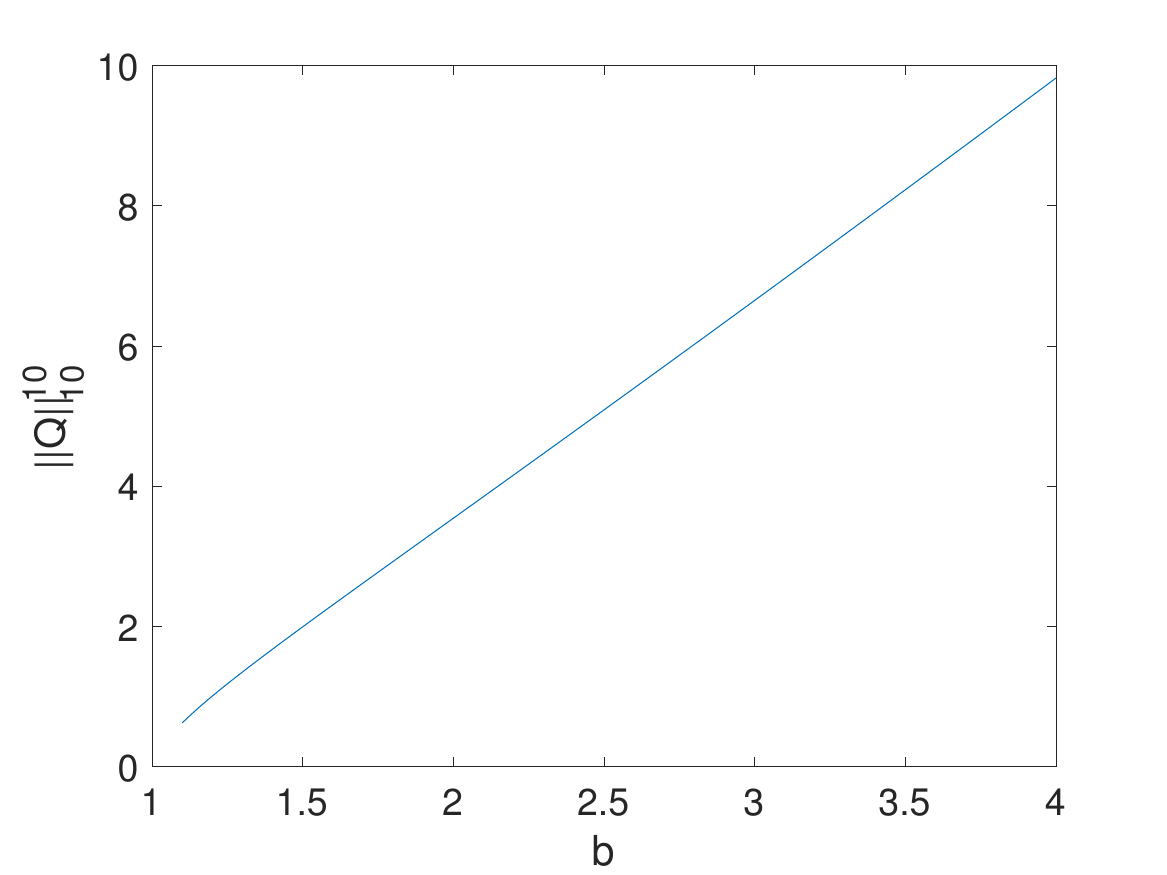} 
\caption{\footnotesize {Dependence of $\|Q^{(1)}\|_{L^{10}}^{10}$ on $b$ in the critical case $\alpha=8$.}} 
\label{F5}
\end{wrapfigure}

For completeness in the critical case, we also provide the dependence on $b$ of the potential term, $L^{10}$-norm, in Fig. \ref{F5}, which shows a linear dependence on $b$. 
\smallskip

Having examined the plethora of ground states, we proceed in the next section onto studying the time evolution of ground states and other solutions to the 
equation \eqref{biNLS}.

\clearpage



\section{Numerical approach for the time evolution}\label{S:evolv}

In this section we present the numerical approach for the study of the time evolution of solutions to \eqref{biNLS} and test it on an example of the stationary 
solutions, constructed in the previous section, and their time evolution. 

For the spatial discretization we use the same approach as in the previous section for equation \eqref{E:groundstate}, a Fourier 
spectral method. As before we consider functions sufficiently rapidly 
decreasing at infinity, i.e., mainly functions from the Schwartz class of rapidly decreasing smooth functions,  considering them on a torus with period $2\pi L$. Here, $L$ is again chosen large enough that the Fourier coefficients, for all initial 
data considered, decrease to machine precision. 
The FFT discretization in $x$ leads to \eqref{E:explicit1} being approximated 
via an equation of the form
\begin{equation}
\hat{u}_{t}=\hat{\mathcal{L}}\hat{u}+F[u]
\label{uhat2},
\end{equation}
where $\hat{\mathcal{L}}=-i(k^{4}-2ak^{2})$ is a diagonal linear 
operator, and where the nonlinear term reads 
$F[u]=i\widehat{|u|^{\alpha}u}$. Due to the fourth derivative in the 
linear term, the system \eqref{uhat2} is \emph{stiff}, which loosely 
speaking means that explicit time integration schemes are inefficient 
for stability reasons, see for instance \cite{HO} for a review of the 
subject and many references. An efficient approach to integrate such 
systems with a diagonal $\hat{\mathcal{L}}$ are so called \emph{exponential 
time differencing schemes}, see \cite{HO}. The idea is to introduce 
equidistant time steps $t_{n}$, $n=0,\ldots,N_{t}$, with 
$t_{n+1}-t_{n}=h$, the same constant $h$ for all $n=0,\ldots,N_{t}$. 
Integrating \eqref{uhat2} from $t_{n}$ to $t_{n+1}$ for some $n$, one 
gets
\begin{equation}
\hat{u}(t_{n+1}) = e^{\hat{\mathcal{L}}h}\hat{u}(t_{n}) + \int_{0}^{h}
e^{\hat{\mathcal{L}}(h-\tau)}F(\hat{u}(t_{n}+\tau))d\tau
	\label{ETD}.
\end{equation}
There are various approximations known in the literature to compute 
the integral in \eqref{ETD}, see \cite{HO}. As in \cite{etna} we apply here the Cox-Matthews scheme \cite{CM}, which is of classical 
order 4 (see discussion about classical order in \cite[p.212]{HO}), since there it was shown that ETD schemes proved to be very efficient for high order dispersive PDEs; we also note that the various schemes produce similar results, see \cite{KR}. The numerical accuracy is controlled as in \cite{etna} via the conserved 
quantities, mass \eqref{MC} and energy \eqref{EC}. These are exactly conserved 
by the equation, but unavoidable numerical errors will lead to a time 
dependence of their numerically computed counterparts. As discussed 
in \cite{etna}, the numerical conservation of these quantities overestimates the numerical error introduced by the time evolution by 
1-2 orders of magnitude and can thus be used to control the resolution in time. 

As a test for the quality of the code we use the solution constructed 
in the previous section for $b=2$, $a=1$ and $\alpha=8$, see 
Fig.~\ref{QA} on the left. Note that this is a non-trivial test, 
since the solution to \eqref{E:explicit1} for the initial data given by the respective $Q$ will have a harmonic time dependence, i.e.,  $e^{2it} Q(x)$. Moreover, the case $\alpha=8$ is critical, and even slight perturbations of $Q$ could cause blow-up in finite time (for example, if the mass is above the mass of $Q$ and $a=0$). If we choose $N_{t}=2000$ time steps for $t\in[0,1]$, then 
the energy is conserved to better  than $10^{-12}$ (similarly, for 
the mass), and the difference between the numerical solution and 
$Qe^{2it}$ is of the order of $10^{-12}$, the expected order of 
accuracy, with which $Q$ was constructed in the previous section. This 
confirms both the numerical 
accuracy of $Q$ and the time evolution code. In addition, it shows 
that the numerical conservation of the energy (and mass) is a valid 
indicator of the resolution in time. 

\section{Near soliton dynamics}\label{S:nearQ} 


In this section, we study small perturbations of the ground states in the mixed dispersion case ($a\neq 0$) to investigate their stability, since as we have seen branching occurs in the graph of energy vs. mass dependence $E = E(M)$ for solutions of certain nonlinearities in \eqref{E:1dGS}, recall the plots in 
Fig.~\ref{F:ME-subcritical} on the right. 
In our simulations we observe that the behavior of ground state perturbations varies significantly, depending on which branch of the $E(M)$ graph they are.
In particular, we identified  {\it stable} and {\it unstable} branches (in those cases where branching exists), which is unexpected, especially in the subcritical case. This means that for some (small) $b$ the ansatz $u(x,t) = e^{itb} Q(x)$ does not have a {\it stable} ground state solution, but rather produces an {\it unstable} state with the same mass as the stable ground state would be, however, with a different (larger) oscillation phase $\tilde b$, i.e., 
$$
\tilde b > b: \quad M[e^{itb}Q] = M[e^{it\tilde b} \tilde Q], \quad \quad E[e^{itb}Q]>E[e^{it \tilde b} \tilde Q].
$$    
Since the energy of the solution $u(t,x) = e^{ibt}Q(x)$ is higher than $u(t,x) = e^{i \tilde b t}\tilde Q(x)$, 
while the mass is the same, one expected an unstable behavior of the first solution. We discuss this in Section \ref{S:jumping}, showing several examples of such behavior. But first (and for comparison later), we discuss the cases, where no branching occurs. 



In simulations here, we use $N=2^{10}$ Fourier modes for $x\in [-50\pi,50\pi]$ and $N_{t}=1000$ time steps for $t\in [0,10]$. 
The Fourier coefficients decrease for all studied examples in this section to 
the order of $10^{-10}$ at least, and the discretized energy is conserved in
relative order to better than $10^{-10}$.

\subsection{Near soliton dynamics, non-branching case}\label{S:nobranch} 

We start with considering the sub-critical case of \eqref{E:explicit1} with  
$\alpha=2$, since this nonlinearity did not show any branching in our investigation of the energy vs. mass $E(M)$ dependence of solutions to \eqref{E:1dGS}, in Fig.~\ref{F:ME-subcritical}.

We fix $b=2$ and take initial data as the ground state perturbations of the form 
\begin{equation}
u(x,0) = A \, Q(x),\quad A \sim 1.
\label{uini}
\end{equation}

$\scriptstyle\blacklozenge$ ({$a=0$}) We begin with the pure quartic case $a=0$, since it is a scaling-invariant case, and hence, the ground state \eqref{Qscaling} scales as well. 
\begin{figure}[htb!]
\begin{subfigure}{.45\textwidth}
\includegraphics[width=1\linewidth,height=0.67\linewidth]{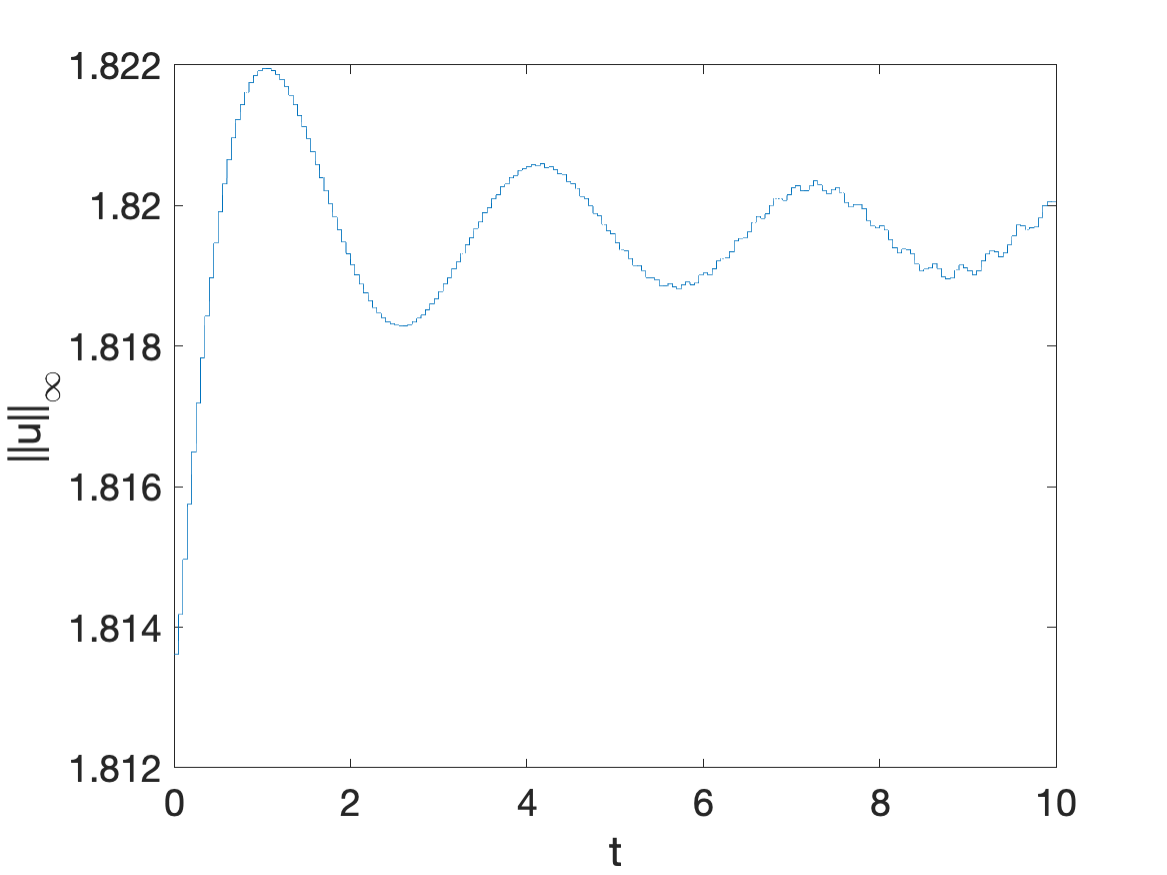} 
\subcaption[]{{\footnotesize $L^\infty$ norm}}
\end{subfigure}
\begin{subfigure}{.45\textwidth}
  \includegraphics[width=1\linewidth,height=0.67\linewidth]{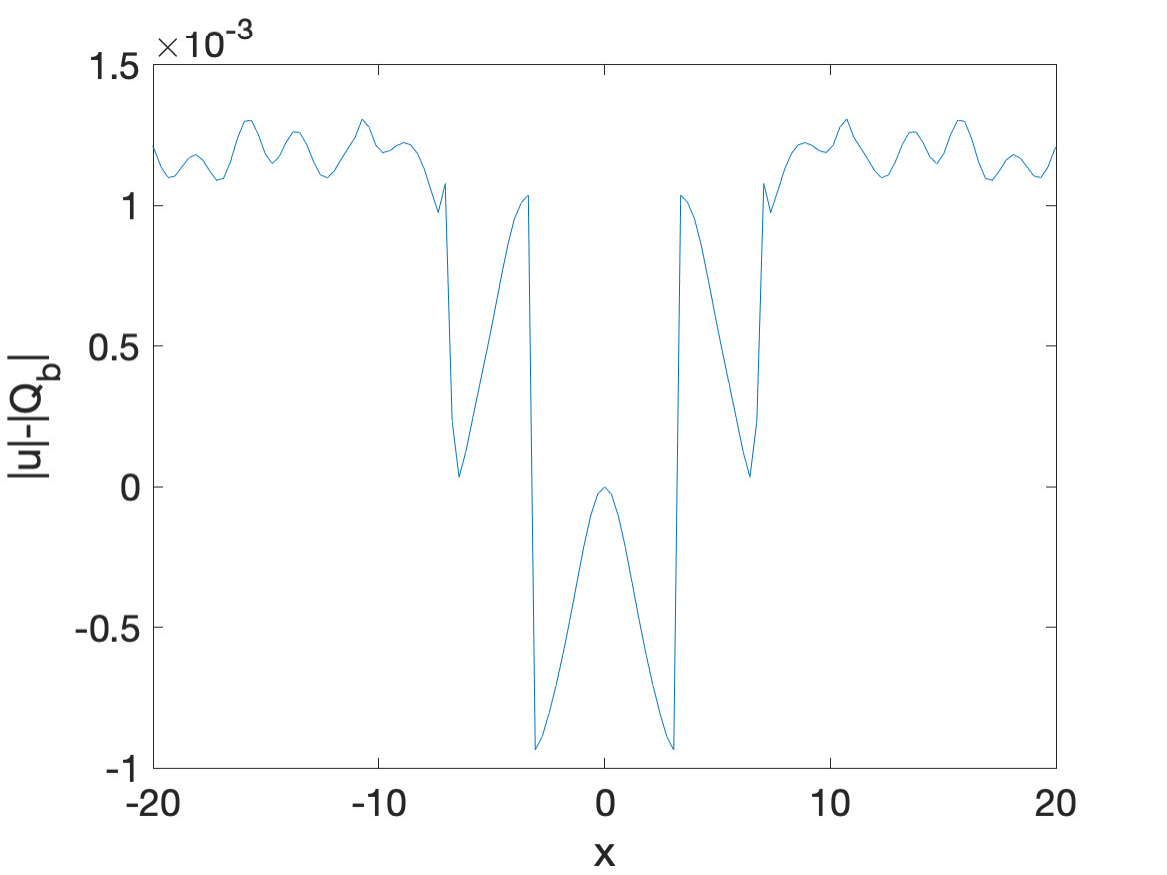}
\subcaption[]{{\footnotesize difference of $|u|$ and rescaled $Q_b$}}
\end{subfigure}\\
\begin{subfigure}{.45\textwidth}
\includegraphics[width=1\linewidth,height=0.67\linewidth]{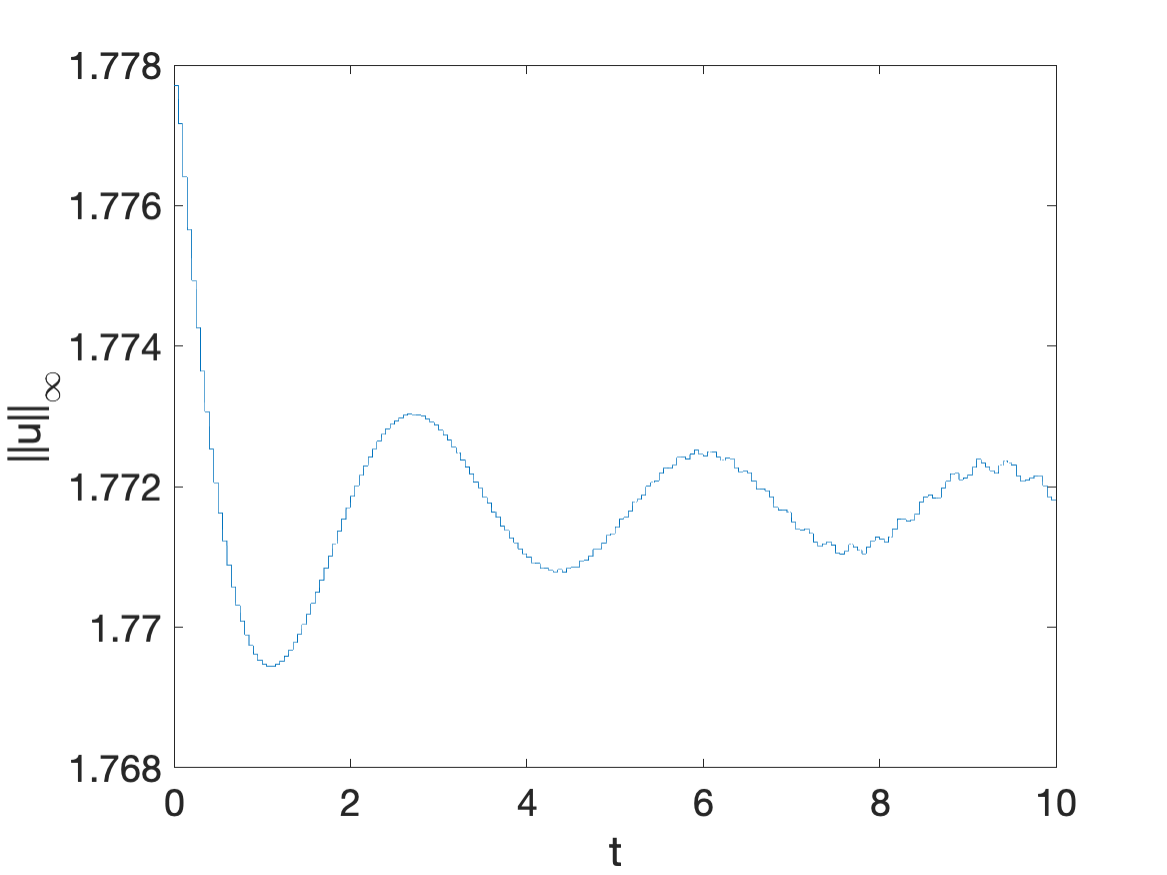}
\subcaption[]{{\footnotesize $L^\infty$ norm}}
\end{subfigure}
\begin{subfigure}{.45\textwidth}
\includegraphics[width=1\linewidth,height=0.67\linewidth]{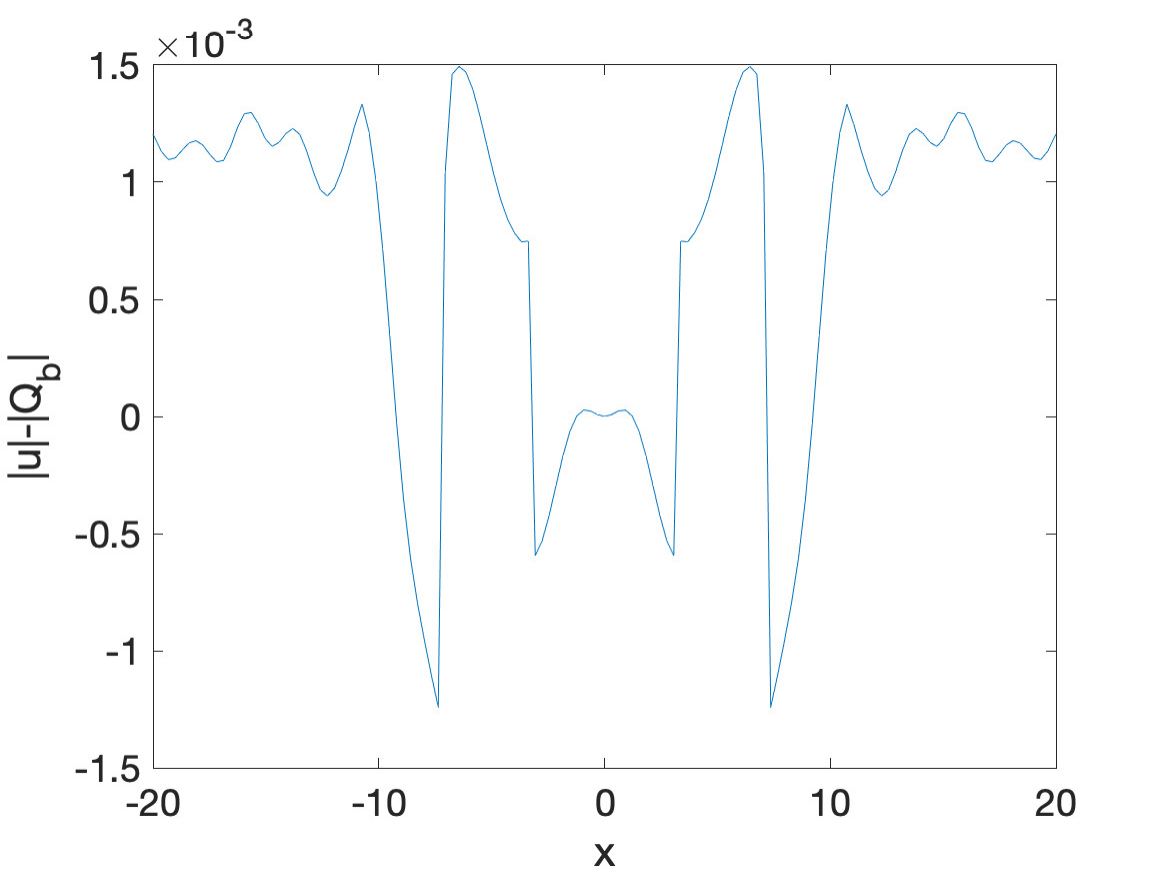}
\subcaption[]{{\footnotesize difference of $|u|$ and rescaled $Q_b$}}
\end{subfigure}
\caption{\footnotesize {Solution to \eqref{E:explicit1}, $\alpha=2$, $a=0$, $b=2$ with $u_0=1.01 Q$ (top) and $u_0=0.99 Q$ (bottom). Left: time dependence of the $L^{\infty}$ norm. Right: the difference of the modulus of the solution at $t=10$  and the rescaled ground state $Q_b$.}}
\label{F:AQ0-alpha2}
\end{figure}
The $L^{\infty}$~norm of the solution with the initial data \eqref{uini} and
$A=1.01$ is shown in plot (A) on the left top of Fig.~\ref{F:AQ0-alpha2}. It can be seen 
that the $L^{\infty}$~norm of the solution grows slightly  and then 
oscillates, approaching in amplitude a level of the value slightly higher than 
the unperturbed ground state. 
This is due to the fact that while analytically the perturbations are infinitesimally small, numerically we have to  
consider a finitely small perturbation in order to observe a visible 
effect of the perturbation (in finite time).
The perturbation leads 
to a ground state, but of a slightly higher $L^\infty$ norm. 
A fit to the ground state 
for a parameter $b$ according to \eqref{Qscaling}, namely,
$Q_{b} = (b/2)^{1/\alpha}Q((b/2)^{1/4}x)$, is done in the following way:
$(b/2)^{1/\alpha}$ is given by the maximum of the modulus of the solution at the 
final time divided by the maximum of $Q$. The difference between 
the modulus of the solution at the final time and the fitted ground state $|Q_b|$ is shown on the right of the same figure. It is on the order of $10^{-3}$, since the final state is not yet reached. 
Note there will always be oscillations around the final state, since there is no dissipation in the systems (see for instance \cite{CKS}), and since we are working on the torus, no radiation can escape to infinity. 
In fact, the  size of the torus has to be chosen large enough in order to limit the effects of radiation 
reentering the computational domain, see also \cite{AKS,KSM} for similar 
perturbations of the NLS ground states. The jumps in the $L^{\infty}$ 
norm are due to the fact that it is determined at the collocation 
points of the FFT, and thus, the result depends on whether the actual 
maximum of $|u|$ is on a grid point or not. 

The situation is similar for a factor $A=0.99$. 
As can be seen on the bottom left of Fig.~\ref{F:AQ0-alpha2}, or plot (C), the $L^{\infty}$~norm of the 
solution decreases slightly at the beginning and then decreasingly 
oscillates in amplitude around a slightly smaller (shorter in height) 
ground state. The difference of the modulus of the solution at the final time with the fitted ground state $Q_b$ (for a fitted $b$ as discussed above) can be seen on the right of the same figure, plot (D). 
The difference is again on the order of $10^{-3}$, which provides a good numerical confirmation about the asymptotic approach to a soliton state, and thus, in a sense a soliton resolution. 
\smallskip

\begin{figure}[!htb]
\begin{subfigure}{.45\textwidth}
\includegraphics[width=1\linewidth,height=0.7\linewidth]{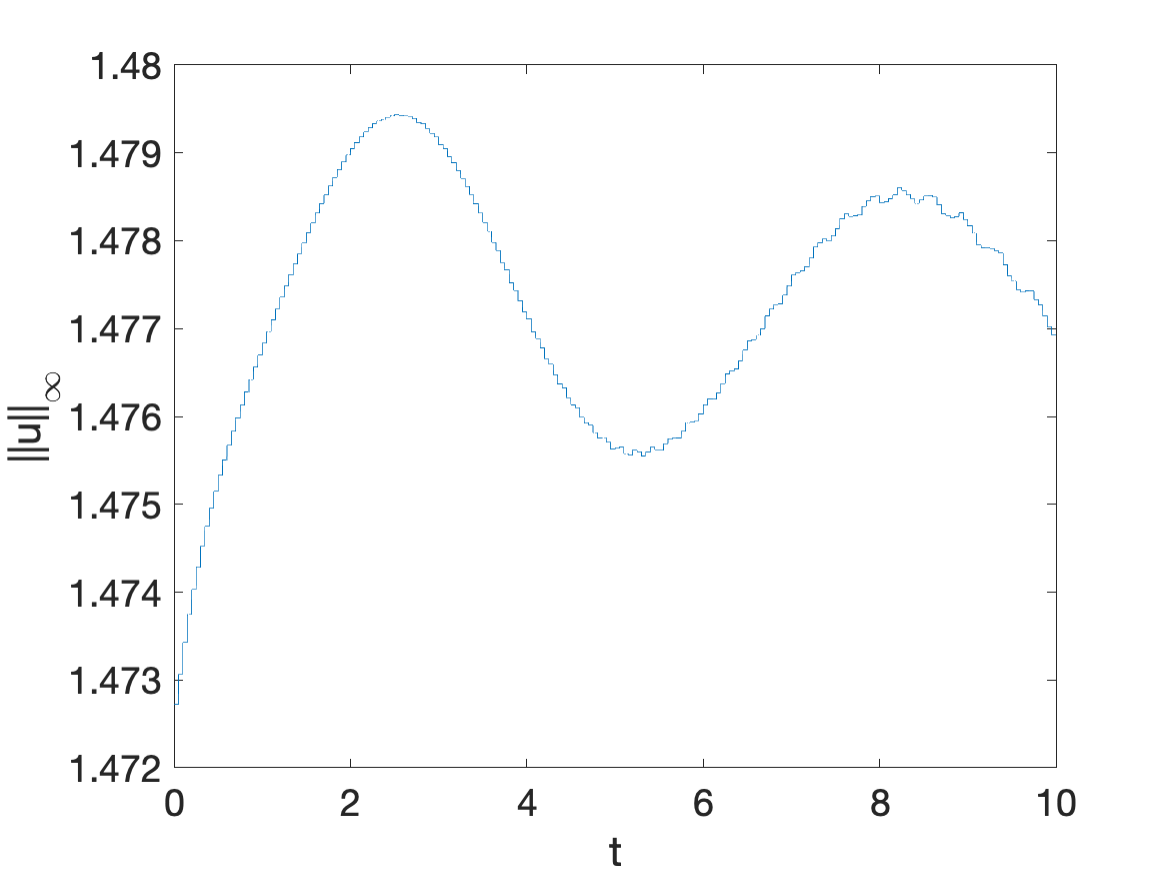} 
\subcaption[]{{\footnotesize $u_0=1.01 Q^{(1)}$ }}
\end{subfigure}
\begin{subfigure}{.45\textwidth}
  \includegraphics[width=1\linewidth,height=0.7\linewidth]{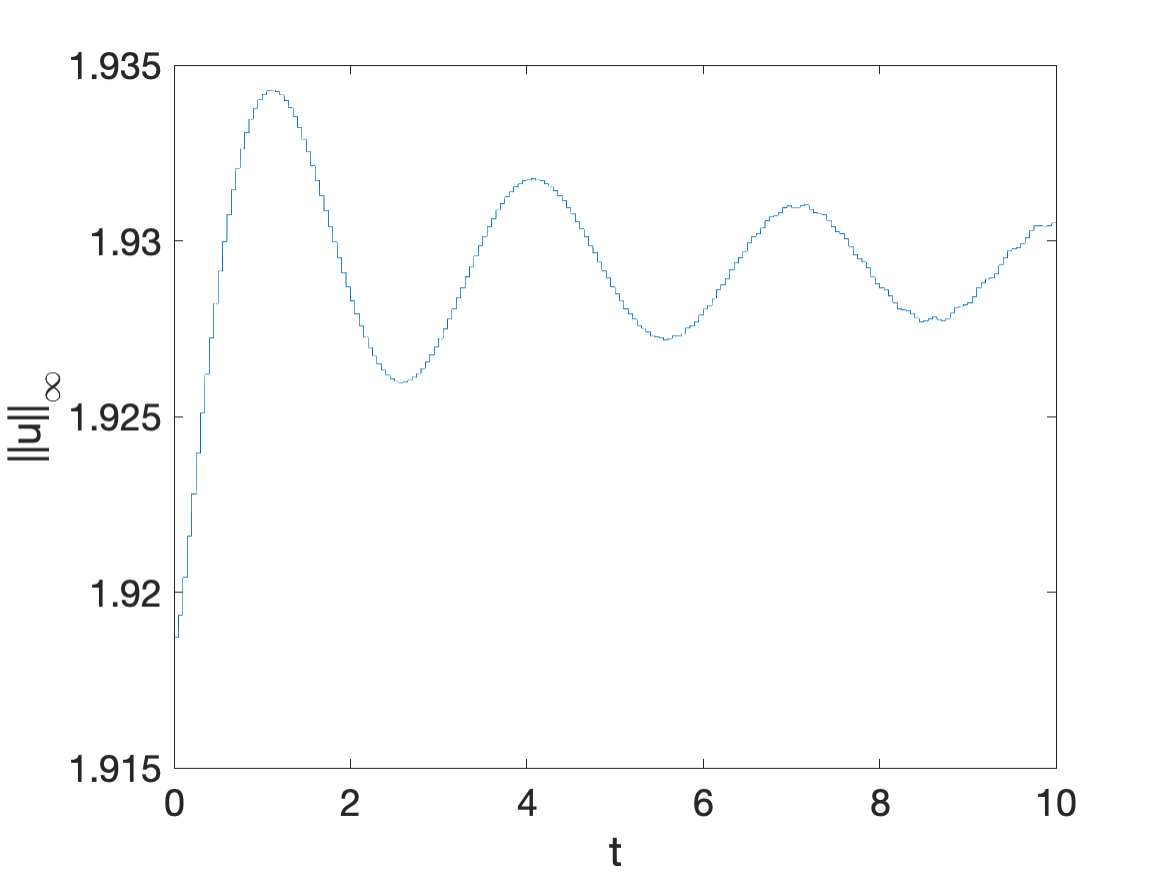}
\subcaption[]{{\footnotesize $u_0=1.01 Q^{(-1)}$}}
\end{subfigure}\\
\begin{subfigure}{.45\textwidth}
\includegraphics[width=1\linewidth,height=0.7\linewidth]{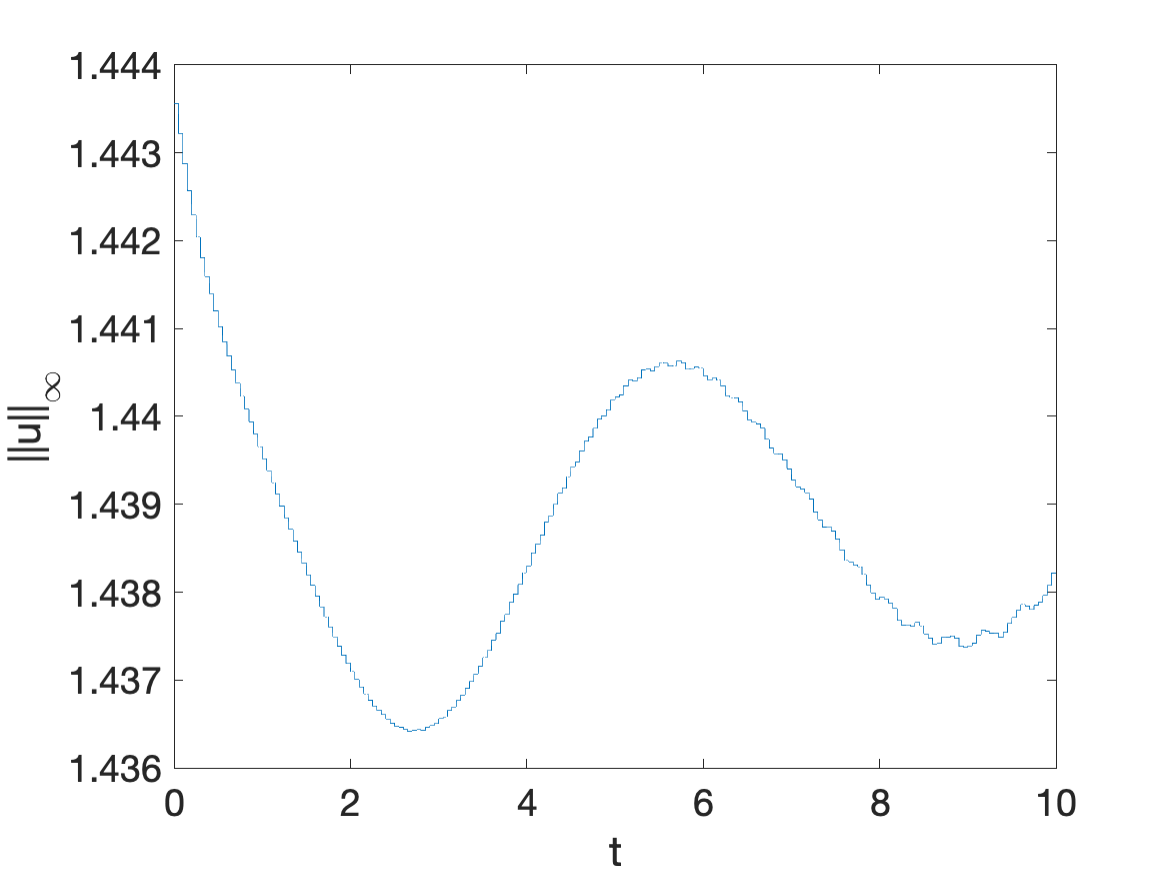}
\subcaption[]{{\footnotesize $u_0=0.99 Q^{(1)}$}}
\end{subfigure}
\begin{subfigure}{.45\textwidth}
\includegraphics[width=1\linewidth,height=0.7\linewidth]{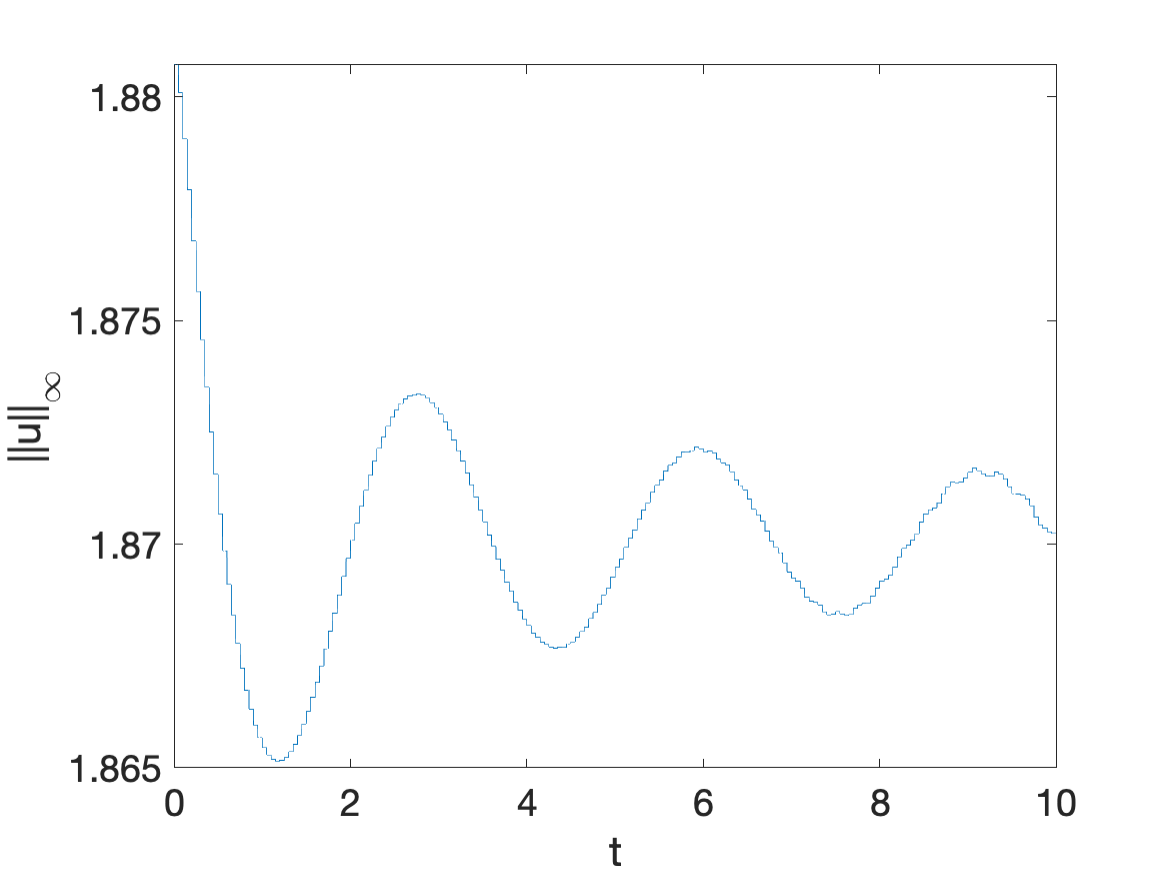}
\subcaption[]{{\footnotesize $u_0=0.99 Q^{(-1)}$}}
\end{subfigure}
\caption{\footnotesize {Time dependence of the $L^\infty$~norms of 
solutions to \eqref{E:explicit1}, $\alpha=2$, $b=2$, with $u_0=A Q^{(a)}$ with $A=1.01$ (top), $A=0.99$ (bottom). Left: $a=1$. Right: $a=-1$.}}
\label{F:AQ0-alpha2-a}
\end{figure}

$\scriptstyle\blacklozenge$ If $a\neq0$, there is no simple scaling 
in $b$ of the form 
\eqref{Qscaling} as it exists when $a=0$. Thus, we can only consider the $L^{\infty}$ norm in these cases. In the top row of Fig.~\ref{F:AQ0-alpha2-a} we show the $L^{\infty}$~norms of the solution $u(t)$ with initial data \eqref{uini}, $u_0 = A\, Q^{(a)}$, $A=1.01$, for $a=\pm1$. The behavior is similar (converging with oscillations) to what is shown on the left of Fig.~\ref{F:AQ0-alpha2}. 

On the bottom row of Fig.~\ref{F:AQ0-alpha2-a} we show solutions with perturbations by $A=0.99$ and $a=\pm1$. 
The behavior is similar to the case of $a=0$ shown on the left of Fig.~\ref{F:AQ0-alpha2}.

We conclude that the ground states in the case of subcritical nonlinearity $\alpha=2$ exhibit a stable behavior for different values of parameter $a$. 
As no branching of energy occurs (at least in the range of parameters that we considered), 
it also shows {\it how a stable ground state behaves under small perturbations} in the case of combined dispersions, 
that is, small amplitude changes in $A$ force the solution to slightly increase or decrease in order to `find' the rescaled version of itself and asymptotically (and in oscillatory manner with very small amplitude oscillations) approach it. We refer to this behavior as the {\it stable} perturbation of the ground state, when considering different branches of ground states below. 




\subsection{Near soliton dynamics, branching case}\label{S:jumping}
Here we consider the nonlinearities where the branching of the energy was observed, recall Fig.~\ref{F:ME-subcritical} and \ref{F:ME-crit+supercrit}. 
\begin{figure}[!htb]
\includegraphics[width=0.45\hsize, height=0.28\hsize]{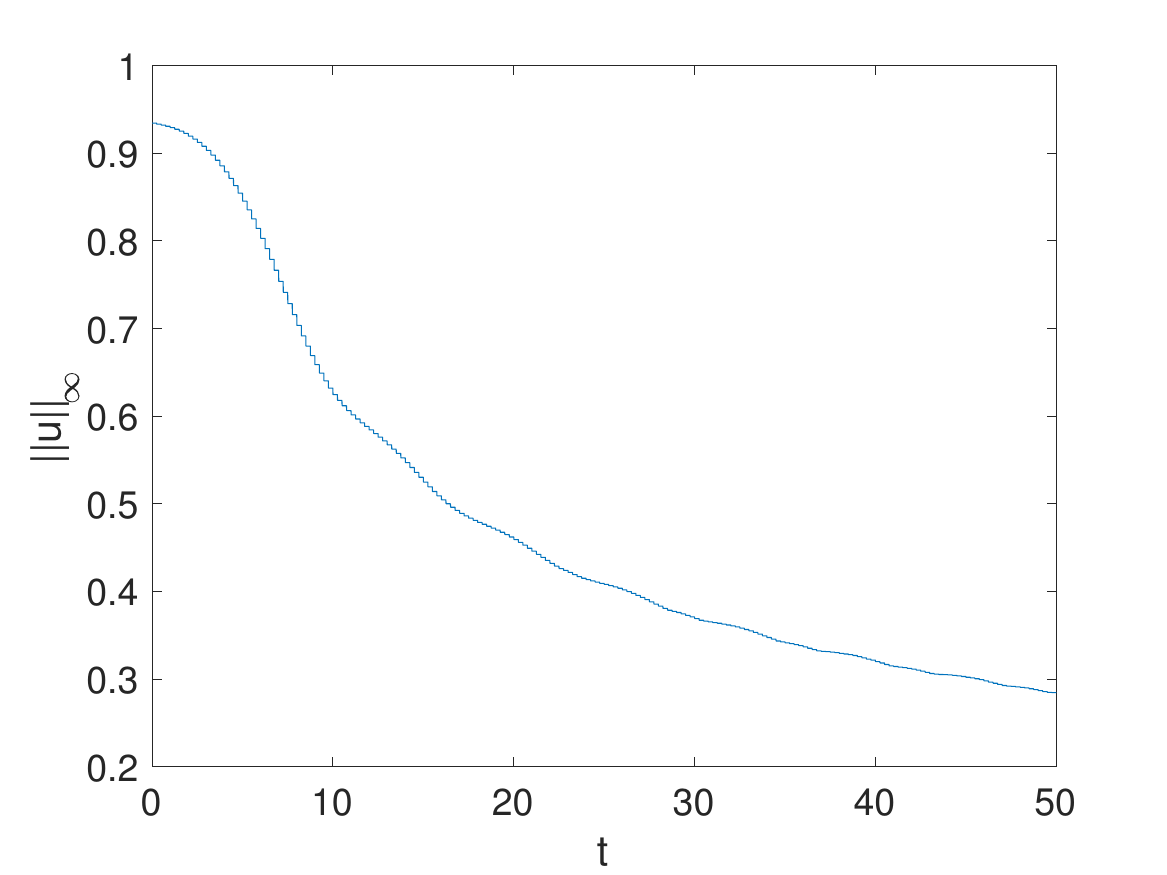}
\includegraphics[width=0.45\hsize, height=0.28\hsize]{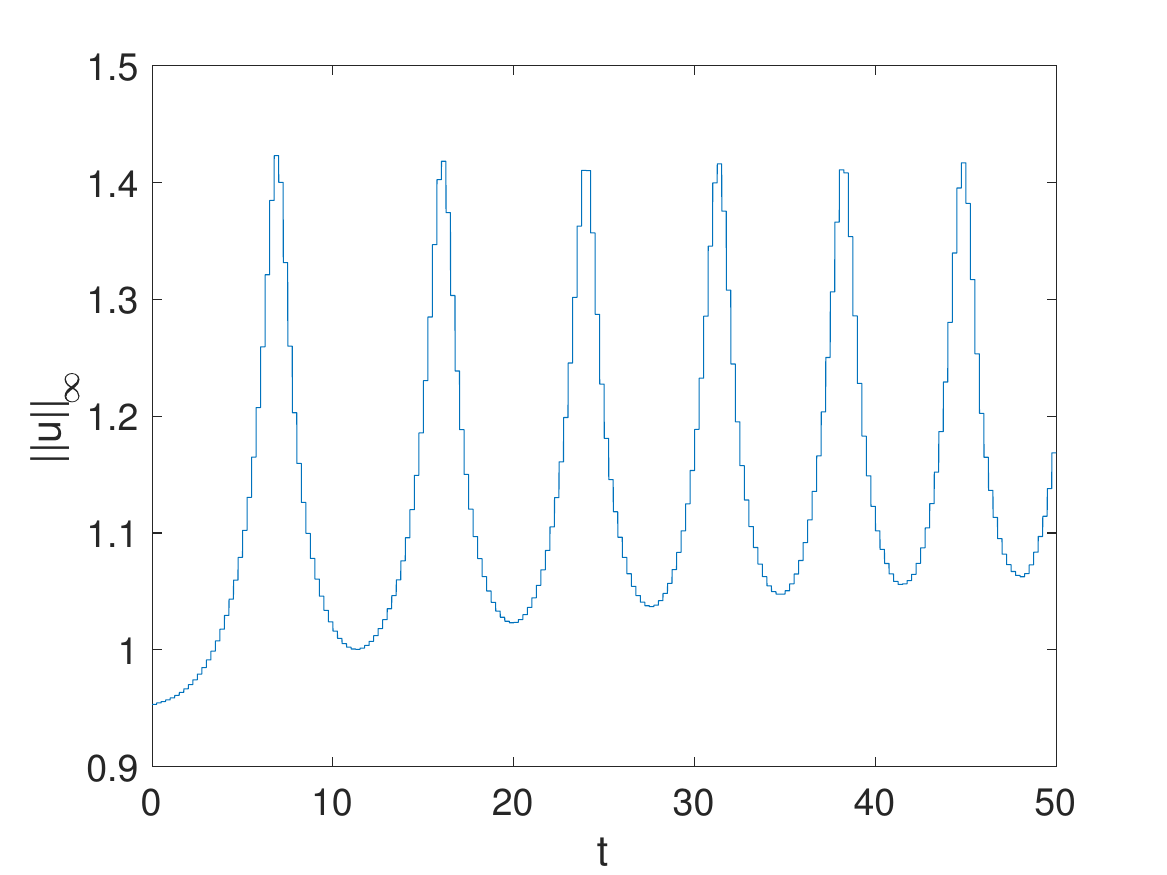}
\caption{\small {Time dependence of the $L^\infty$ norm of solutions to \eqref{E:explicit1}, $\alpha=6$, $a=1$ with $u_0=.99Q^{(1)}$ (left) and  $u_0=1.01Q^{(1)}$ (right) for $b=1.1$. }}
\label{F:near-alpha6}
\end{figure}

\begin{figure}[!htb]
\begin{subfigure}{.45\textwidth}
\includegraphics[width=1\linewidth,height=0.65\linewidth]{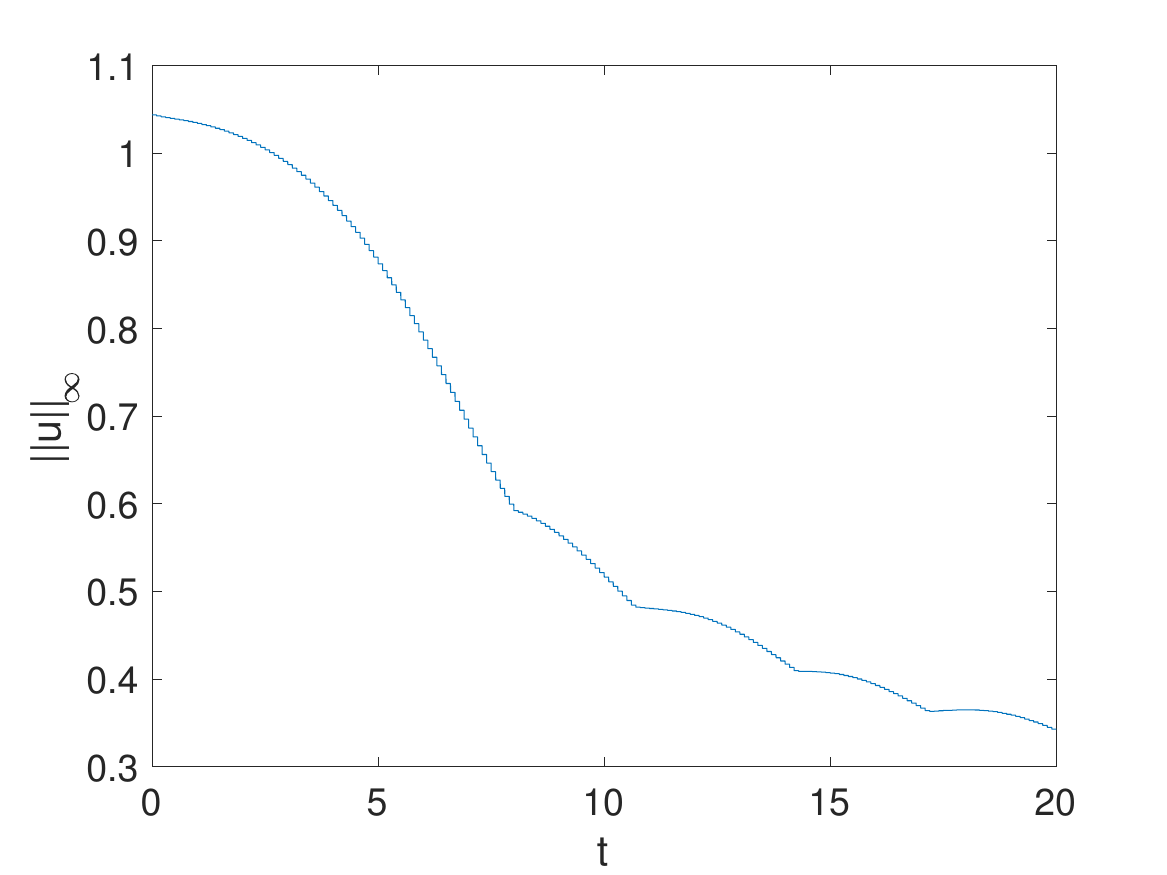} 
\subcaption[]{{\footnotesize $u_0=0.99Q^{(1)}$, $b=1.3$}}
\end{subfigure}
\begin{subfigure}{.45\textwidth}
  \includegraphics[width=1\linewidth,height=0.65\linewidth]{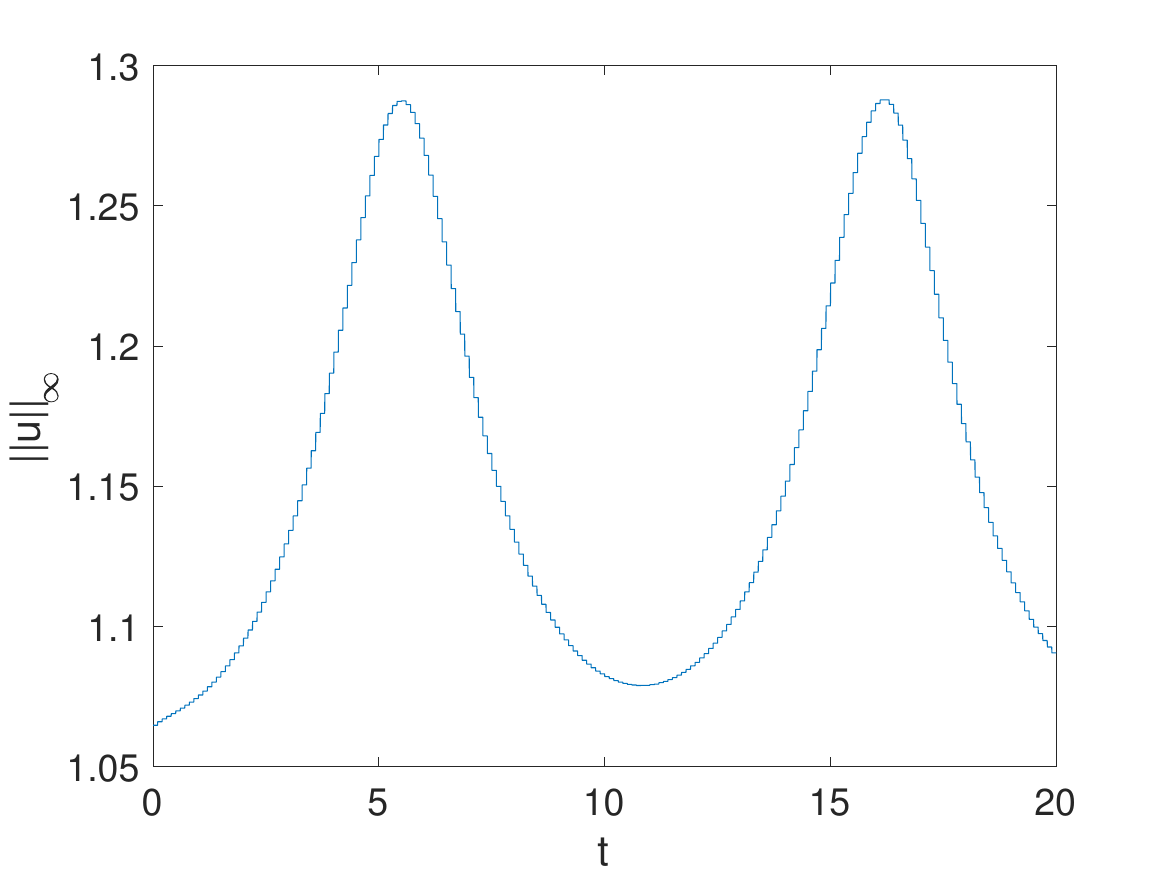}
\subcaption[]{{\footnotesize $u_0=1.01Q^{(1)}$, $b=1.3$}}
\end{subfigure}\\
\begin{subfigure}{.45\textwidth}
\includegraphics[width=1\linewidth,height=0.65\linewidth]{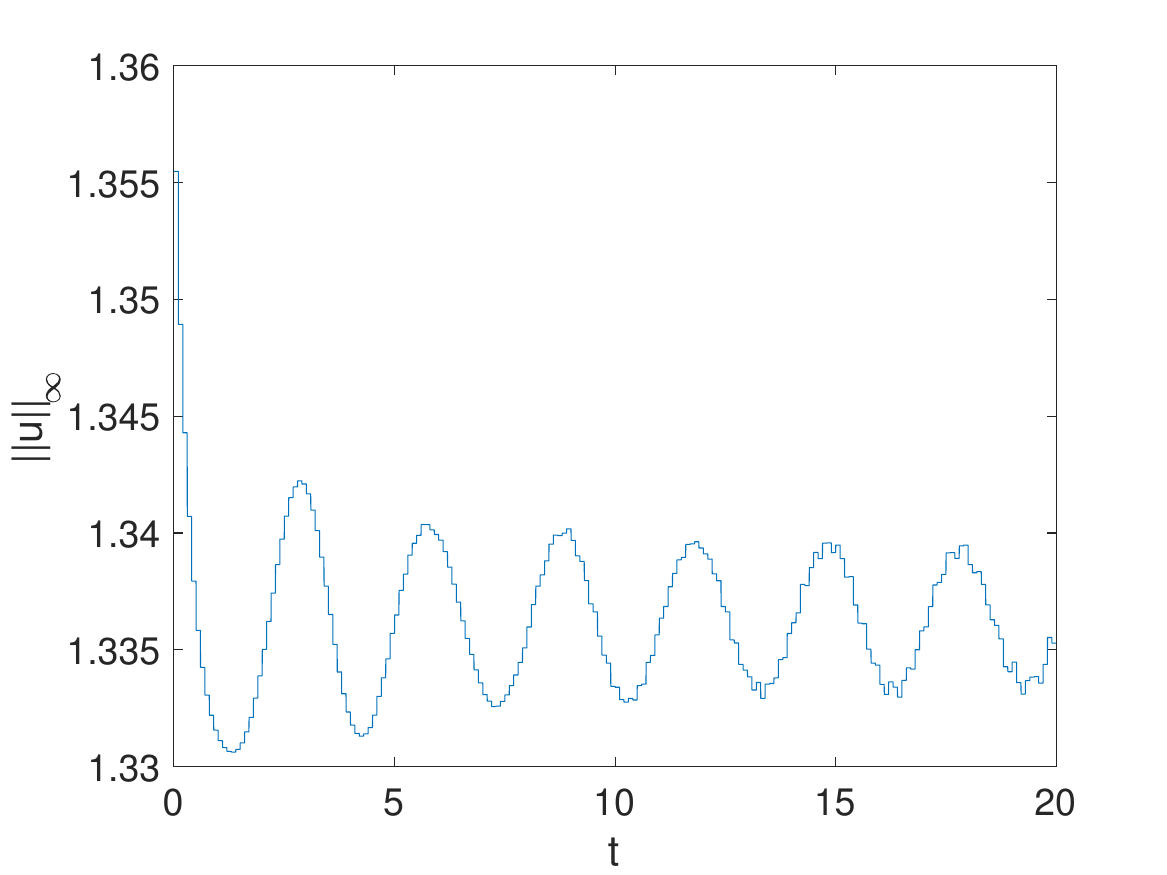}
\subcaption[]{{\footnotesize $u_0=0.99Q^{(1)}$, $b=3.5$}}
\end{subfigure}
\begin{subfigure}{.45\textwidth}
\includegraphics[width=1\linewidth,height=0.65\linewidth]{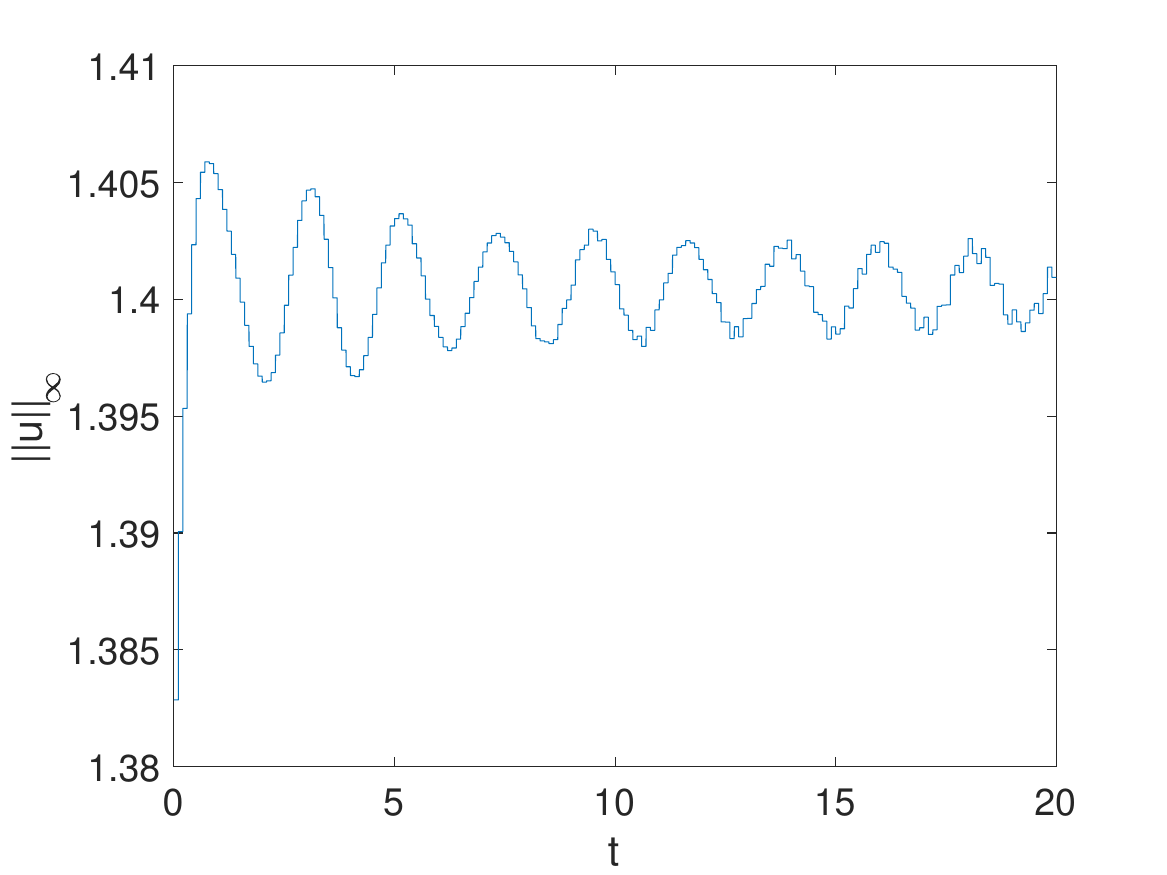}
\subcaption[]{{\footnotesize $u_0=1.01Q^{(1)}$, $b=3.5$}}
\end{subfigure}
\caption{\footnotesize {Solution to \eqref{E:explicit1}, $\alpha=6$, $a=1$ with $b=1.3$ (top), $b=3.5$ (bottom), 
and $u_0=0.99 Q^{(1)}$ (left), $u_0=1.01 Q^{(1)}$ (right). 
The initial data are chosen such that the mass is approximately the same, but the ground state in the top row has higher energy, thus, unstable.}}
\label{F:AQ0-alpha6-mass2.9}
\end{figure}

$\scriptstyle\blacklozenge$ {\it Subcritical case.}
We take $\alpha=6$, fix $a=1$, and consider the initial data $u_0=A\, Q^{(1)}$ with different values of $b$ and $A \sim 1$.  
In our first example, we consider $b=1.1$, noting that the mass is relatively large, see Fig.~\ref{F:ME-subcritical}(G), 
therefore, we expect this ground state solution to be unstable, which is indeed the case, shown in Fig.~\ref{F:near-alpha6}. 
We point out that the behavior of these perturbations is qualitatively different from those in Fig.~\ref{F:AQ0-alpha2} and \ref{F:AQ0-alpha2-a}, 
the perturbations with $A<1$ 
scatter to zero, those with $A>1$ create large oscillations, of different nature, than those shown in Fig.~\ref{F:AQ0-alpha2}(A), \ref{F:AQ0-alpha2-a}(A)\&(B), and seem to approach a different state (a stable ground state). 

To further confirm the unstable behavior and the unstable branch, we take the initial data, of approximately the same mass, 
and compare two solutions. For example, fixing the mass around $2.8$, 
we find that two values of $b$ can produce the ground state solutions, 
namely, $b\sim 1.3$ and $b\sim 3.5$. We fix these values, compute their energies, and then run the simulations to compare the behavior of perturbations. 
We have for $b=1.3$, $M[Q^{(1)}] = 2.7$, $E[Q^{(1)}]=-1.23$ and for $b=3.5$, $M[Q^{(1)}] = 2.9$, $E[Q^{(1)}]=-1.59$. Since the energy is larger for smaller $b$, we conclude that the ground state for a smaller $b$ value should be unstable. We confirm this in Fig.~\ref{F:AQ0-alpha6-mass2.9}.  
\smallskip

\begin{figure}[!htb]
\begin{subfigure}{.32\textwidth}
\includegraphics[width=1\linewidth,height=0.7\linewidth]{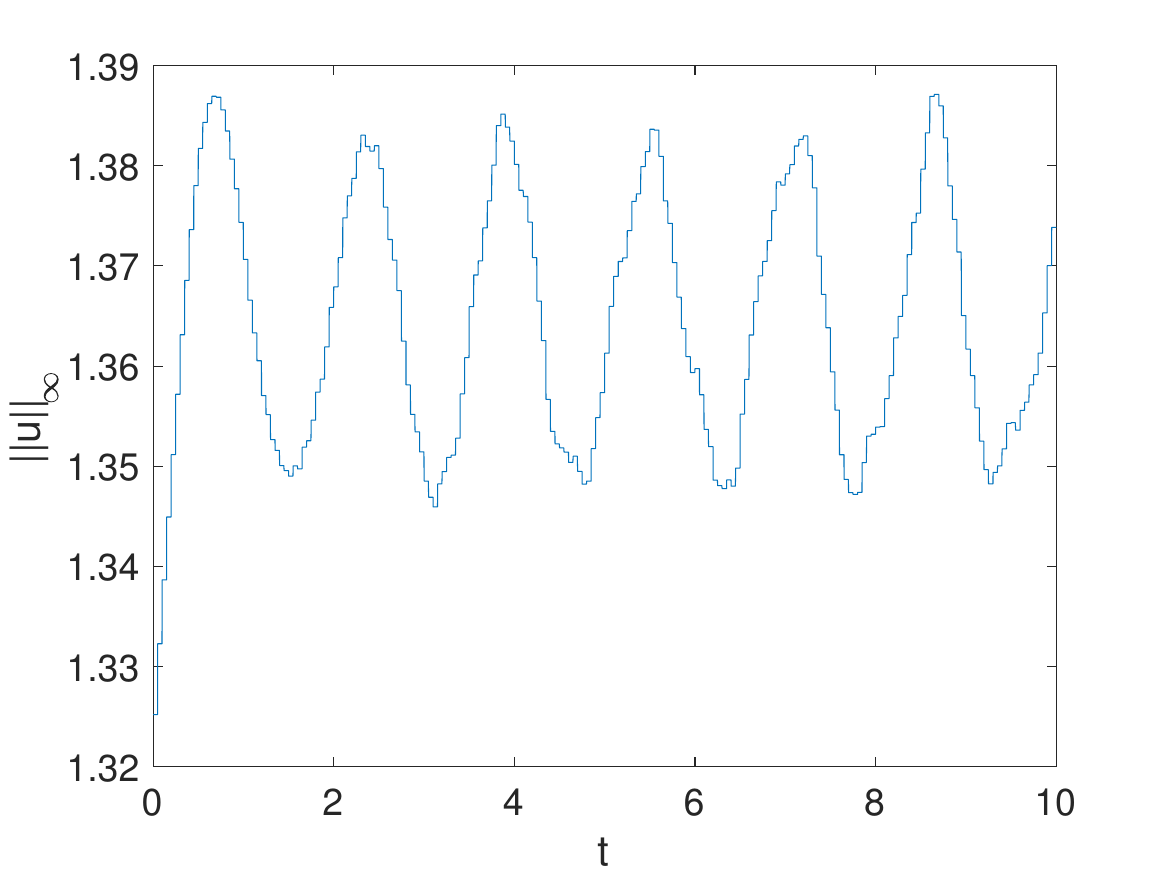}
\subcaption[]{{\footnotesize $u_0=1.01Q^{(1)}$, $b=4$}}
\end{subfigure}
\begin{subfigure}{.32\textwidth}
  \includegraphics[width=1\linewidth,height=0.7\linewidth]{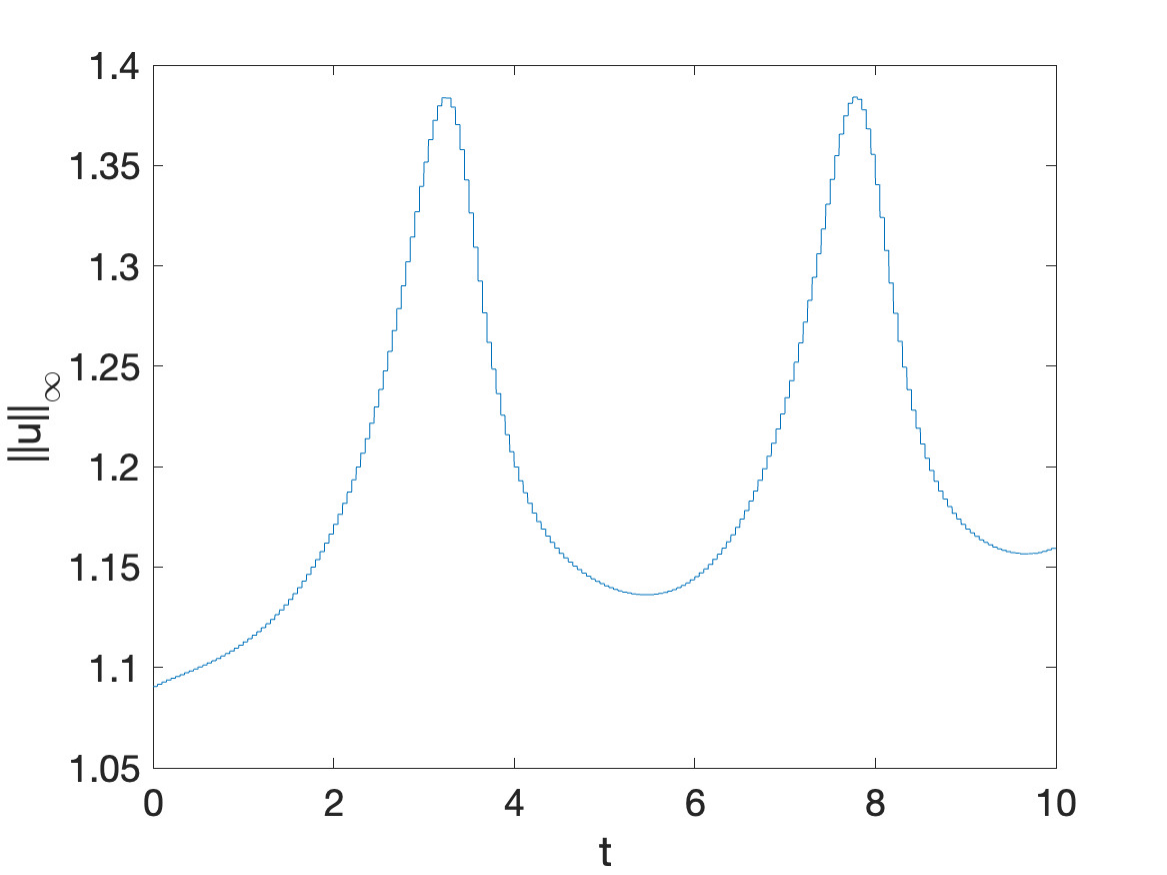}
\subcaption[]{{\footnotesize $u_0=1.01Q^{(1)}$, $b=1.4$}}
\end{subfigure}
\begin{subfigure}{.32\textwidth}
\includegraphics[width=1\linewidth,height=0.7\linewidth]{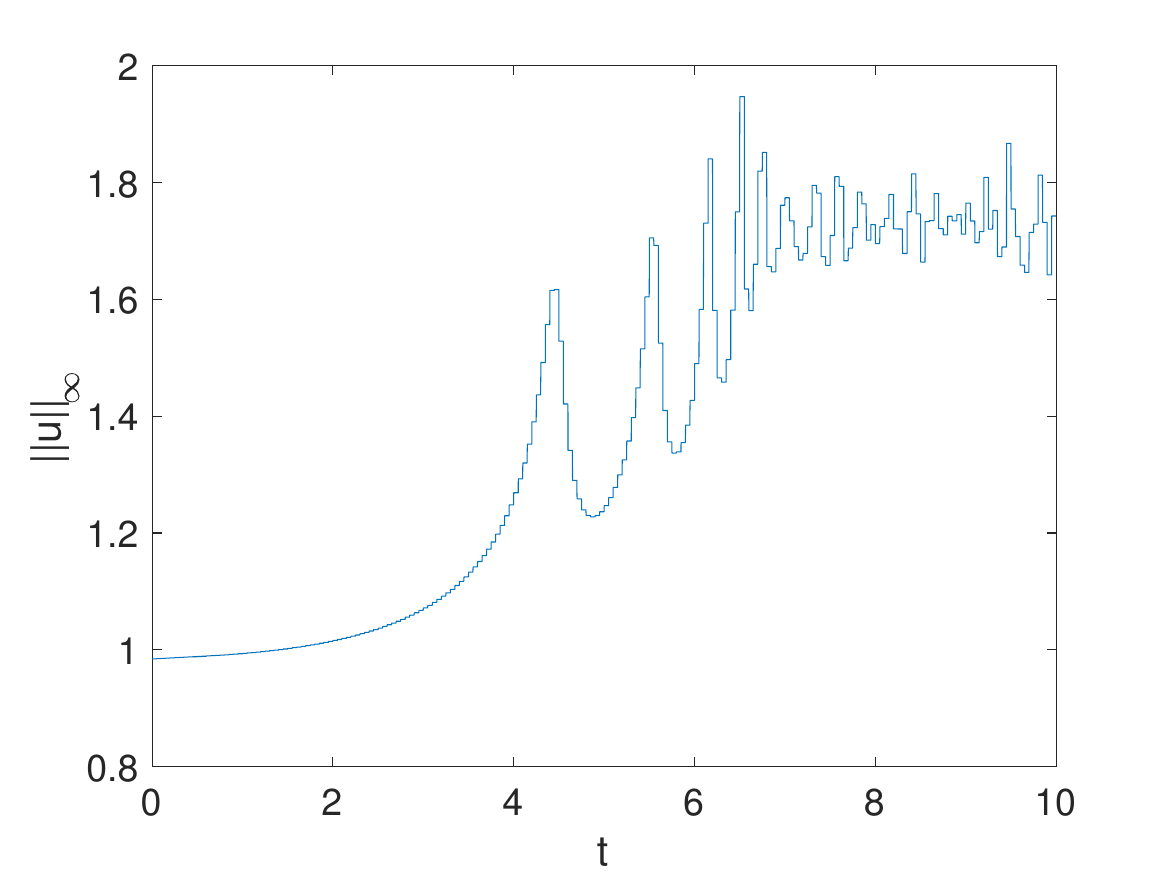}
\subcaption[]{{\footnotesize $u_0=1.01Q^{(1)}$, $b=1.1$}}
\end{subfigure}
\begin{subfigure}{.32\textwidth}
\includegraphics[width=1\linewidth, height=0.7\linewidth]{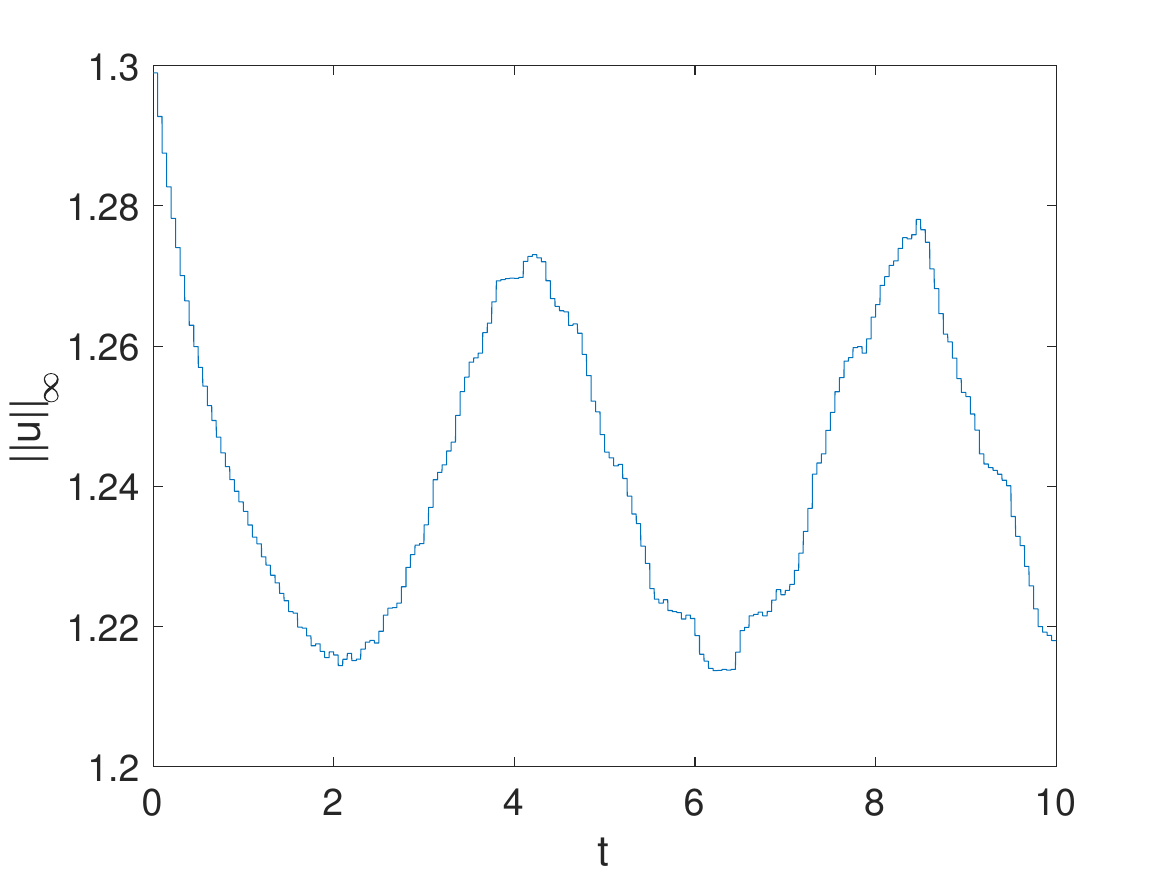}
\subcaption[]{{\footnotesize $u_0=0.99Q^{(1)}$, $b=4$}}
\end{subfigure}
\begin{subfigure}{.32\textwidth}
\includegraphics[width=1\linewidth, height=0.7\linewidth]{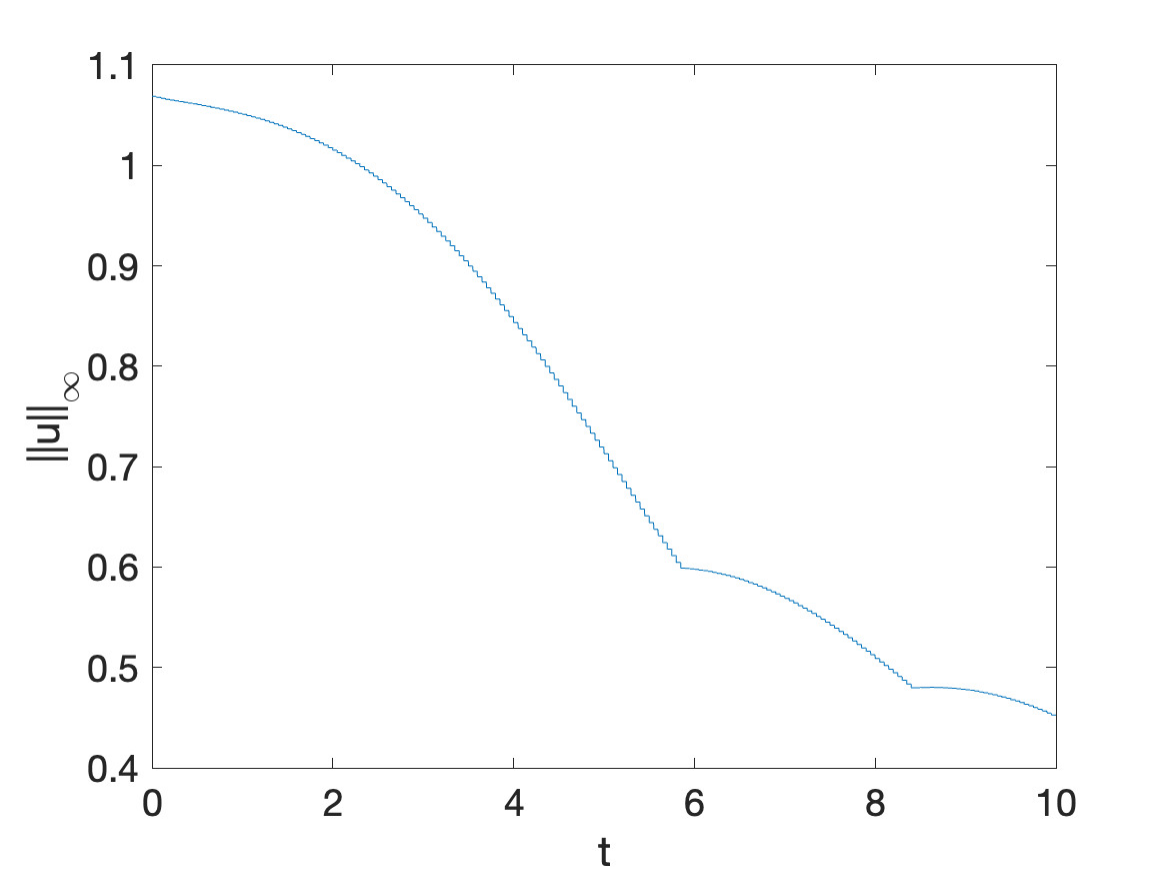}
\subcaption[]{{\footnotesize $u_0=0.99Q^{(1)}$, $b=1.4$}}
\end{subfigure}
\begin{subfigure}{.32\textwidth}
\includegraphics[width=1\linewidth, height=0.7\linewidth]{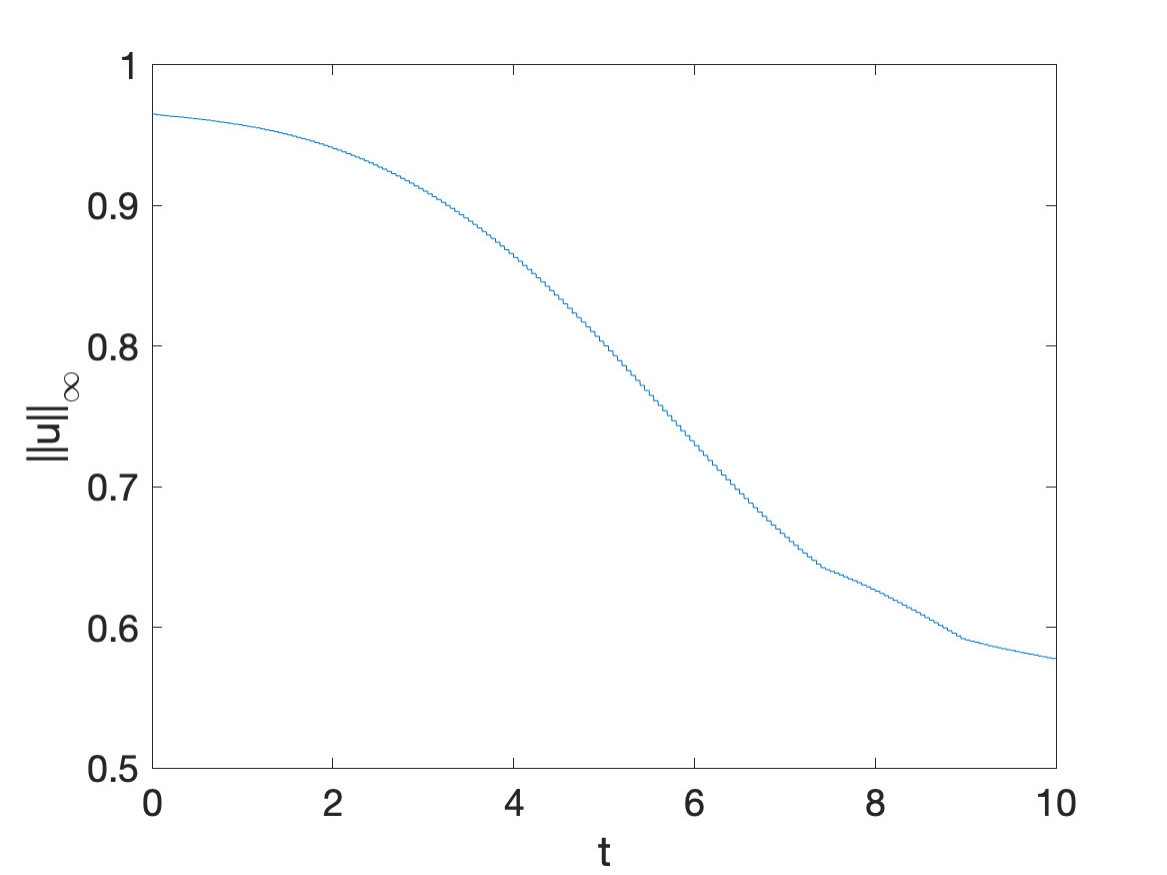}
\subcaption[]{{\footnotesize $u_0=0.99Q^{(1)}$, $b=1.1$}}
\end{subfigure}
\caption{\footnotesize {Solution to \eqref{E:explicit1}, $\alpha=8$, $a=1$, with 
initial data $u_0$ and the parameter $b$ as indicated.
}}
\label{F:alpha8branch}
\end{figure}

$\scriptstyle\blacklozenge$ {\it Critical case.} 
We also study the critical case, $\alpha=8$, to confirm existence of stable and unstable branches of ground states. 
We continue with $a=1$ and consider three cases of the parameter $b = 1.1, 1.3, 4$. The first two values produce the ground state with the energies on the lower branch in Fig.~\ref{F:ME-crit+supercrit}(C), and thus, are expected to be unstable. We confirm that in Fig.~\ref{F:alpha8branch}.

\smallskip

One can observe a stable behavior for $b=4$ in the left column of Fig.~\ref{F:alpha8branch}, with small oscillations (eventually converging to some final state); however, the other two values of $b$ show the unstable behavior with dispersing down to zero for amplitude perturbations $A<1$ and having large oscillatory behavior for $A>1$. 
In particular, one can notice the convergence in the right top of Fig.~\ref{F:alpha8branch} to a higher level, indicating the convergence to a stable ground state.  
We also note that in the critical case, we do not see formation of blow-up right after the value $A=1$, 
which is different from the classical case. Thus, this is a qualitatively different behavior of ground states from the classical (as in the NLS, e.g., see \cite{HR2008}) ground state sharp threshold for blow-up vs. scattering behavior. This happens in the non-scaling invariant cases ($a \neq 0$), since the lower order dispersion (with either sign) cannot be scaled away.

\section{Dynamics of generic solutions in the subcritical case}\label{S:subcritical-general} 

In this section, we further investigate the subcritical case of \eqref{E:explicit1}, in particular, 
we study behavior of generic initial data of Schwartz class of rapidly decreasing smooth functions, including different Gaussian-type data. 

We show that solutions with sufficiently large amplitude of Gaussian data confirm the soliton resolution conjecture, splitting into a final state rescaled soliton and radiation; for smaller amplitudes solutions disperse to zero. 

\subsection{Behavior of solutions for Gaussian initial data}
We start with  
the following  Gaussian 
initial data 
\begin{equation}\label{E:gaussian}
u_0(x)=A e^{-x^{2}}, \quad A>0.
\end{equation}
(Note that $M[u_0]=A^2 \sqrt{\frac{\pi}{2}}$.)

The computation is done on large tori in order to avoid emission of radiation 
towards one boundary of the computational domain and its reappearance on the other 
side (due to the periodicity of the setting), which may have a strong effect on the results. But since, as can be seen below, there is a large amount of radiation, this cannot be completely avoided even on large tori. 
We consider $x\in[-100\pi,100\pi]$ with $N^{12}$ Fourier modes and 
$N_{t}=10^{4}$ time steps for $t\in[0,10]$. The Fourier coefficients 
decrease to machine precision during the whole computation, and the discretized
energy is conserved in relative error 
to the order of $10^{-9}$.

\begin{figure}[htb!]
 \includegraphics[width=0.36\textwidth,height=0.26\textwidth]{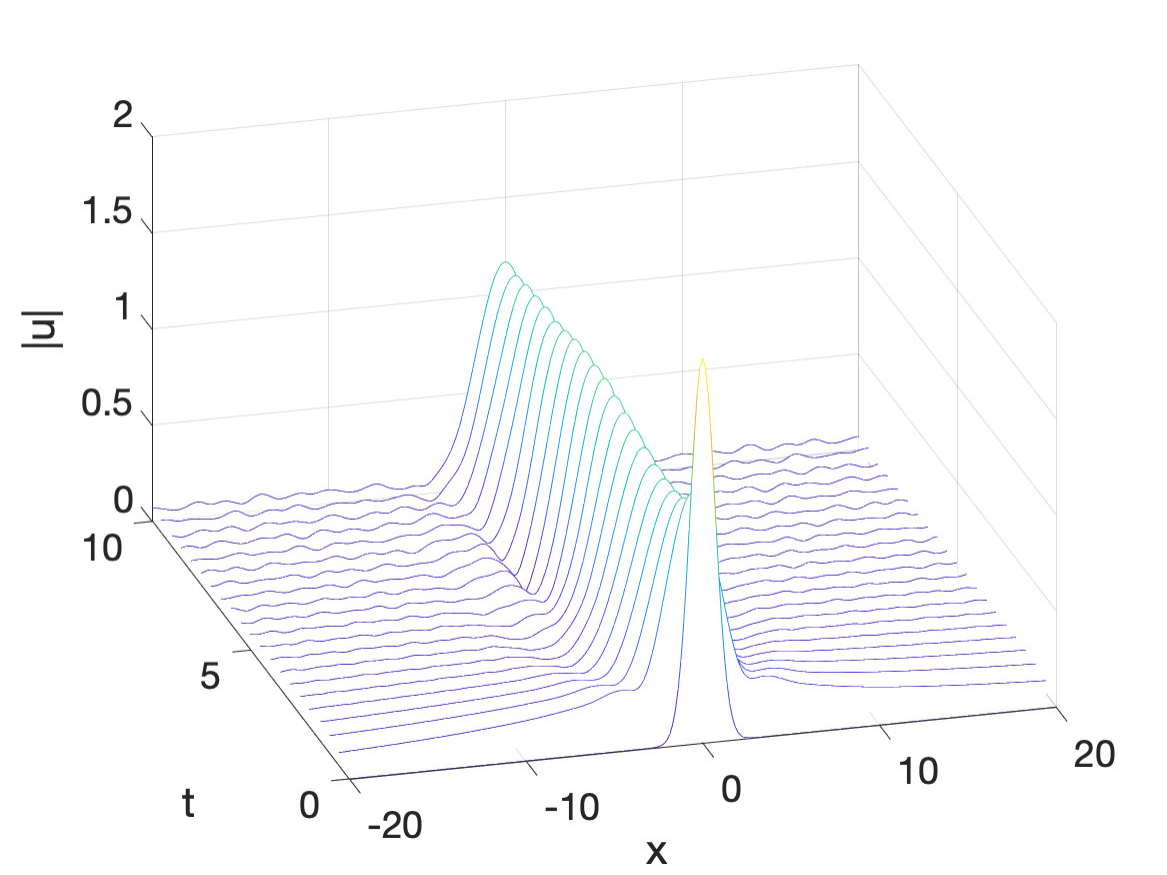}
\includegraphics[width=0.31\textwidth,height=0.25\textwidth]{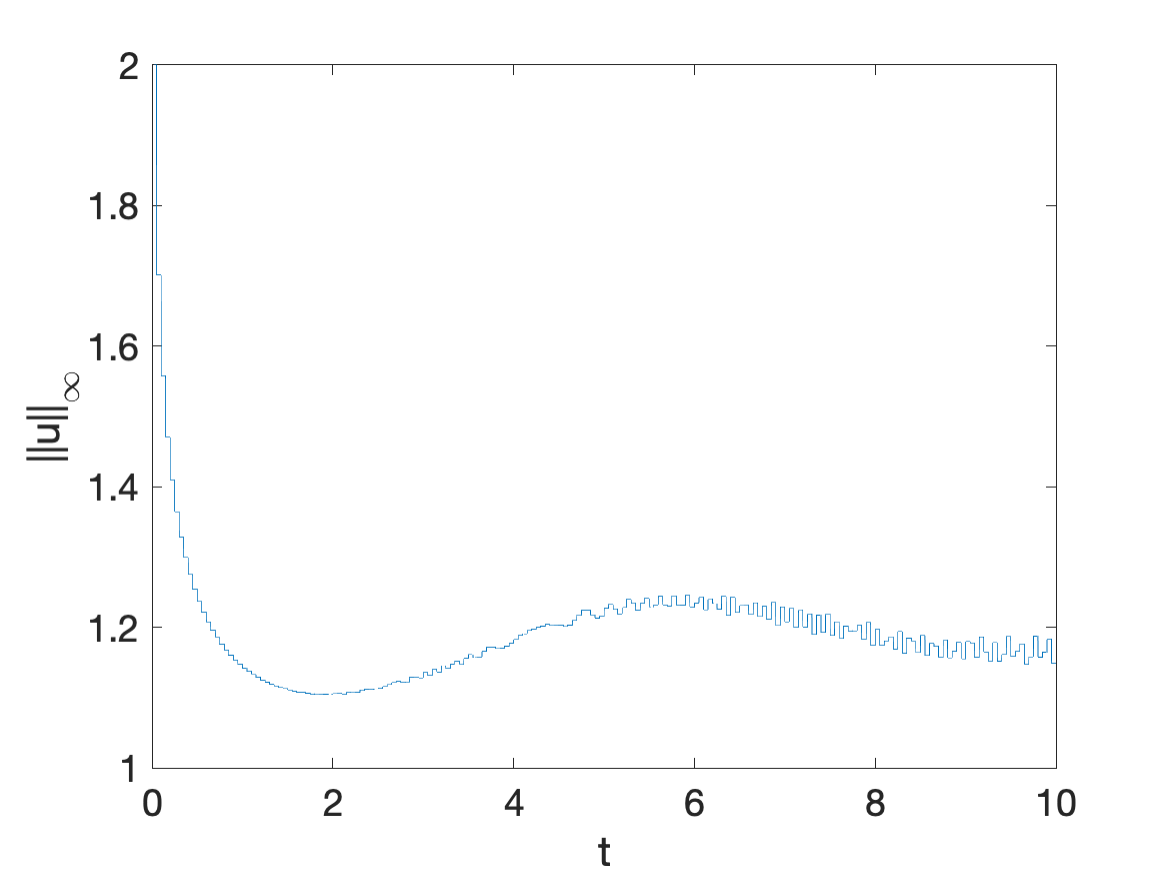}
\includegraphics[width=0.31\textwidth,height=0.25\textwidth]{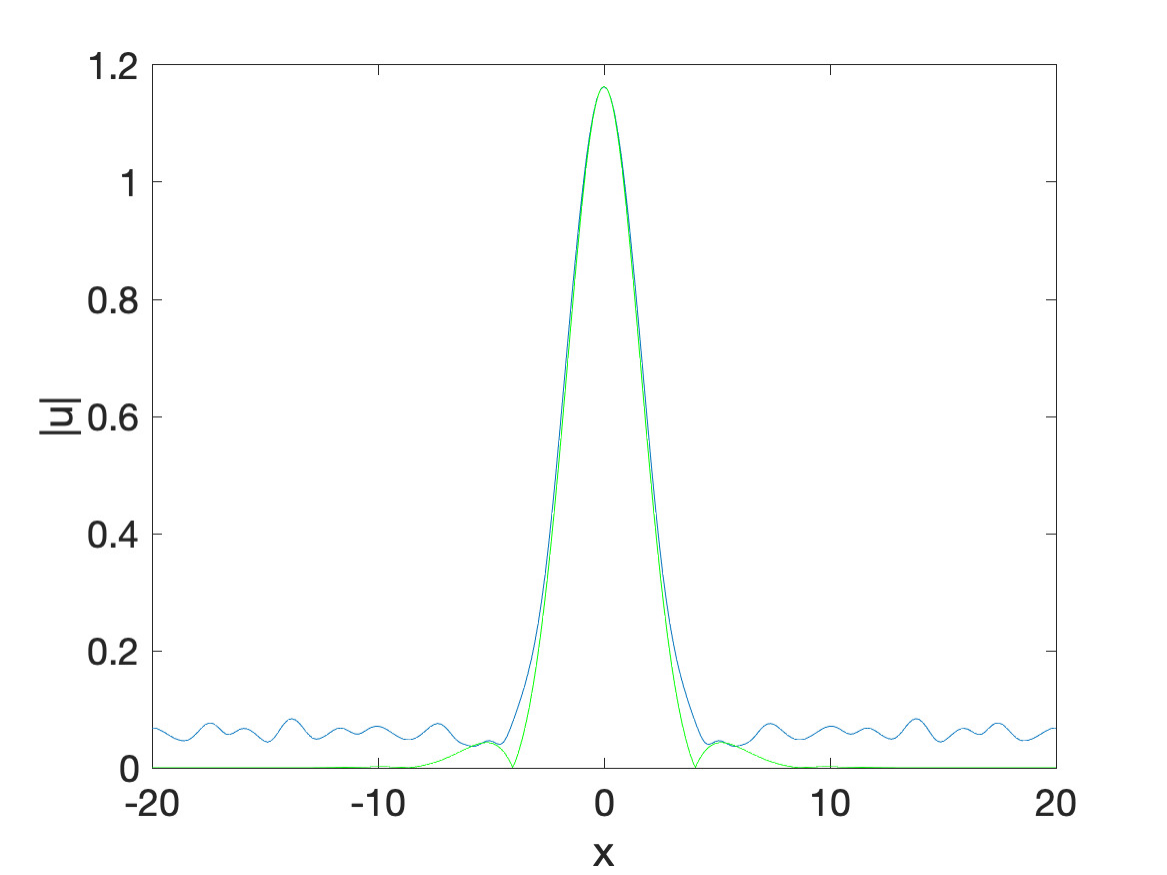}
\caption{\footnotesize Left: time evolution of the solution to \eqref{E:explicit1}, $\alpha=2$, $a=0$, with $u_0=2\,e^{-x^2}$.
Middle: time dependence of the $L^{\infty}$~norm of the solution. Right: solution at the final computational time (blue) and a fitted rescaled ground state (green).}
 \label{binlsA0alpha22gauss}
\end{figure}

{\it $\bullet$ Gaussian, $a=0$.} Fixing $a=0$ and $\alpha=2$ in \eqref{E:explicit1}, we show its solution with the initial condition \eqref{E:gaussian} and $A=2$
in Fig.~\ref{binlsA0alpha22gauss}. The initial maximal peak height drops down from $2$ to about $1.2$ and an asymptotic profile appears to emerge, although on a non-negligible radiation background. 
The $L^{\infty}$~norm of the solution is shown in the middle of Fig.~\ref{binlsA0alpha22gauss}.

The norm appears to oscillate asymptotically approaching some final state, though the 
radiation in the background becomes visible after some time (see left and middle plots after about $t=5$). 
The interpretation of the asymptotic final state as {\it a ground 
state} is supported by the right plot of the same figure, where we 
show a fit to a rescaled soliton according to \eqref{Qscaling}. This 
is one of the biggest advantages of the pure quartic case, having the 
scaling. This also confirms the soliton resolution, where localized smooth generic data splits into a rescaled and shifted soliton(s) and radiation. 

\begin{figure}[htb!]
\includegraphics[width=0.36\textwidth,height=0.26\textwidth]{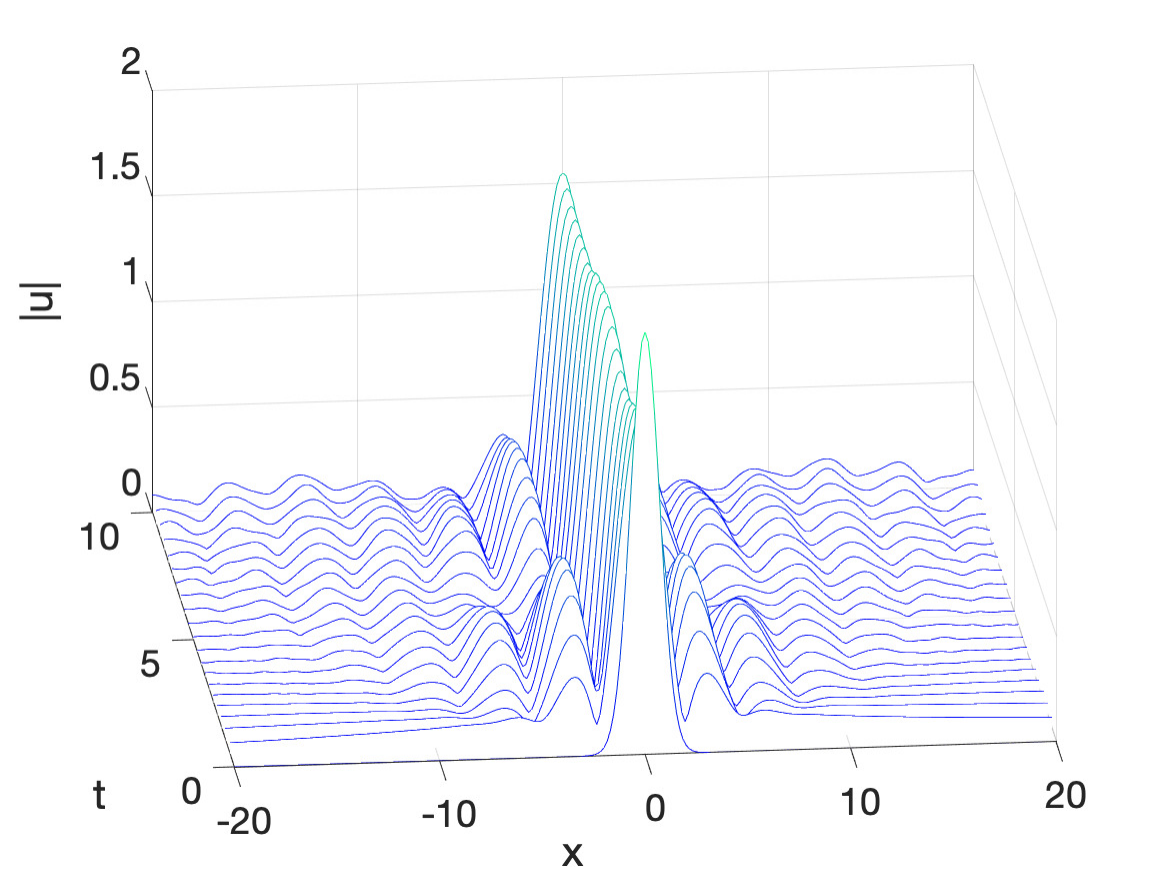}
\includegraphics[width=0.31\textwidth,height=0.25\textwidth]{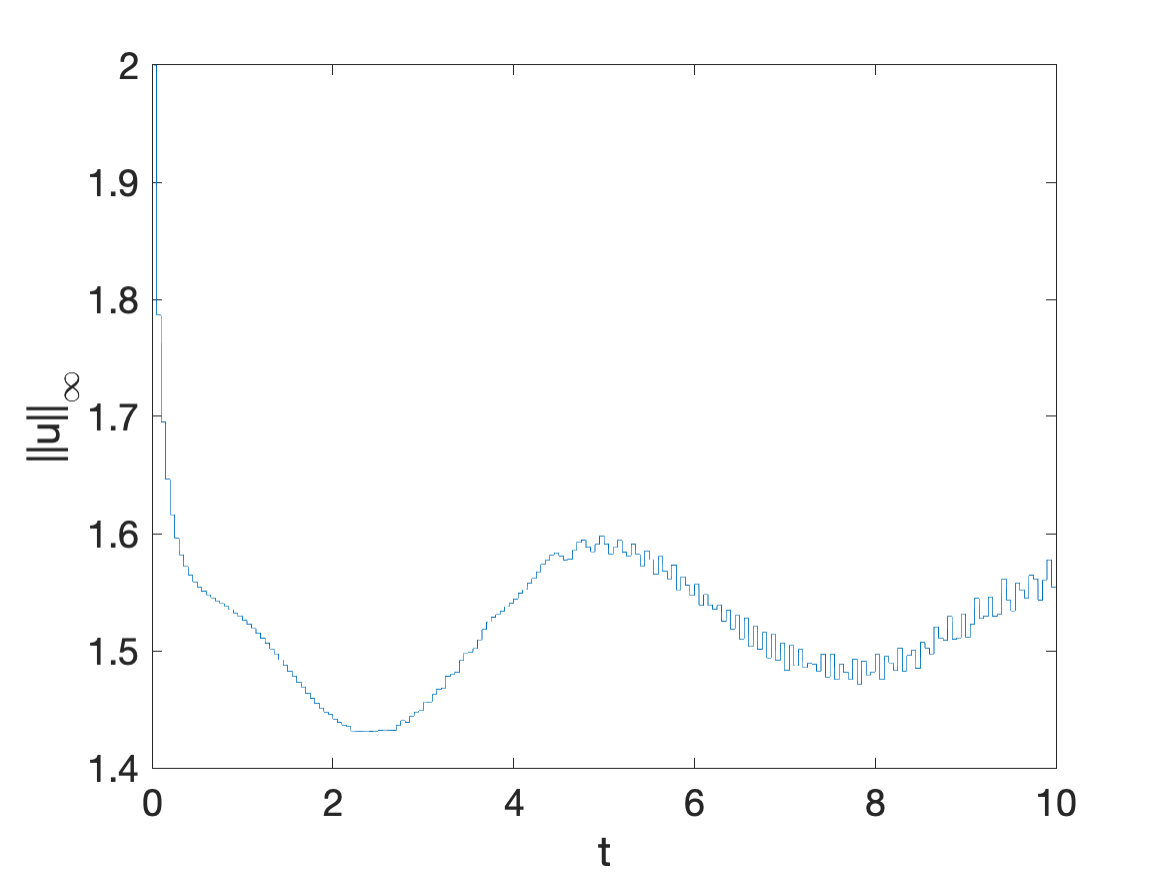}
\includegraphics[width=0.31\textwidth,height=0.25\textwidth]{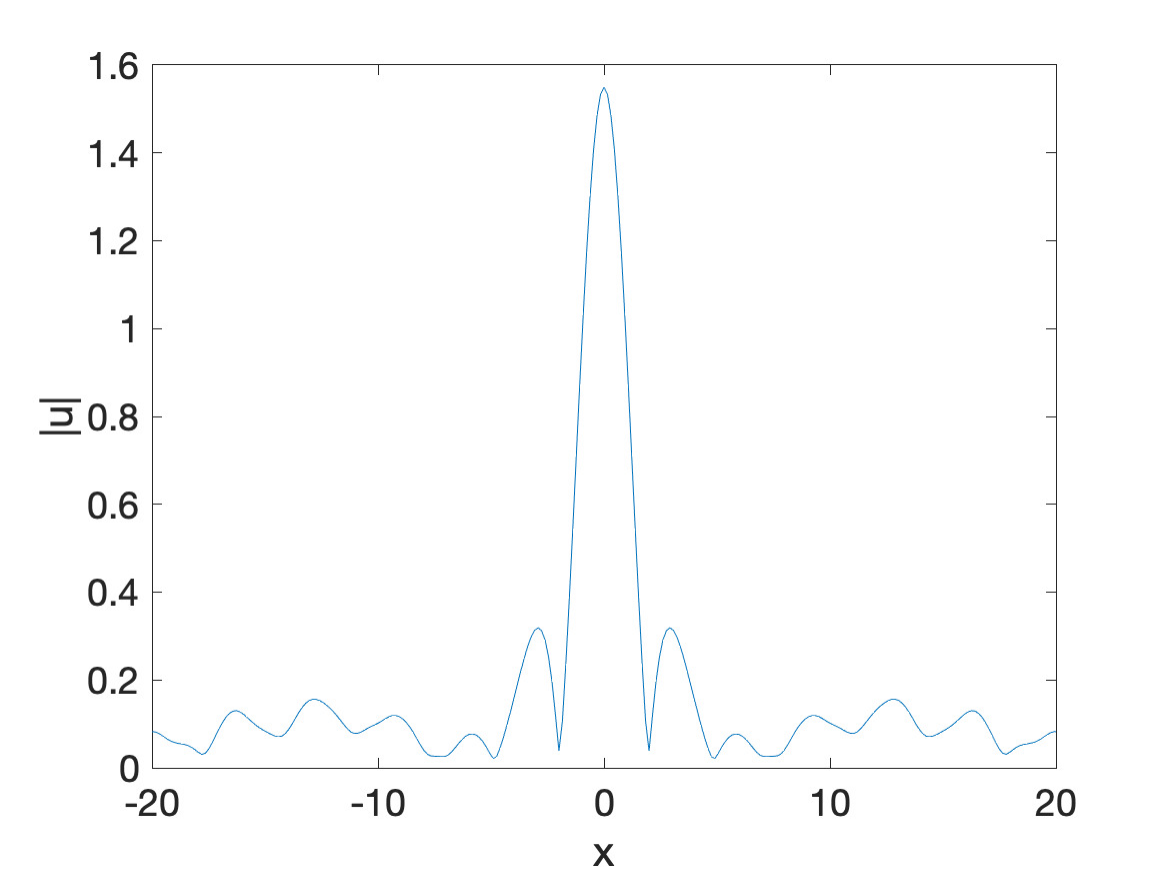}
\caption{{\footnotesize Left: time evolution of solution to \eqref{E:explicit1}, $\alpha=2$, $a=1$, with $u_0=2\,e^{-x^2}$.
Middle: Time dependence of $L^{\infty}$~norm.
Right: solution at the final computational time.
}} 
\label{binlsA1alpha22gauss}
\end{figure}

{\it $\bullet$ Gaussian, $a=1$.}
Considering the same Gaussian initial data \eqref{E:gaussian} for the equation \eqref{E:explicit1} but with $a=1$, we get the solution shown on the left of 
Fig.~\ref{binlsA1alpha22gauss}. 
The height of the main peak drops down from $2$ to around $1.5$ and then oscillates to the asymptotic final state while emitting non-trivial amount of radiation.
This interpretation is in accordance with the $L^{\infty}$~norm of the solution on the left of Fig.~\ref{binlsA1alpha22gauss}. 
The norm appears to oscillate around the height of the final state, but since the computation is done on a large torus, the radiation is clearly visible for large times (see middle plot after about $t=4.5$). The solution at the final computational time also suggests that it can be considered as a ground state, though no scaling as in \eqref{Qscaling} is available, and thus, we do not perform any fitting.  

{\it $\bullet$ Gaussian, $a=-1$.}
We next consider the opposite sign of the parameter $a$, namely, $a=-1$, where the ground state $Q^{(-1)}$ has a larger maximum in amplitude than in the cases $a=0$ or $a=1$, as can be seen on the right of Fig.~\ref{QA}. The solution to the bi-harmonic NLS equation \eqref{E:explicit1} with the same Gaussian initial condition \eqref{E:gaussian} and $A=2$,  
appears to radiate away (or scatter), no stable structures emerge. 
\begin{figure}[htb!]
\includegraphics[width=0.36\textwidth,height=0.26\textwidth]{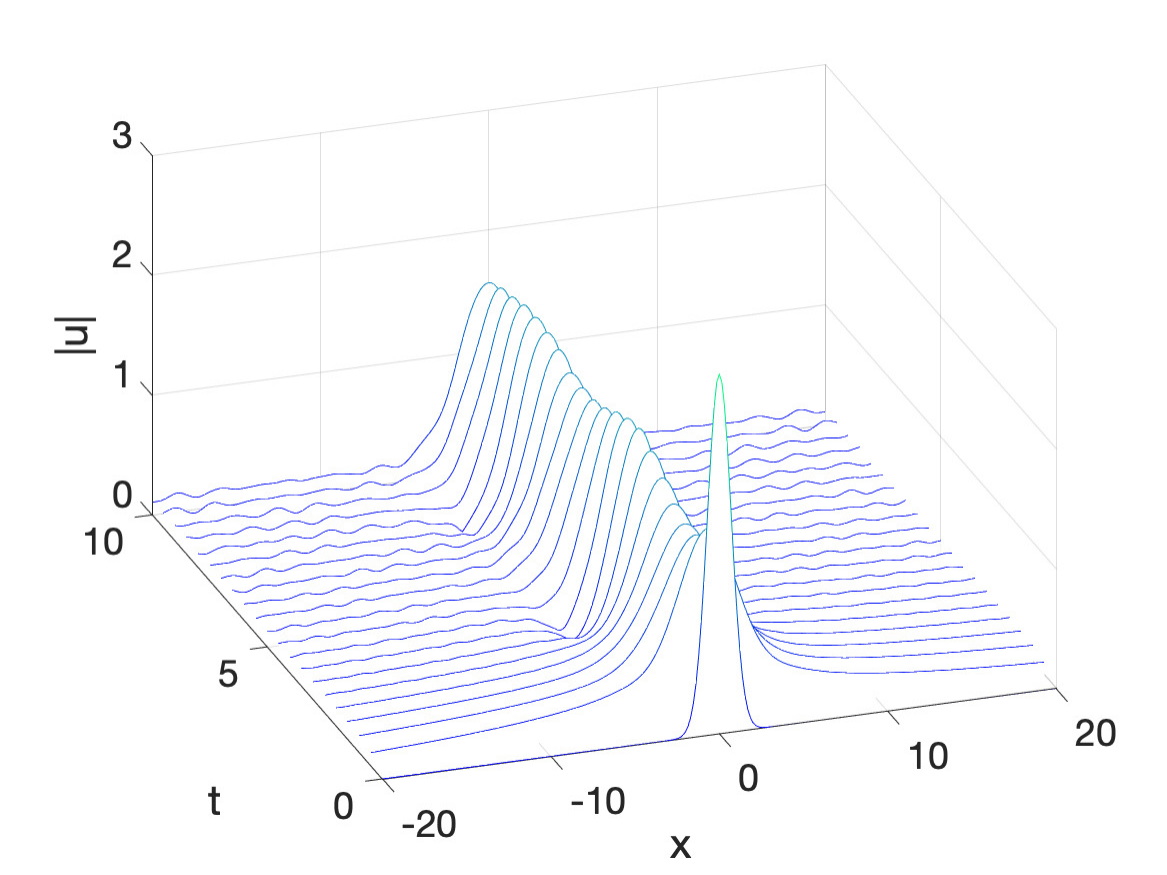}
  \includegraphics[width=0.31\textwidth,height=0.25\textwidth]{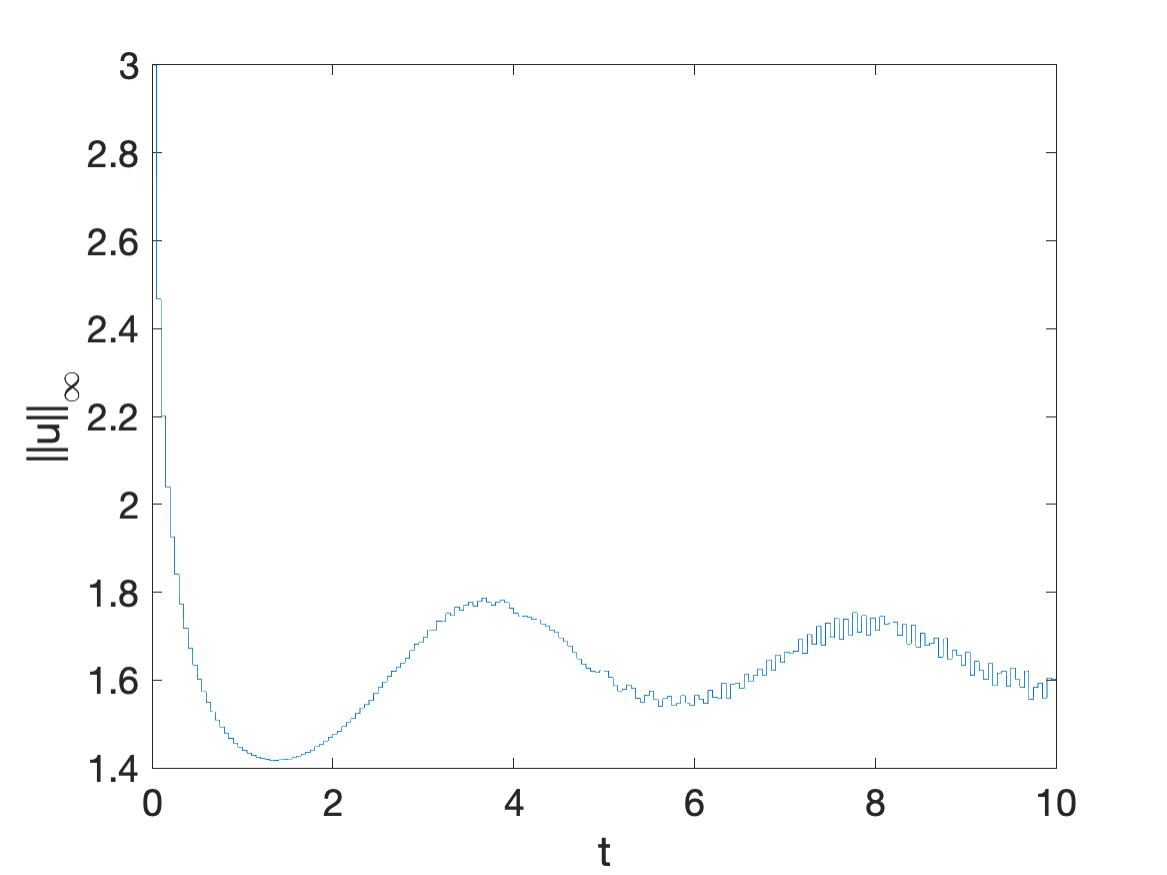}
  \includegraphics[width=0.31\textwidth,height=0.25\textwidth]{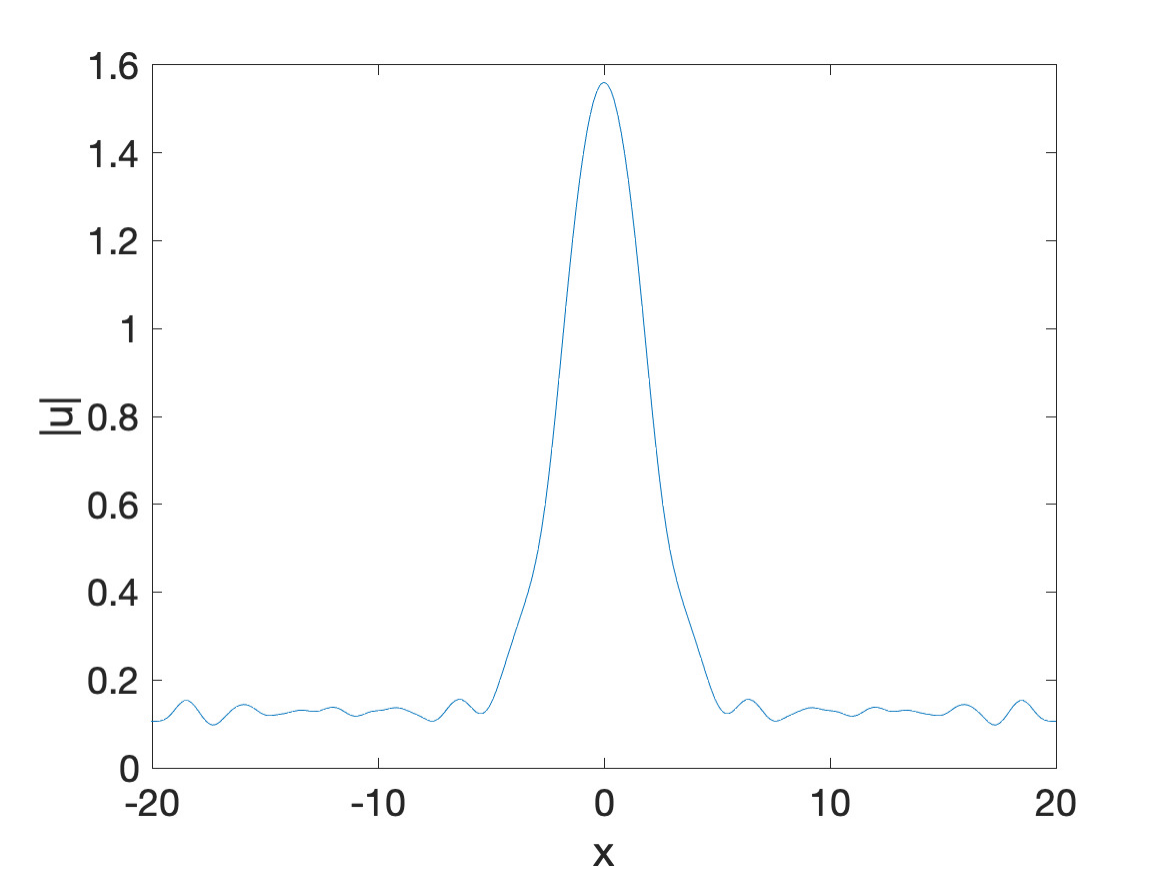}
\caption{{\footnotesize Left: time evolution of solution to \eqref{E:explicit1}, $\alpha=2$, $a=-1$, with $u_0=3\,e^{-x^2}$.
Middle: Time dependence of $L^{\infty}$~norm. 
Right: solution at the final computational time.}} 
 \label{binlsAm1alpha23gauss}
\end{figure}
However, if we consider the initial condition \eqref{E:gaussian} with a larger amplitude, for example, $A=3$, 
an asymptotic final 
state plus radiation is observed as in previous examples, see Fig.\ref{binlsAm1alpha23gauss}. 
The solution approaching an asymptotic final state
is confirmed by the $L^{\infty}$ norm of the solution in the middle of Fig.~\ref{binlsAm1alpha23gauss}, where as before, the solution is oscillatory approaching the asymptotic profile; we plot the solution at the final computational time on the right. Since there is no scaling in this case either, there is no matching to the ground state provided.

We conjecture that the scattering, or dispersion of the solution down to zero, happens if the mass of the initial data is below some threshold that is less than $\min \{Q^{(0)}, Q^{(a)} \}$ for a given power $\alpha < 8$ (and $a < \sqrt{b}$), which would be interesting to investigate further. Above this threshold, localized smooth solutions confirm the soliton resolution.


\section{Dynamics of solutions in the critical case}\label{S:critical}

When $\alpha=8$, in the critical case of \eqref{biNLS} in 1d, 
one expects to observe blow-up behavior for some 
solutions, see \cite{Fibich2002}. In dimensions two and higher the existence of finite time blow-up  
was proved in \cite{BL2017} for certain cases of \eqref{biNLS}. 
In the critical case, once the mass threshold $M[Q]$ is reached by {\it sufficiently localized} initial data $u_0$ (i.e., $M[u_0] \geq M[Q]$), 
the solution is expected to blow up in some analogous fashion with the standard NLS case, see \cite{Fibich2002}, \cite{BFM2010}. However, a more careful investigation shows that in the non-pure case, this is not necessarily the case. In fact, we show for various types of data, that there is a gap above the mass of ground states, above which the solutions do not blow-up.

In this section we first look at a possible threshold for blow-up, 
and then investigate the blow-up behavior as a self-similar solution, studying its rate and profile.

\subsection{Behavior above the mass of ground states.}
\begin{figure}[htb!]
\includegraphics[width=0.45\hsize, height=0.33\hsize]{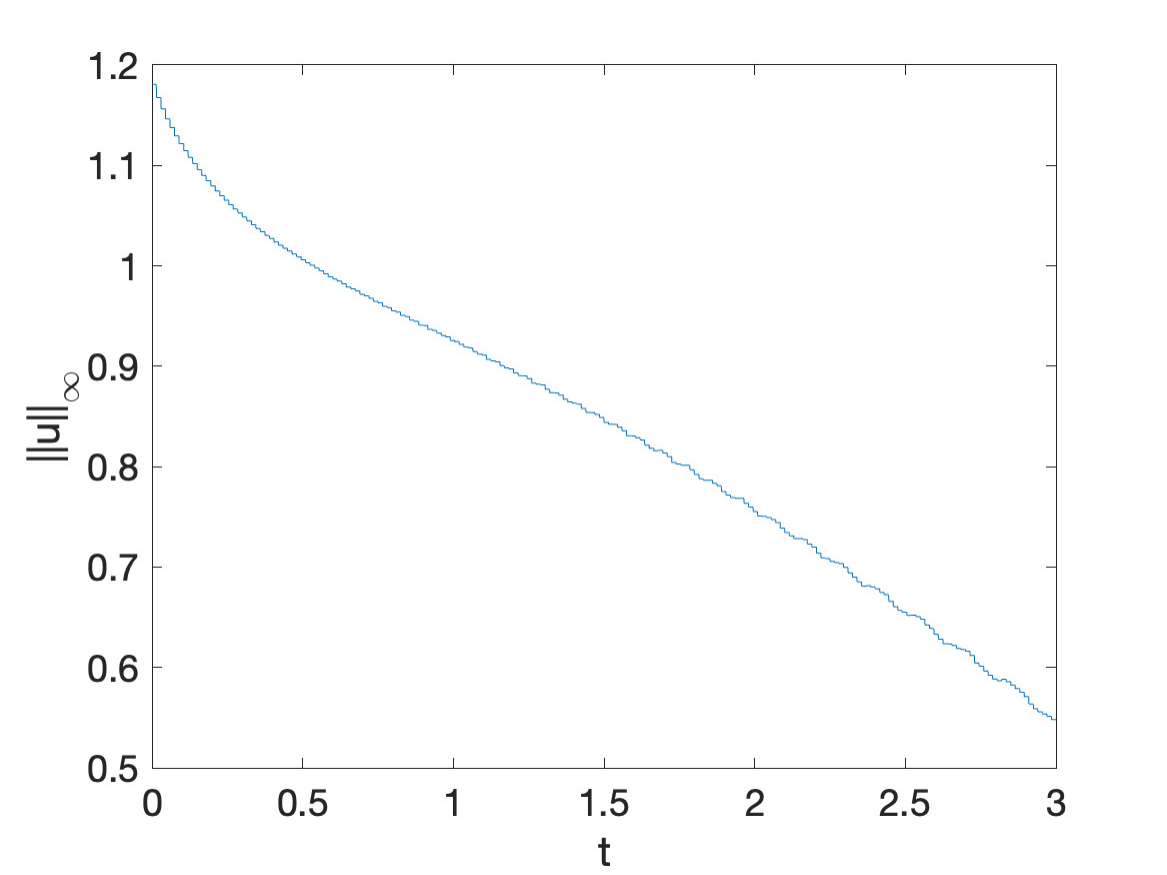}
\includegraphics[width=0.45\hsize, height=0.33\hsize]{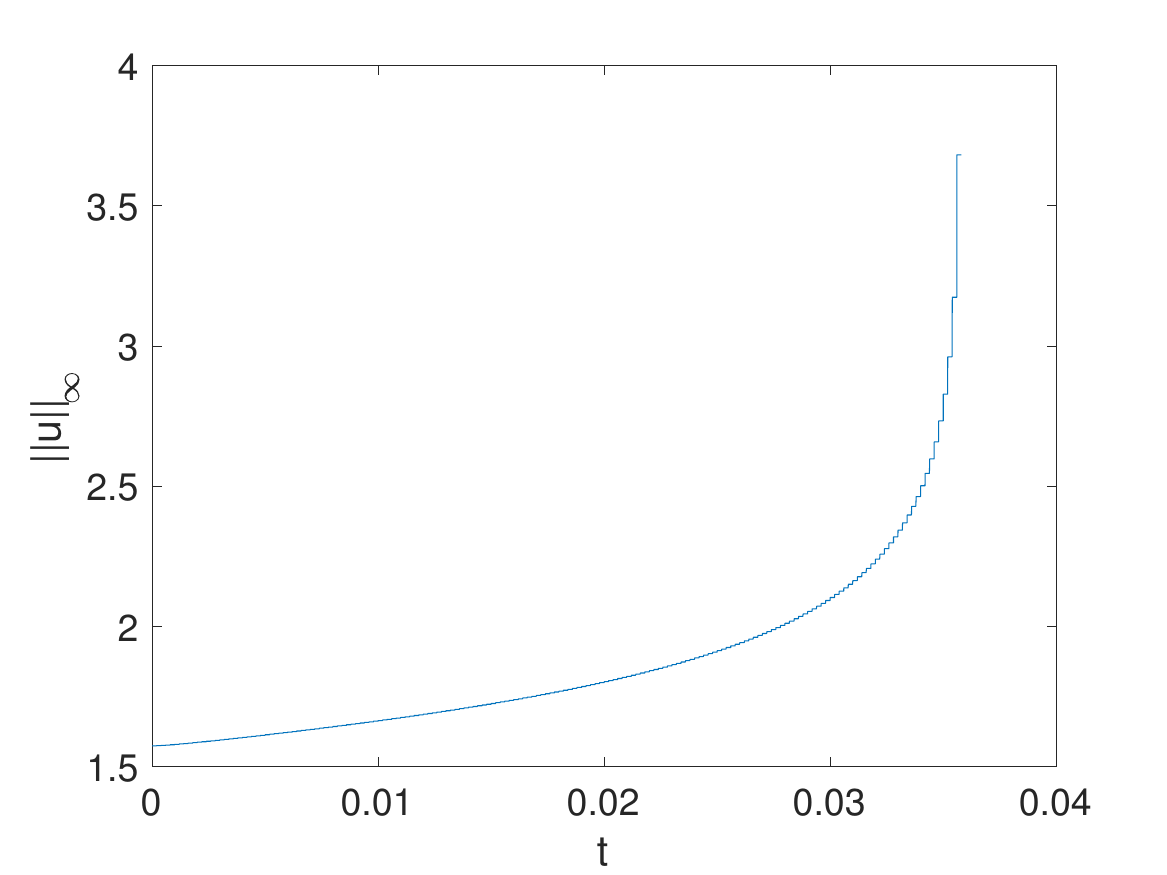}
\caption{\footnotesize {Time dependence of the $L^\infty$ norm of the solution $u(t)$ for $b=4$, $\alpha=8$, $a=1$ with $u_0=AQ^{(1)}$, $A=0.9$ (left) and $A=1.2Q^{(1)}$ (right).}}
\label{F:alpha8-blowupQ}
\end{figure}

We look at different solutions in the critical case, 
recalling from Section \ref{S:jumping} (Fig.~\ref{F:alpha8branch}, top row, and Fig.~\ref{F:alpha8-blowupQ}, right plot) 
that perturbations of ground states $Q^{(a)}$ have a somewhat different behavior than in the standard NLS case:  
for very small perturbations of $A Q^{(a)}$, $a \neq 0$, with $1<A < 1+\epsilon$, the solutions  
do not blow up, but asymptotically (in oscillatory manner) approach another stable state. 
For the first several examples of ground state perturbations we fix $b=2$, later we vary this value.  

We first consider cases $a \leq 0$ and then $a>0$, since for small $\epsilon>0$ the mass $M[Q^{(0+\epsilon)}] < M[Q^{(0)}] < M[Q^{(0-\epsilon)}]$ as shown in Figure \ref{F:alpha8}.
\begin{figure}[!htb]
 \includegraphics[width=0.45\textwidth,height=0.31\textwidth]
 {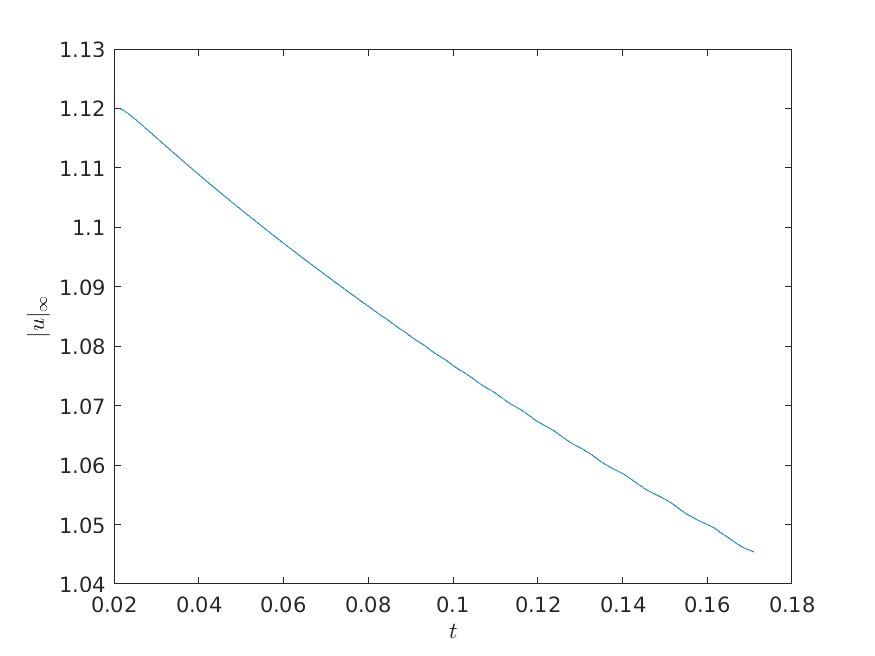}
 \includegraphics[width=0.45\textwidth,height=0.31\textwidth]{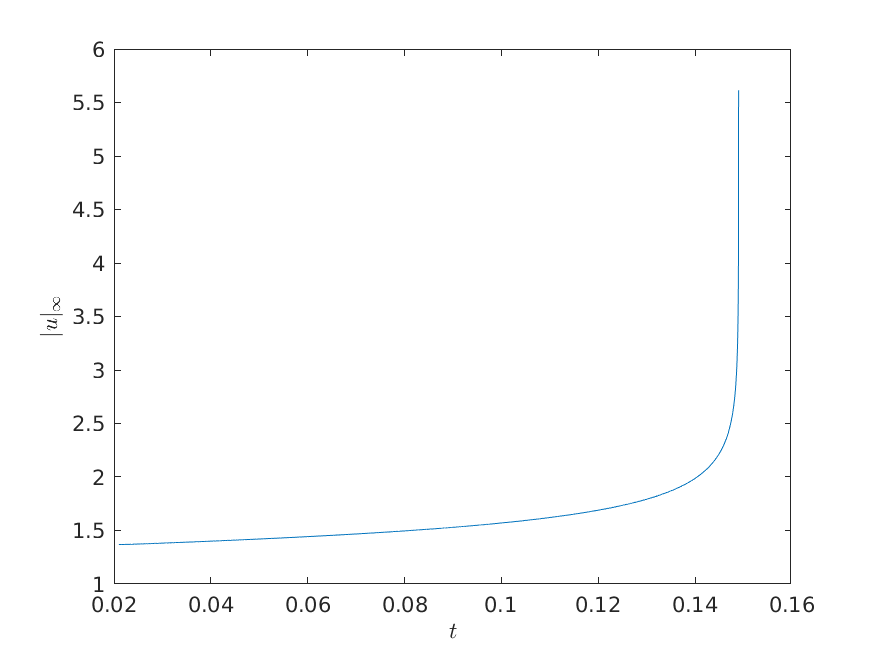}
\caption{\footnotesize {
Time dependence of the $L^\infty$ norm of solutions to \eqref{E:explicit1}, $\alpha=8$, $b=2$, $a=0$, with $u_0 = A\, Q^{(0)}$, $A=0.9$ (left) and $A= 1.1$ (right).}}
\label{F:Q0critical}
\end{figure}

$\scriptstyle\blacklozenge$ \underline{Case {$a=0$.}}
We consider 
$u_0(x)=A\, Q^{(0)}(x)$, with $A \approx 1$.
The time evolution of the $L^\infty$ norm of solutions with $A=0.9$ or $1.1$ is shown in Fig.~\ref{F:Q0critical} on the left and right, correspondingly. 
In the first case the solution is dispersed, whereas in the second case it blows up in finite time. We also recall that the energy $E[Q^{(0)}] = 0$ and the value of the ground state mass does not depend on $b$ in this pure quartic dispersion case, as the solution can be rescaled. 

Furthermore, since the scaling invariance holds in the $a=0$ case, we are able to fit the final computational state with the (numerical) rescaled ground state profile of $Q^{(0)}$, 
see details in Section \ref{S:profile-critical}.

$\scriptstyle\blacklozenge$ \underline{Case $a<0$.} We fix $a=-1$ and observe that the perturbation of $u_0 = A\, Q^{(-1)}$ with $A = 0.95$ produces a bounded (oscillatory) behavior and with $A=1.1$ the solution blows up in finite time. 

Next, we consider $a=-2$ and compute the corresponding ground state solution $Q^{(-2)}$ from \eqref{E:1dGS}, which has no oscillations in this case (since $a<-\sqrt 2$). 
We consider perturbations $u_0 = A Q^{(a)}$ with $A=1.1$. This initial data produces a blow up in finite time, see Fig.~\ref{F:Soliton11_A-2} (we discuss this case in more details below in terms of the blow-up rate and profile). 

\smallskip

$\scriptstyle\blacklozenge$ \underline{Case $a>0$.} 
We fix $a=1$ and consider $u_0=A\, Q^{(1)}$, recalling that for $A=0.99$ the solution dispersed and for $A=1.01$ it  asymptotically (in oscillatory manner) approached some stable state, see Fig.~\ref{F:alpha8branch}. This behavior was observed for small perturbations  
with $1<A \lesssim 1.1$. Only when $A \gtrsim 1.2$ we observe that the perturbations 
blow up in finite time, see Fig.\ref{F:alpha8-blowupQ}. 
We note that $M[Q^{(1)}]<M[Q^{(0)}]$ (see Table \ref{T:1}), and when 
$A=1.10055$, we get $M[Q^{(0)}] = M[A\, Q^{(1)}]$. Thus, in this 
case, we observe that 
\begin{itemize}
\item 
scattering happens below the mass of $Q^{(1)}$ (and in this case, $M[Q^{(1)}]<M[Q^{(0)}]$), 
\item
there is a gap between $M[Q^{(1)}]$ and $M[Q^{(0)}]$, where solutions with perturbed $A\, Q^{(1)}$ data do not scatter but instead approach, in oscillatory manner, some final state,
\item 
when the mass of the initial data gets above $M[Q^{(0)}]$, solutions blow up in finite time. 
\end{itemize}

Fixing $a = 0.5$ and considering the ground state $Q^{(0.5)}$ of \eqref{E:1dGS} and its perturbations $u_0 = A\, Q^{(0.5)}$, 
we observe that $A = 1.05$ produces a bounded (oscillatory) behavior and with $A=1.1$ the solution blows up in finite time. 
Noting again that $M[1.1 Q^{(0.5)}] = 3.23 > M[Q^{(0)}]$ and $M[1.05 Q^{(0.5)}] = 2.95 < M[Q^{(0)}]$ gives a positive confirmation to the statement that 
the behavior of solutions is determined by the size compared to either the ground state $M[Q^{(1)}]$ or the scaling-invariant ground state $Q^{(0)}$ and its mass.  

Fixing $a=1.35$, we first note that $M[Q^{(1.35)}] \approx 3.41 > 
M[Q^{(0)}]$, and recalling the  left plot of Fig.~\ref{F:alpha8}, which indicates that this ground state is on the unstable branch, we investigate the perturbations of  $u_0 = A\, Q^{(1.35)}$.  
For $A = 1.2$ the solution produces a bounded (oscillatory) behavior and with $A=1.5$ the solution blows up in finite time. 
While $M[1.2Q^{(1.35)}]$ is greater than that of the scaling invariant ground state mass, $M[Q^{(0)}]$, we point out that the observed behavior could be an indication that we are dealing with the unstable branch of ground states. At larger magnitudes of $A$ we do observe blow-up in finite time.

\subsubsection{Behavior of other types of data above the threshold} 
We consider different types of initial data to investigate further the behavior of solutions above the corresponding ground state's mass in the critical case. 

$\scriptstyle\blacklozenge$ We take {\it super-Gaussian} initial data, compute its mass, and investigate its evolution, 
$$
u_0 = A\, e^{-x^4}, \qquad M[u_0] = 2^{\frac34} \Gamma(\tfrac54) A^2.
$$  
We note that the mass of super-Gaussian data is smaller than the mass of the pure quartic ground state $Q^{(0)}$, that is, $M[u_0] < M[Q^{(0)}]$, when
(approximately) $A < 1.4$. 

\begin{figure}[!htb]
\includegraphics[width=0.49\hsize, height=0.32\hsize]{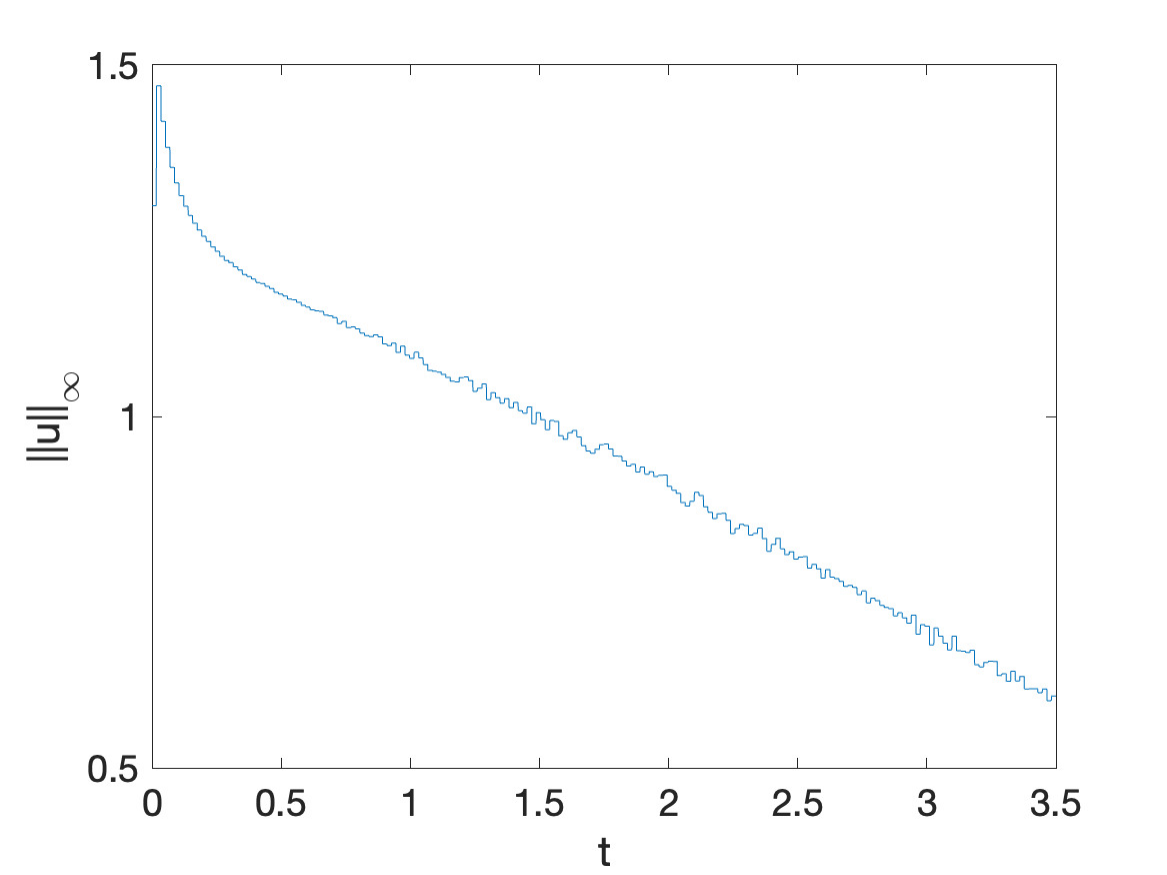}
\includegraphics[width=0.49\hsize, height=0.32\hsize]{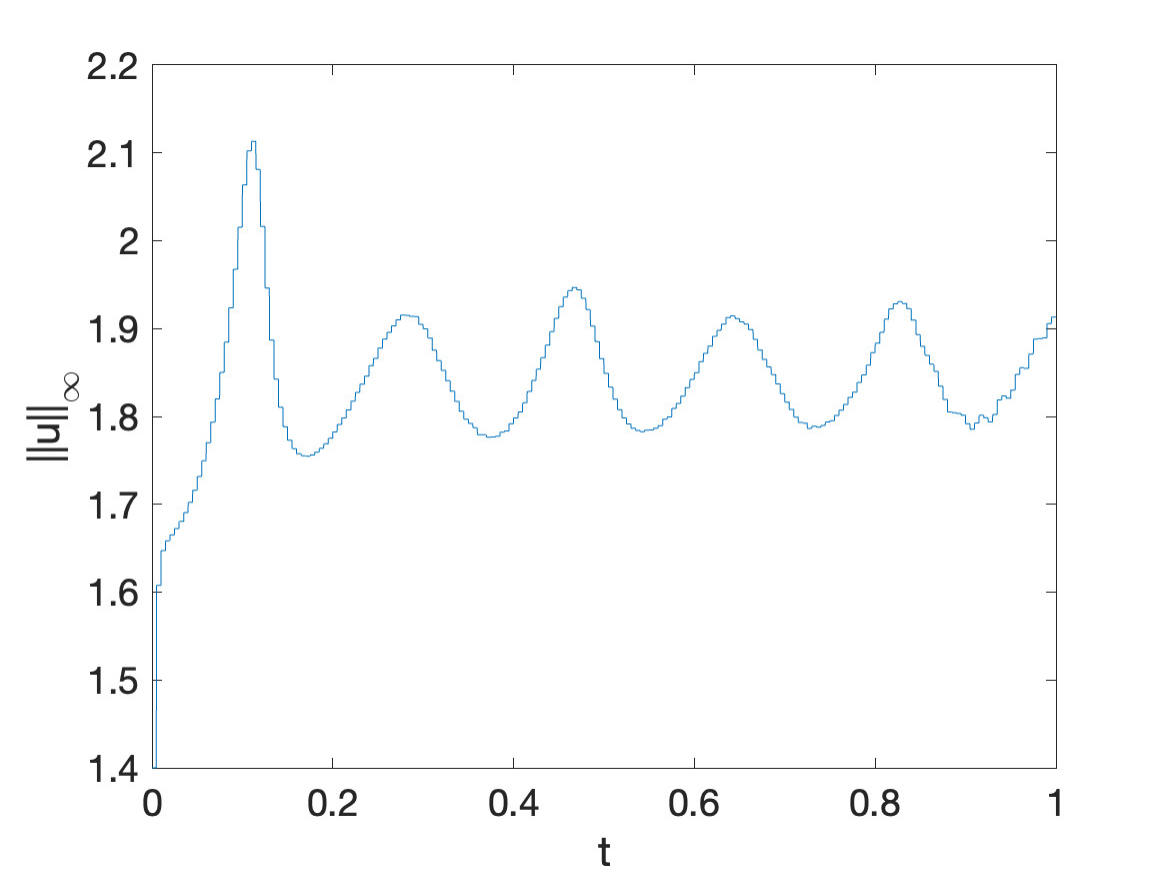}
\caption{\small {Left column: time dependence of the $L^\infty$ norm of solution to \eqref{E:explicit1}, $\alpha=8$, $a=1$, with 
$u_0 = A e^{-x^4}$, $A=1.3$ (left) and $A=1.4$ (right).}}
\label{F:superG}
\end{figure}

In Fig.~\ref{F:superG} we show the $L^\infty$ norm of solutions to \eqref{biNLS}, $\alpha=8$, $a=1$, 
for the super-Gaussian initial data with $A=1.3$ on the left and with $A=1.4$ on the right. 
For values of $A$ larger than $1.4$, we observe finite time blow-up. 
Thus, the blow-up occurs for the initial data with the mass above that of the scaling-invariant ground state $M[Q^{(0)}]$ (though the equation has the second order dispersion term with $a=1$), then there is a gap where the solution does not seem to neither blow-up nor scatter down to zero, but instead approach some final state, and below certain value the solutions scatters to zero.

It seems to be plausible to conjecture that if the mass of the non scaling-invariant ground state is less than the mass of $Q^{(0)}$, i.e., if $M[Q^{(a)}] < M[Q^{(0)}]$, then the solution scatters to zero in that case. 
There seems to be a small gap where solutions do not disperse to zero, nor blow-up, but rather approach another non-zero state. 
Eventually, the perturbations with larger $A$ blow-up in finite time, starting from a certain threshold above that mass of $Q^{(0)}$ (or above $M[Q^{(a)}]$, if that value is larger).  
\smallskip

\begin{figure}[!htb]
\includegraphics[width=0.49\hsize, height=0.32\hsize]{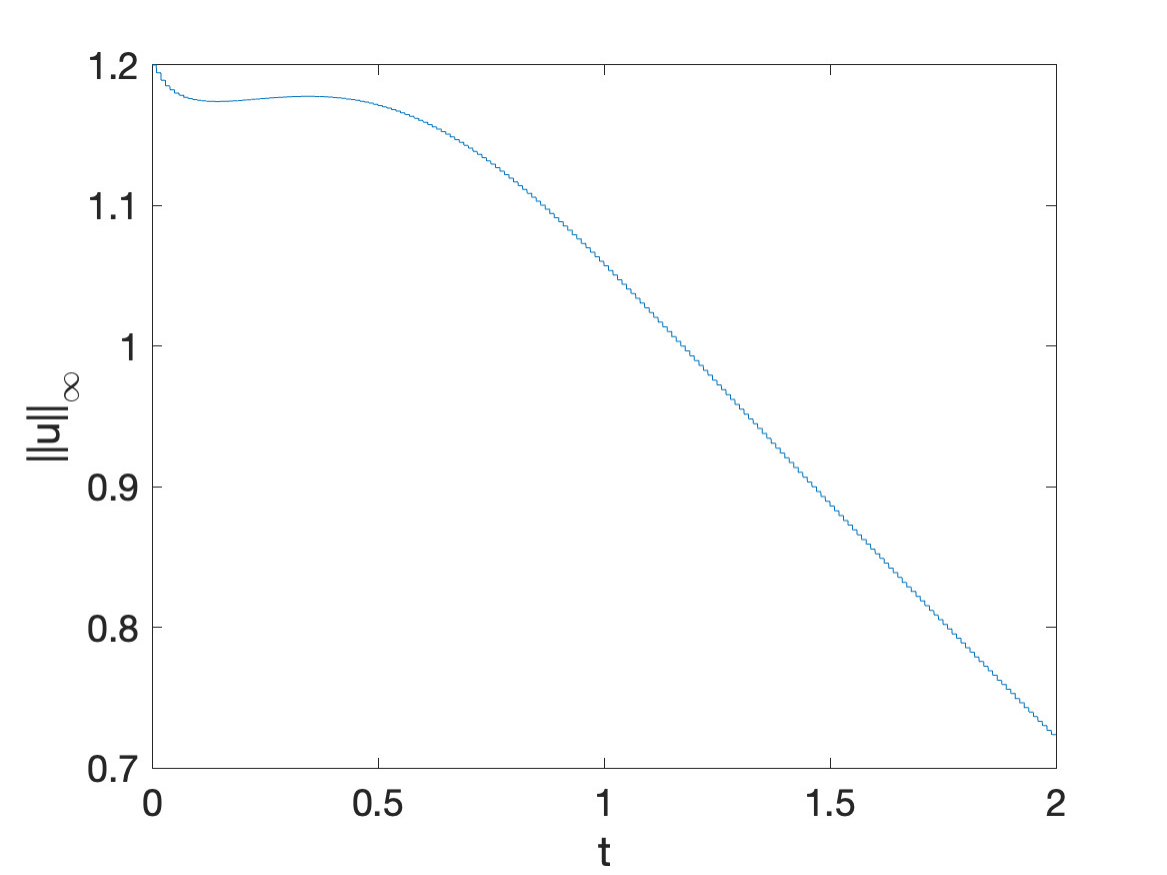}
\includegraphics[width=0.49\hsize, height=0.32\hsize]{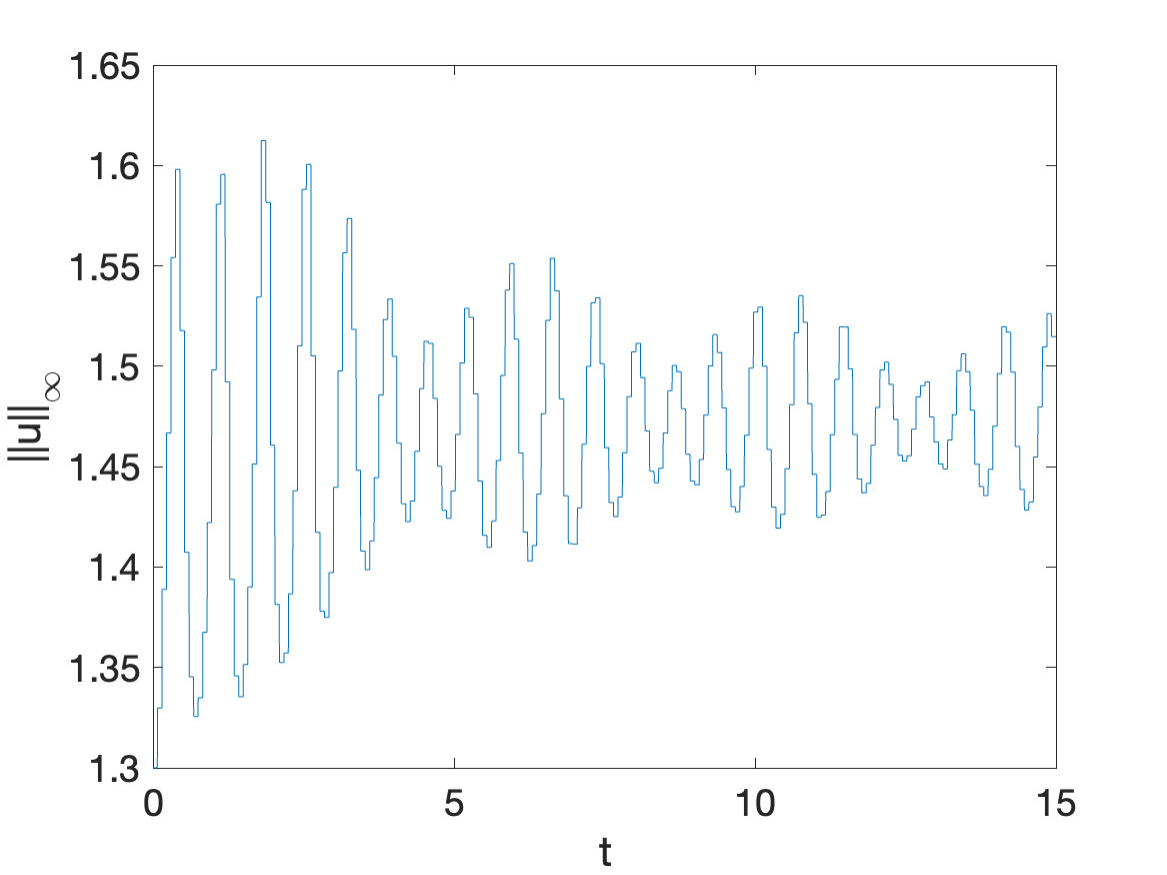}
\caption{\small {Left column: time dependence of the $L^\infty$ norm of solution to \eqref{E:explicit1}, $\alpha=8$, $a=1$, with 
$u_0 = A \, \sech x$, $A=1.2$ (left) and $A=1.3$ (right).}}
\label{F:sech}
\end{figure}

$\scriptstyle\blacklozenge$
Another example we look at is of slower decay than Gaussian or super-Gaussian, namely, 
$$
u_0=A\, \sech\, x, \quad 
M[u_0] = 2 A^2.
$$ 
The condition $M[u_0]< M[Q^{(0)}]$ for such data holds when $A < 1.23$. 
We show the dynamics of solutions for such data in Fig.~\ref{F:sech} with $A=1.2$ on the left, which disperses, and with $A=1.3$ on the right, which oscillates toward some final state. For values much larger than 1.3, $A \gg 1.3$, we observe blow-up.  

To summarize the results of this part, we find that in the critical case, the thresholds for scattering (down to zero) or blow-up are not easily identified for different parameters of lower dispersion $a$ or initial magnitudes and profiles (and $b$); 
the values of the mass of corresponding ground states, and especially the scaling-invariant value $M[Q^{(0)}]$ should be taken into account, as well as the values (and possibly other properties) of the ground states $Q^{(a)}$. Furthermore, the stable and unstable branches of the ground states play a significant role in the solutions behavior, possibly identifying thresholds for scattering down to zero or approaching some final states, or blow-up in finite time. 
We do show a new behavior in the critical case, which is significantly different from classical dichotomy of scattering vs. blow-up, as we observe a `gap' when solutions neither scatter nor blow-up but rather approach some final state.  
\smallskip 

We next investigate self similarity in blow up solutions.  

\subsection{Self-similar blow-up}\label{S:blowup} 
A typical approach to study the self-similar blow-up  solutions in equations with scaling  
symmetry, for example, as \eqref{scaling}, is to consider the method of dynamical rescaling (see \cite[\S 6.1.2]{SS1999}) with a time dependent factor $L(t)$. To be precise, set
\begin{equation}
\xi=\frac{x}{\lambda(t)},\quad \frac{d\tau}{dt}=\frac{1}{\lambda(t)^{4}},\quad 
u(t,x) = \frac1{\lambda(t)^{4/\alpha}} \,U (\tau,\xi),
\label{dyn1}
\end{equation}
with the factor $\lambda(t)$ to be chosen such that it tends to zero as $t\to t^{*}$ (finite blow-up time $t^{*}<\infty$), 
whereas the rescaled time variable $\tau$ tends to infinity in this limit. 
This time dependent change of variables replaces the equation \eqref{E:explicit1} with 
\begin{equation}
iU_{\tau}-i(\ln \lambda)_{\tau}(\tfrac4{\alpha}\, U+\xi U_{\xi})-U_{\xi\xi\xi\xi}
-2a\lambda^{2}U_{\xi\xi}+|U|^{\alpha}U=0.
\label{dyn2}
\end{equation}

In order to look for singular self-similar type solutions near the blow-up, the following ansatz 
$$
U(\tau,\xi) = e^{ib\tau} \mathcal R^{(\tau)}(\xi)
$$ 
is used to obtain the profile equation (we drop the superscript $\tau$ in the profile $\mathcal R^{(\tau)}$ for brevity):
\begin{equation}
-b\, \mathcal R -i(\ln \lambda)_{\tau}(\tfrac4{\alpha}\, \mathcal R + \xi \mathcal R_{\xi}) - \mathcal R_{\xi\xi\xi\xi}
-2 a \lambda^{2} \mathcal R_{\xi\xi}+|\mathcal R|^{\alpha} \mathcal R=0.
\label{E:profile-super}
\end{equation}

In the limit $\tau\to\infty$ (i.e., $t \to t^*$), the blow-up profile $\mathcal R^{(\tau)}$
is expected to become $\tau$-independent, denote it by $\mathcal R^\infty$. Setting  
$A^{\infty}:=\lim_{\tau\to\infty}(\ln \lambda)_{\tau}$ and observing that, 
since $\lambda(t) \to 0$ as $\tau \to \infty$ ($t \to t^*$), the term proportional to $a$ 
vanishes 
and we obtain the profile equation from 
\eqref{E:profile-super} as follows
\begin{equation}\label{E:Q1}
-b \, \mathcal R^{\infty}-i\,A^{\infty}(\tfrac4{\alpha} \, \mathcal R^{\infty}+\xi 
\mathcal R^{\infty}_{\xi}) - \mathcal R^{\infty}_{\xi\xi\xi\xi}
+|\mathcal R^{\infty}|^{\alpha} \mathcal R^{\infty}=0.
\end{equation}
In the critical case ($\alpha=8$) the coefficient $A^{\infty}$  becomes $0$ (recall that in the NLS equation the rate of how fast this coefficient decreases to zero determines the correction in the stable blow-up, so called `log-log' correction of the square-root rate), and the profile equation \eqref{E:Q1} reduces to the equation 
$b\, \mathcal R^{\infty} + \mathcal R^{\infty}_{\xi\xi\xi\xi}
-|\mathcal R^{\infty}|^{8} \mathcal R^{\infty}=0$,
the same as in \eqref{E:groundstate} with $a=0$, and hence, instead of $\mathcal R^\infty$ from now on, we can simply write $Q$ in the critical case as in \eqref{E:groundstate}, namely, in the one dimensional case:
\begin{equation}\label{E:Q2}
b\, Q + \partial^4_{\xi}Q -|Q|^{8} Q=0.
\end{equation} 
Solving for $Q$ as we did in Section \ref{S:GS-num}, we obtain profiles for different values of $b$. 
Observe that since the critical equation preserves the $L^2$~norm, the mass of the ground states for different values of $b$ is the same, and thus, we simply denote it by $M[Q]$. 
We obtain (computed for different values of $b$ and also listed in Table 1)
\begin{equation}\label{E:mass}
M[Q] = 2.986792978326142,
\end{equation}
which is consistent, for example, with the results in \cite[(41)]{Fibich2002}.  

\begin{figure}[!htb]
 \includegraphics[width=0.32\hsize,height=0.28\hsize]{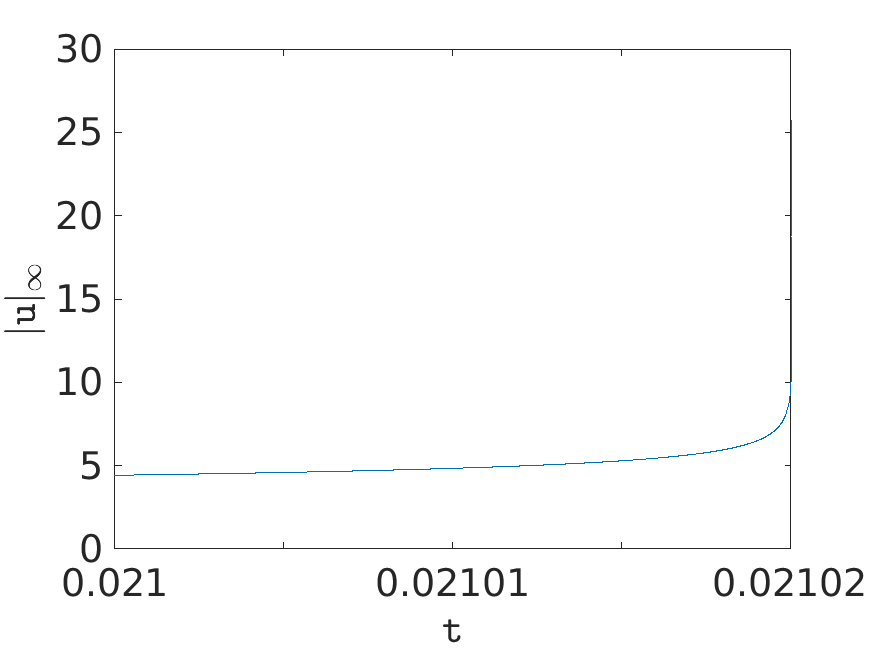}
 \includegraphics[width=0.32\hsize,height=0.28\hsize]{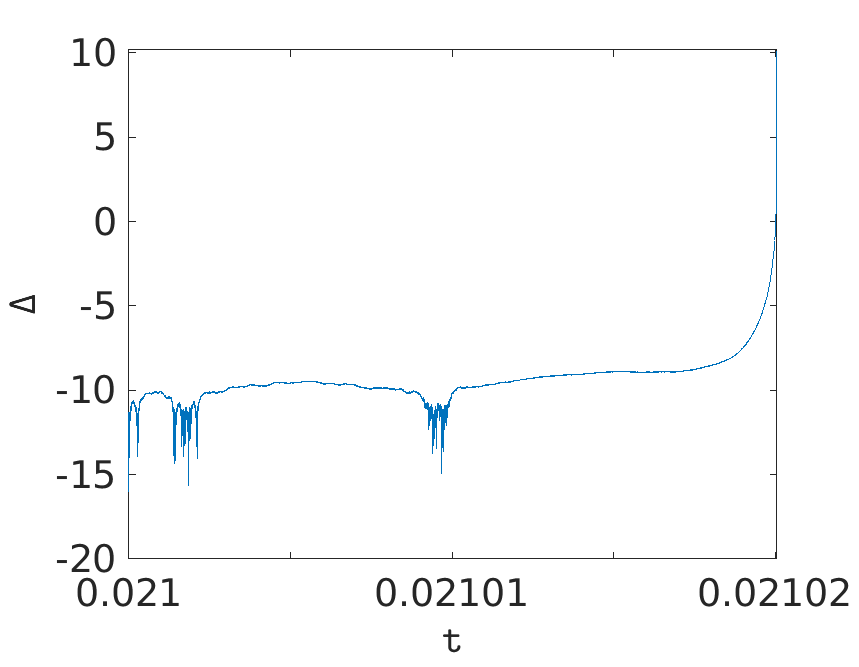}
\includegraphics[width=0.32\hsize,height=0.28\hsize]{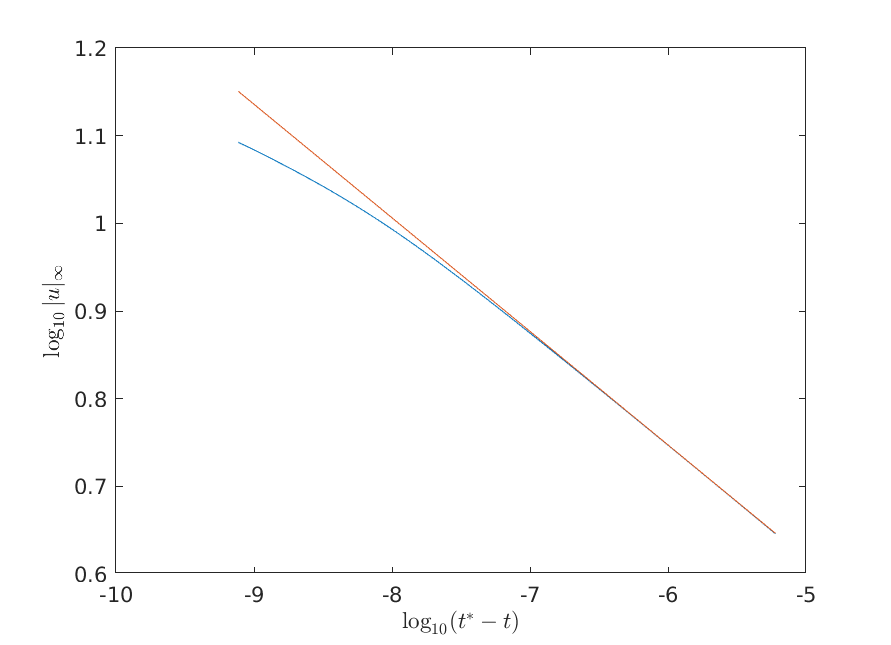}
\caption{\footnotesize {Solution to \eqref{E:explicit1}, $\alpha=8$, $a = 0$, and $u_0=1.7\, e^{-x^2}$. 
Left: time dependence of the $L^\infty$~norm. Middle: energy conservation $\Delta = \log_{10}|E - E_0|$.
Right: Blow up rate and linear fitting on a log scale, 
the blow-up occurs at $t_* = 0.0210375$, 
the fitting on a log scale up to the time step $1.5 \times 10^{-10}$ 
is too delicate to give a conclusive value for the rate parameter.  
}}
\label{F:GaussRate}
\end{figure}

To investigate self-similar blow-up numerically, we 
fix $a=0$ in \eqref{E:explicit1} 
and take the Gaussian initial data 
\eqref{E:gaussian}, $u_0 = A\, e^{-x^2}$ with $A=1.7$. 
The reason for this choice of the amplitude $A$ is that the mass of the initial data $M[u_0] = A^2 \, \sqrt{\tfrac{\pi}{2}} $ 
is about $1.1$ times the mass of the solitary wave $M[Q^{(0)}]$ (computed in \eqref{E:mass}), and thus, blow-up is expected. 
Since the mass ratio is sufficiently close to 1, it allows us to track the numerical solution 
for some non-trivial finite time before the amplitude gets too large to be handled numerically, see Fig.~\ref{F:GaussRate}. 

A note about the blow-up time: an estimated value for the blow-up time is determined experimentally. Then the time interval up to the (numerical) blow-up time is subdivided into several intervals 
with increasingly smaller time steps as it gets closer to the blow-up. 
For the last interval the time step $h\approx 10^{-10}$, which is a reasonable implementation in a double precision setting. Such small time steps are needed because of the 
dependence of the dynamically rescaled time $\tau$ on $t$ in 
\eqref{dyn1}. The $L^\infty$~norm of the solution is shown on the left 
of  Fig.~\ref{F:GaussRate}. Fitting the $L^\infty$ norm of the 
numerical solution close to the blow-up core as discussed above, we obtain the numerical 
blow-up time $t^* = 0.2102004$.
In the middle of Fig.~\ref{F:GaussRate}, we track the conservation of 
the energy during this simulation, namely, the quantity $\delta = \log_{10}|E - E_0|$  
(the computation is terminated once the relative energy conservation drops below 
$10^{-3}$).

\subsubsection{Rate of self-similar blow-up}\label{S:rate-critical}
Using the fact that \eqref{biNLS} is locally well-posed in $H^2$, a similar argument from the classical existence theory as in the critical NLS case (see \cite[Thm 5.9]{SS1999}) gives the lower bound on $\Vert \Delta u(t)\Vert_{L^2}$ close to the blow-up time (in the scaling invariance case with $a=0$), see details in \cite[Theorem 5]{BFM2010}, 
\begin{equation}\label{E:rate-lb}
 \Vert \Delta u(t)\Vert_{L^2(\R)}\geq\frac{C(u_0)}{\sqrt{t^*-t}},
\end{equation}
where $C(u_0)>0$ is some constant that depends on the initial data. A further conjecture on the blow-up rate in the critical case (pure quartic dispersion) is given in \cite[Conjecture 8]{BFM2010}.

To investigate the blow-up rate, we numerically track the $L^\infty$~norm of the solution, 
see Fig.~\ref{F:GaussRate}. For that we examine the numerical time 
dependence of the factor $\lambda(t)$, from the ansatz \eqref{dyn1}, via the time dependence of 
$\|u(t)\|_{L^\infty}$. As discussed above, the decay of $A^\infty \equiv (\ln L)_\tau \searrow 0$ controls the blow-up rate (as well as the convergence to the blow-up profile). 
Specifically, in the $L^{2}$-critical case for a stable blow-up an algebraic decrease $L\propto 1/\tau$ is expected, 
which implies that near the blow-up time, the time dependence is 
\begin{equation}
\lambda(t) \propto (t^{*}-t)^{1/3},
\label{Lcrit}
\end{equation}
whereas an exponential decrease $\lambda\propto e^{-\beta \tau}$ is 
expected for a stable blow-up in the $L^2$-supercritical case, leading to 
\begin{equation}
\lambda(t)\propto (t^{*}-t)^{1/4}.
\label{Lsupcrit}
\end{equation}
Using the rescaling \eqref{dyn1} and $\alpha=8$, we obtain 
\begin{equation}\label{E:infty-norm}
\|u(t)\|_{L^\infty} = \frac1{\lambda(t)^{1/2}} \|U(\tau)\|_{L^\infty} \cong \frac{1}{\big( (t^*-t)^{1/3}\big)^{1/2}} = \frac1{(t^*-t)^{1/6}},
\end{equation}
and 
\begin{equation}\label{E:L2-norm}
\|\partial_x^2 u(t)\|_{L^2} = \frac{c}{\lambda(t)^{2}} \cong \frac{1}{\big( (t^*-t)^{1/3}\big)^{2}} = \frac1{(t^*-t)^{2/3}}.
\end{equation}

\begin{figure}[!htb]
\includegraphics[width=0.45\hsize,height=0.32\hsize]{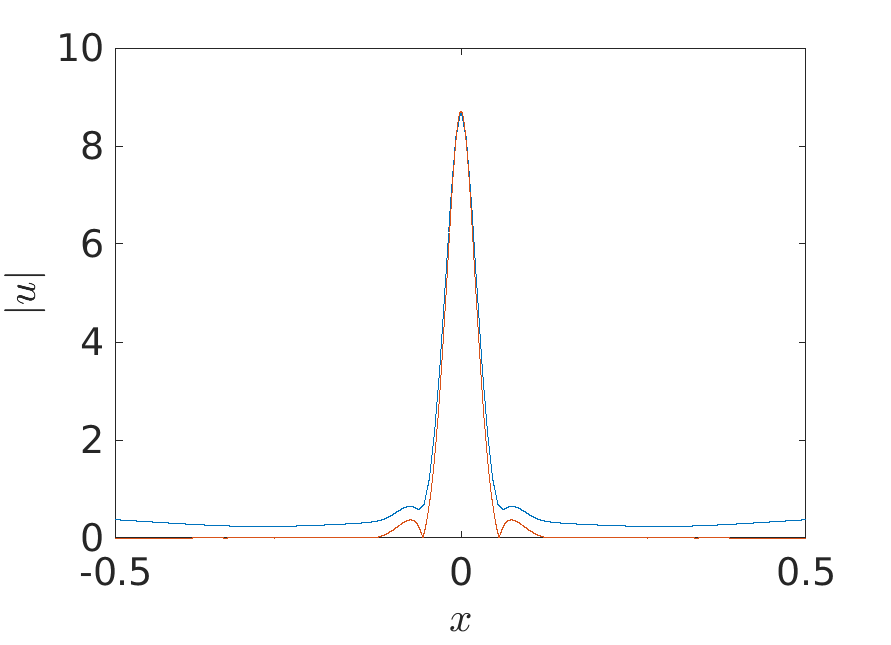}
\includegraphics[width=0.45\hsize,height=0.32\hsize]{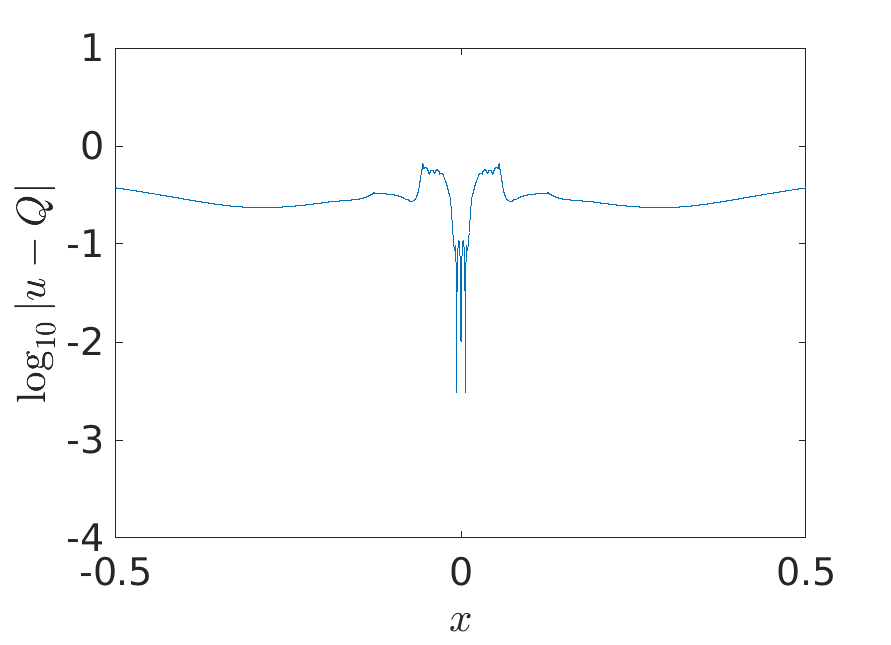}
\caption{\footnotesize {Blow-up profile and the fitting error for the solution of \eqref{E:explicit1}, $\alpha=8$, $a=0$, with $u_0=1.7 e^{-x^2}$, evolution of which is shown in Fig.~\ref{F:GaussRate}. 
Left: profile $\|u(t)\|_{L^\infty}$ (blue) at time $t_m: t^*-t_m = 1.8028\times 10^{-5}$, fitted to a rescaled ground state $Q$ (red) from \eqref{E:Q2}. Right: difference on a log scale between the solution and the fitted ground state.}}
\label{F:GaussProfile}
\end{figure}

Confirmation towards the rate in \eqref{E:infty-norm} can be seen on the right of 
Fig.~\ref{F:Soliton11_A-2} instead, whereas  
Fig.~\ref{F:GaussRate} does not provide a conclusive answer. 
This rate appears to be different from the conjectured rates in \cite{BFM2010}.

\subsubsection{Asymptotic profile of self-similar blow-up}\label{S:profile-critical}
Recalling the example of Gaussian data discussed in Fig.~\ref{F:GaussRate}, 
we show the asymptotic blow-up profile on the left of Fig. \ref{F:GaussProfile}.
In our simulations we use the profile $Q=Q^{(0)}$ from \eqref{E:Q2}, which is 
the same as \eqref{E:groundstate}. We rescale it and fit with the final 
computational state to check the matching of the asymptotic profile in the case of 
the finite time blow-up solution. The fitting is 
done as follows: the maximum of the modulus of the solution 
to equation \eqref{biNLS} is divided by the maximum of $Q$ to give 
according to \eqref{Qscaling} the value of $b^{1/\alpha}$. With this 
value of $b$ relation \eqref{Qscaling} gives $Q_{b}$ for the computed $Q$.

\smallskip

$\scriptstyle\blacklozenge$ ($a=0$) 
When $a=0$ (hence, scaling invariance holds), we are able to fit the final computational state with the (numerical) rescaled ground state profile of $Q^{(0)}$, 
for the initial datum $u_0 = 1.1Q^{(0)}$, we omit the figure for brevity, as it is similar to the next example. 

\begin{figure}[!htb]
\includegraphics[width=0.45\hsize,height=0.3\hsize]{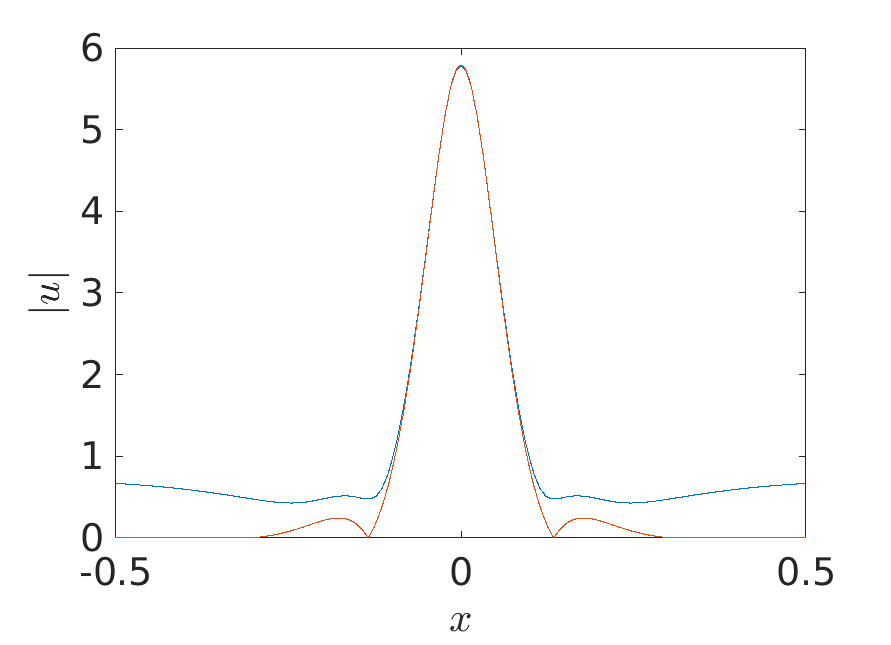}
 \includegraphics[width=0.45\hsize,height=0.3\hsize]{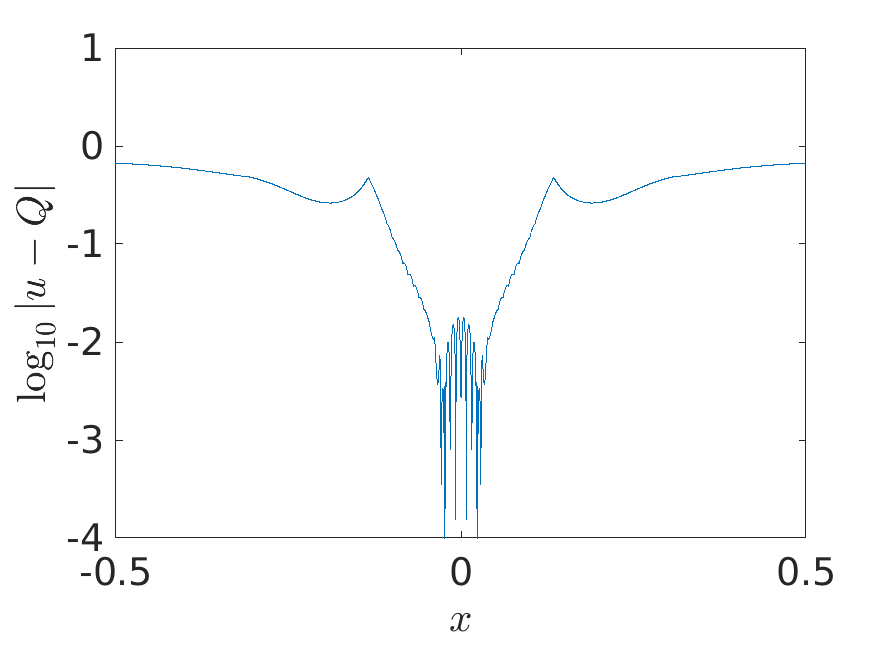}\\
\includegraphics[width=0.45\hsize,height=0.33\hsize]{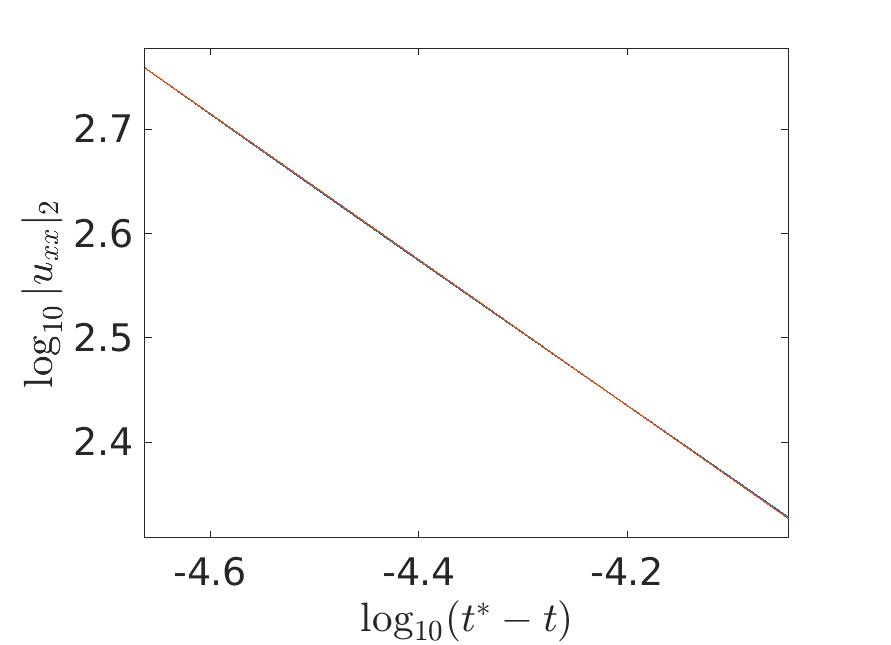}
 \includegraphics[width=0.45\hsize,height=0.33\hsize]{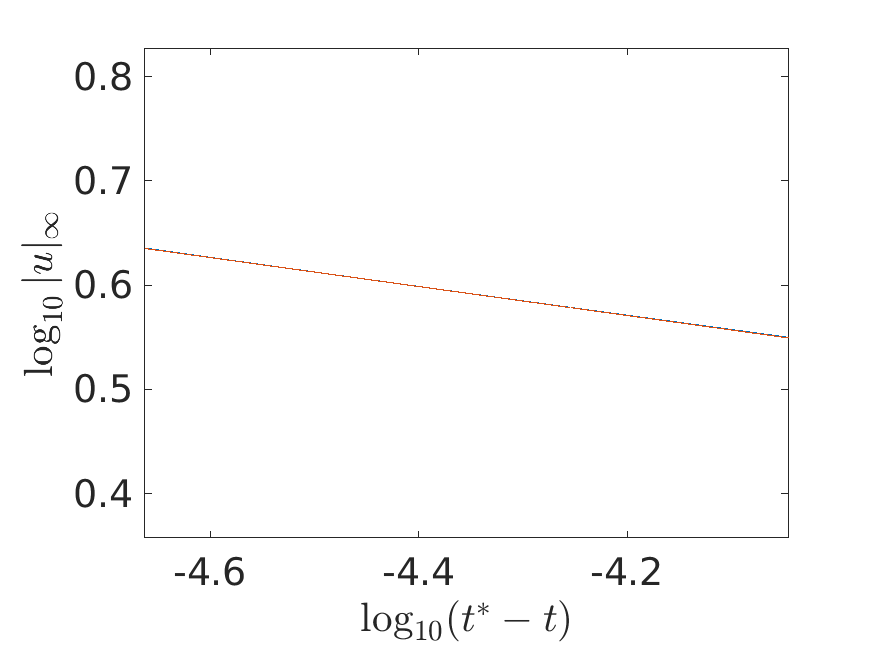}
\caption{\footnotesize
{Blow-up profile for the solution of \eqref{biNLS}, $\alpha=8$, $a = -2$, with $u_0=1.1\, Q^{(-2)}$ and fitting it with the rescaled ground state $Q^{(0)}$.
The solution blows up at $t^* =  0.19177$. 
Top left: Profile $|u(t)|$ at time $t_m: t^*-t_m = -4.1132\times 10^{-6}$ (blue) fitted to a re-scaled ground state $Q^{(0)}$ (red). 
Top right: Difference on a log scale between the solution and the fitted ground state.  
Bottom left: Blow-up of $\|u_{xx}\|_{L^2}$ experimentally fitted to the slope $0.7$. 
Bottom right: Blow-up of $\|u\|_{L^\infty}$ experimentally fitted to the slope of $0.139$.
}}
\label{F:Soliton11_A-2}
\end{figure}
\smallskip

$\scriptstyle\blacklozenge$ ($a\neq 0$)
We take $a=-2$ and compute the corresponding ground state solution $Q^{(-2)}$ (which has no oscillations in this case, since $a<-\sqrt 2$) from \eqref{E:1dGS}. The mass of this ground state is given in Table \ref{T:1}.
Then we consider perturbations $u_0 = A\, Q^{(a)}$ for different $A$. 
In Fig. ~\ref{F:Soliton11_A-2} we show the blow-up solution with the initial condition $u_0 = 1.1 Q^{(-2)}$, where we plot the solution and the fitting by the rescaled ground state $Q^{(0)}_{b'}$ at the blow-up time (namely, the last computable time before the blow-up with the step time difference on the order of $10^{-6}$) on the left plot and the error in the log scale on the right. The bottom row shows our fitting of the blow-up rate on the log scale, where we fit $(t^*-t)$ vs. $\|u_{xx}\|_{L^2}$ or $\|u\|_{L^\infty}$ norms close to the blow-up time. 
We obtain the power for the rate $\|u_{xx}\|_{L^2}$ to be around 0.7, 
and for  $\|u\|_{L^\infty}$ around 0.139, which provides some 
positive confirmation towards \eqref{E:infty-norm} and 
\eqref{E:L2-norm}. (We do stress that it is computationally very challenging 
to obtain more refined and more 
accurate rates in a double precision approach.)


We next consider the equation \eqref{E:explicit1} with $a=1$ (thus, the ground state $Q^{(1)}$ from \eqref{E:1dGS} has more oscillation than the scaling-invariant ground state $Q^{(0)}$ in the pure quartic case, $a=0$) and we study the blow-up behavior in that case. Note that the mass of the ground state $Q^{(1)}$ is smaller than that in the pure quartic case, $M[Q^{(1)}] < M[Q^{(0)}]$, see Table \ref{T:1}.

\begin{figure}[!htb]
\includegraphics[width=0.45\hsize,height=0.3\hsize]{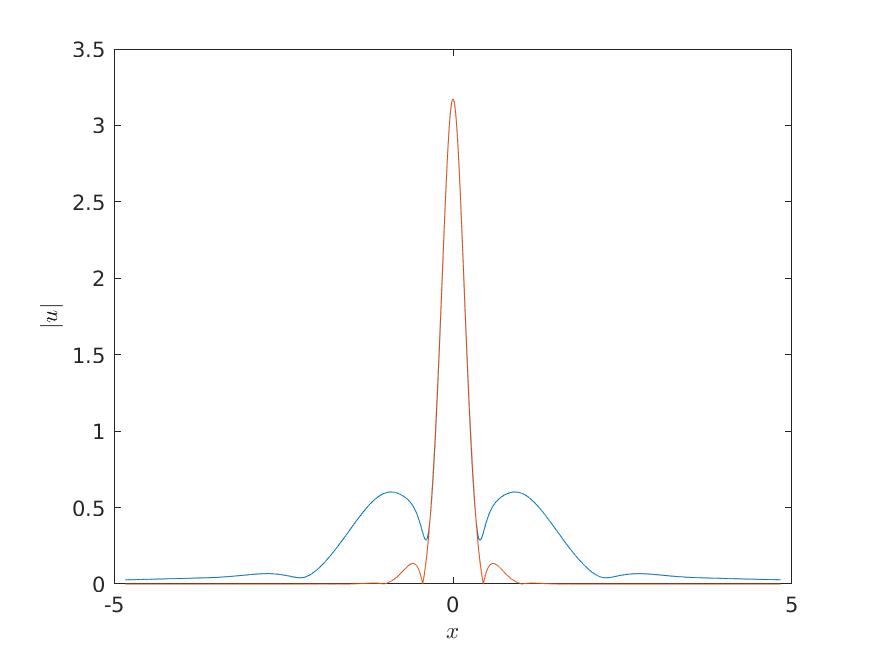}
\includegraphics[width=0.45\hsize,height=0.3\hsize]{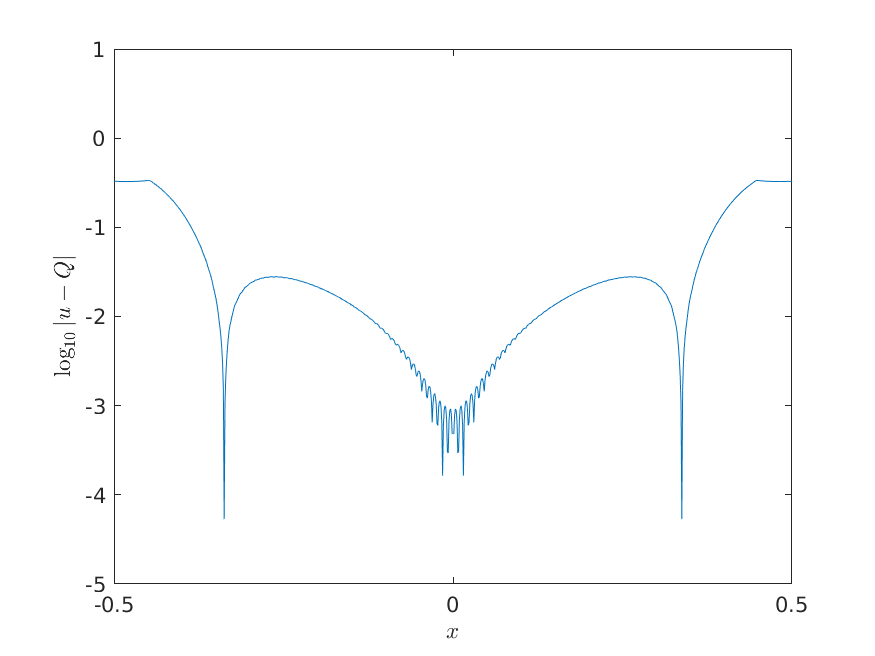}\\
\includegraphics[width=0.45\hsize,height=.33\hsize]{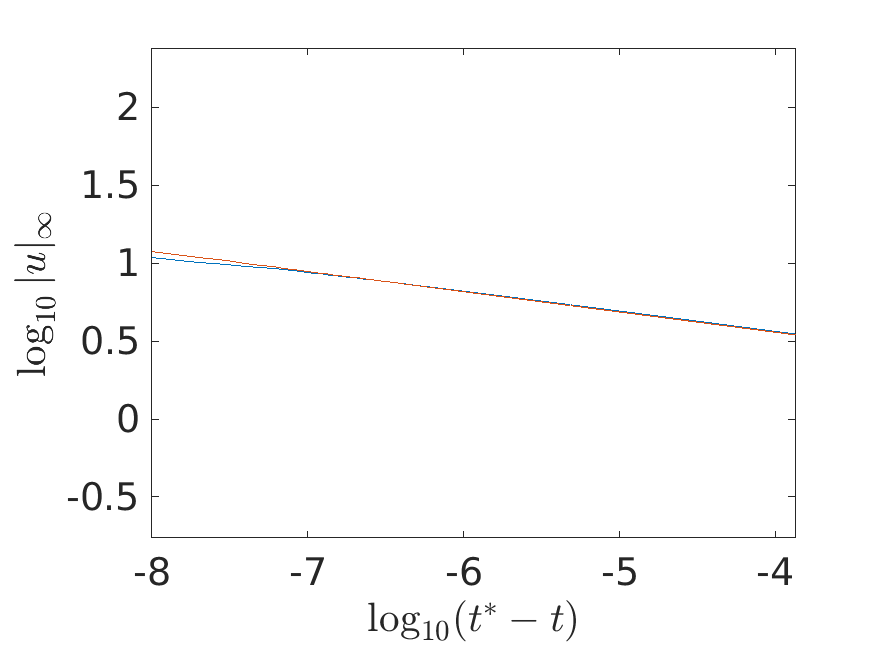}
\includegraphics[width=0.45\hsize,height=.33\hsize]{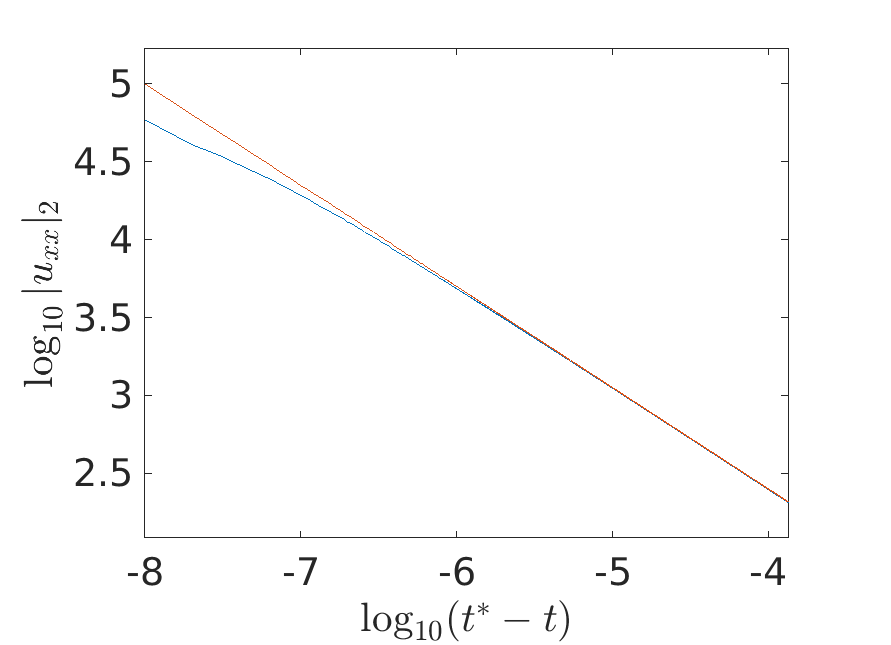}\\
\caption{\footnotesize Blow-up profile and fitting with the rescaled ground state for the solution of \eqref{biNLS} with $a = 1$ and  Gaussian initial data $u(x,0)=1.7 e^{-x^2}$. The solution blows up at $t^* =  0.021933$.  Snapshot of the profile at  $t_m: t^*-t_m = 1.3\times 10^{-6}$.
Top left: Profile $|u(t)|$ at time $t_m$ (blue) fitted to a rescaled ground state $Q^{(0)}$ (red). 
Top right: Difference on a log scale between the solution and the fitted ground state. 
Bottom left: Blow-up of $\|u_{xx}\|_{L^2}$ experimentally fitted to the slope 
$0.65$. Bottom right: Blow-up of $\|u\|_{L^\infty}$ experimentally fitted to the slope of $0.13$.}
\label{F:Gauss_A+1}
\end{figure}

In Fig. \ref{F:Gauss_A+1} instead of a slightly perturbed ground state data (which has similar results), we show the behavior of the Gaussian initial condition $u_0 = 1.7 e^{-x^2}$, which has mass above $M[Q^{(0)}]$. We show that it blows up in finite time with the profile converging to the rescaled ground state $Q^{(0)}$, see top left plot in Fig. \ref{F:Gauss_A+1}. The difference between the solution (at the final computed time) and the rescaled soliton $Q^{(0)}$ is shown in the top right plot. We also show the fitted rates in the bottom row, as in the case of $a=-2$, we get similar values.  
\smallskip

\begin{figure}[!htb]
\includegraphics[width=0.45\hsize,height=0.33\hsize]{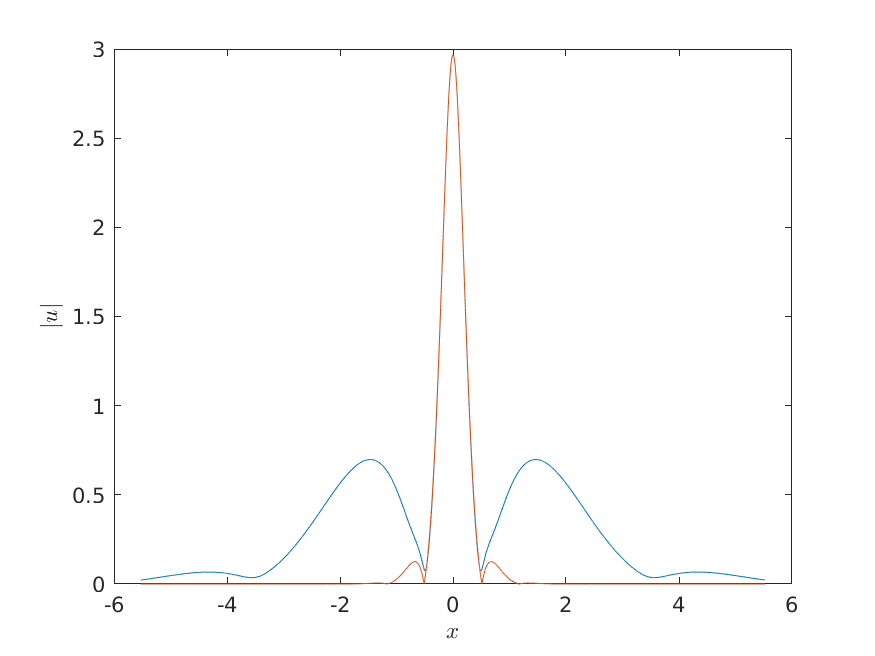}
 \includegraphics[width=0.45\hsize,height=0.33\hsize]{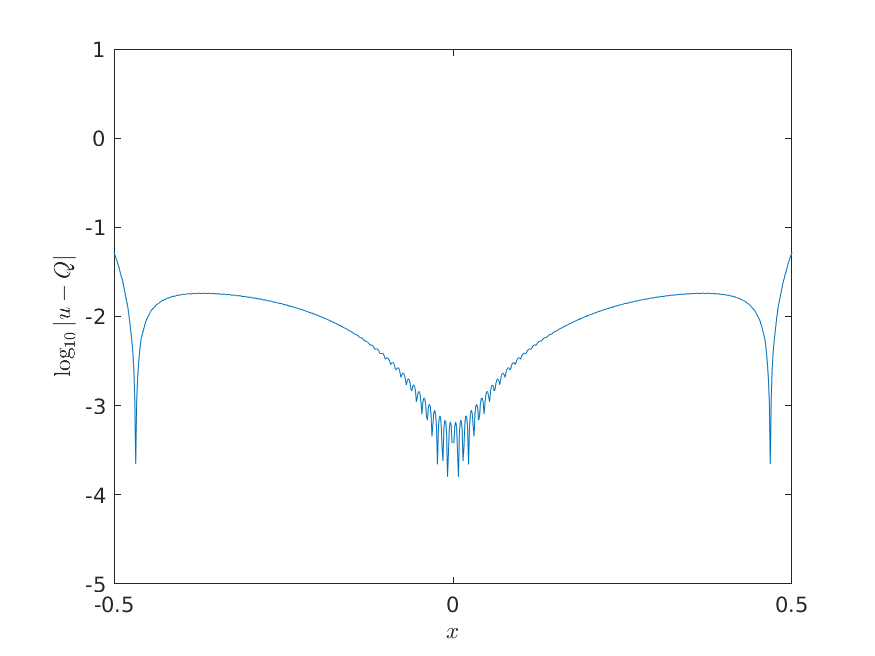}
\caption{\footnotesize Blow up profile and fitting with the rescaled ground state for the solution of \eqref{biNLS} with $a = 1$, $\alpha=8,$ and  initial data $u(x,0)=1.45 \,\sech(x)$. 
Top left: Profile $|u(t)|$ at time $t_m$ (blue) fitted to a rescaled ground state $Q^{(0)}$ (red). 
Top right: Difference on a log scale between the solution and the fitted ground state. 
}
\label{F:Sech_A+1}
\end{figure}

In Fig. \ref{F:Sech_A+1} we consider hyperbolic secant initial condition $u_0 = A \, \sech \, x$ with the magnitude $A = 1.45$, so that the mass is above $M[Q^{(0)}]$. The plot on the left shows that the solution blows up in finite time with the profile converging to the rescaled ground state $Q^{(0)}$. 

\clearpage

\section{Dynamics of solutions in the supercritical case}\label{S:supercritical}

When $\alpha>8$, in the supercritical case of \eqref{biNLS} in 1d, 
blow-up is also expected for some sufficiently localized data, see \cite{Fibich2002}. Similarly to the critical case, in dimensions two and higher the existence of finite time blow-up was proved in \cite{BL2017} for certain cases of \eqref{biNLS}. 
\begin{figure}[!htb]
\includegraphics[width=0.45\hsize,height=0.28\hsize]{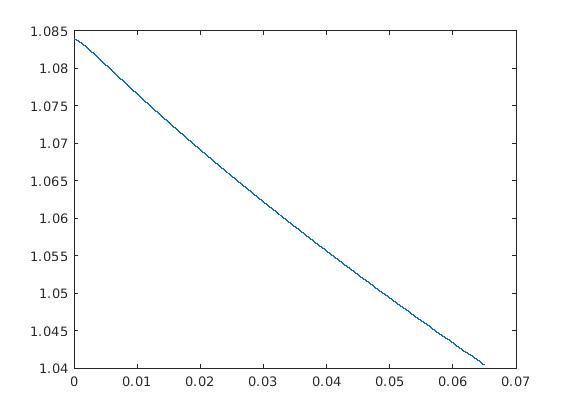}
\includegraphics[width=0.45\hsize,height=0.28\hsize]{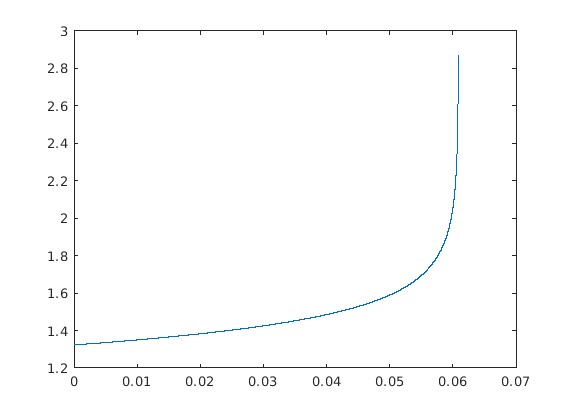}
\caption{\footnotesize {Supercritical case $\alpha = 10$. Dichotomy behavior in the pure quartic bi-NLS \eqref{E:explicit1}, 
$a=0$, for $u_0 = A\, Q$ with $A=0.9$ (left) and $A=1.1$ (right).}}
\label{F:super1}
\end{figure}

In the supercritical case, there is no obvious threshold as in the 
critical case, at least for the standard NLS, however, in 2d and 
higher the results on the long term behavior of solutions in the spirit of the original dichotomy of Holmer-Roudenko \cite{HR2007, HR2008} (global existence and scattering vs. blow-up) have been shown in \cite{BL2017}, see also \cite{Dinh2021}. 

In the $L^2$-supercritical case, we consider the scaling-invariant case first ($a=0$) and investigate the time evolution of various perturbations of ground states 
$u_0=A\, Q^{(0)}$ with $A \approx 1$. We find that the solutions with $A<1$ scatter and with $A>1$ blow up in finite time, which we show in Fig.~\ref{F:super1} for $\alpha=10$.
We note that since in this supercritical case we did not observe any branching in the energy-mass curve, see right bottom plot in Fig. \ref{F:ME-crit+supercrit}, we do not expect any changes in dichotomy behavior of solutions, which is confirmed in Fig. \ref{F:super1}.

For other types of data and comparison, we record the mass of the ground state $Q^{(0)}$ (in the case of $\alpha=10$), which  is
$M[Q^{(0)}] = 2.75816089721144$.

\subsection{Rate and Profile}\label{S:rate-supercritical}

To confirm the blow up rate, 
we use the rescaling \eqref{dyn1} to deduce 
\begin{equation}\label{E:infty-norm-super}
\|u(t)\|_{L^\infty} = \frac1{\lambda(t)^{4/\alpha}} \|U(\tau)\|_{L^\infty} \cong \frac{1}{\big( (t^*-t)^{1/4}\big)^{4/\alpha}} = \frac1{(t^*-t)^{1/\alpha}}.
\end{equation}
\begin{figure}[!htb]
\includegraphics[width=0.43\hsize,height=0.3\hsize]{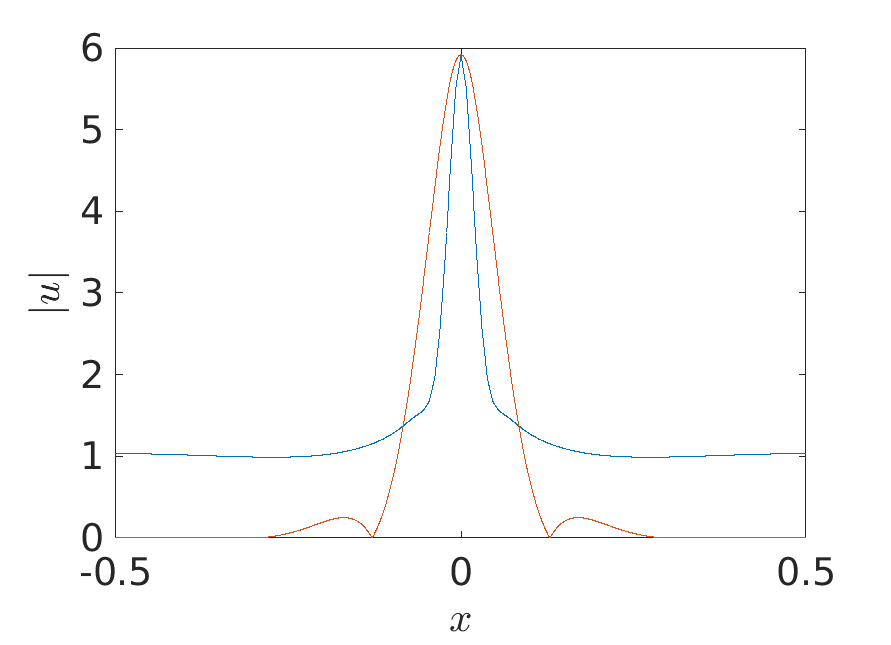}
 \includegraphics[width=0.43\hsize,height=0.3\hsize]{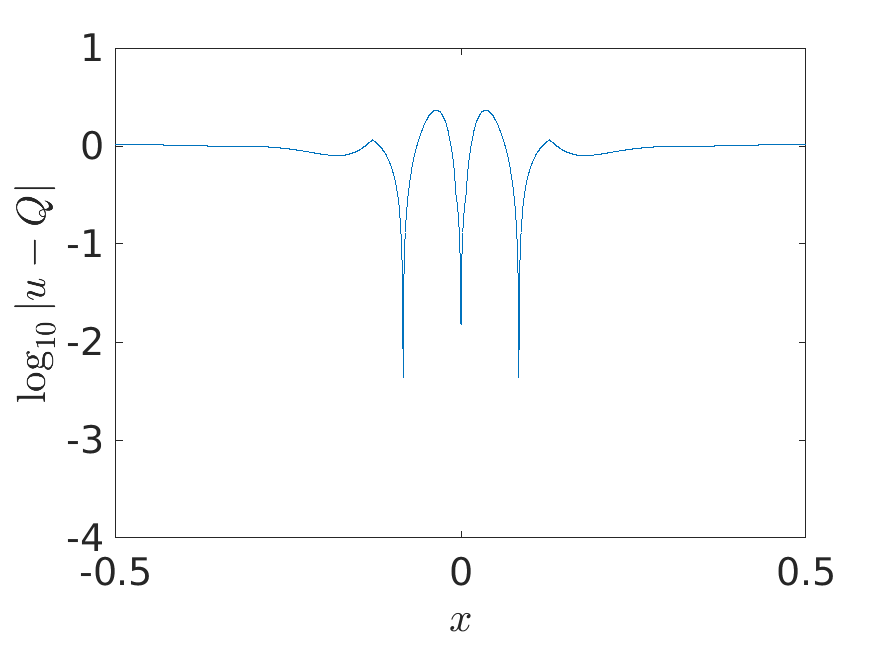}\\
\includegraphics[width=0.43\hsize,height=0.3\hsize]{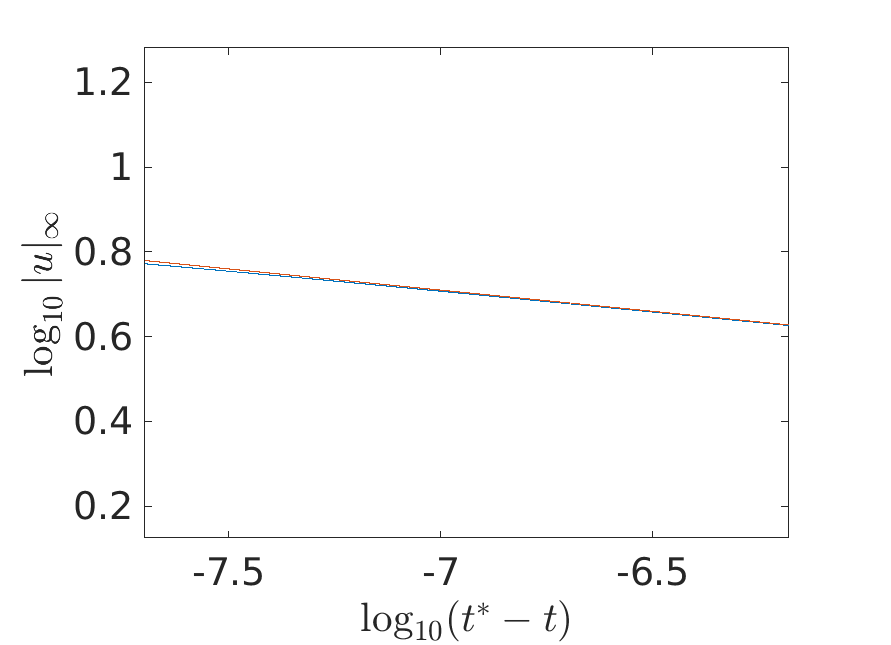}
\includegraphics[width=0.43\hsize,height=0.3\hsize]{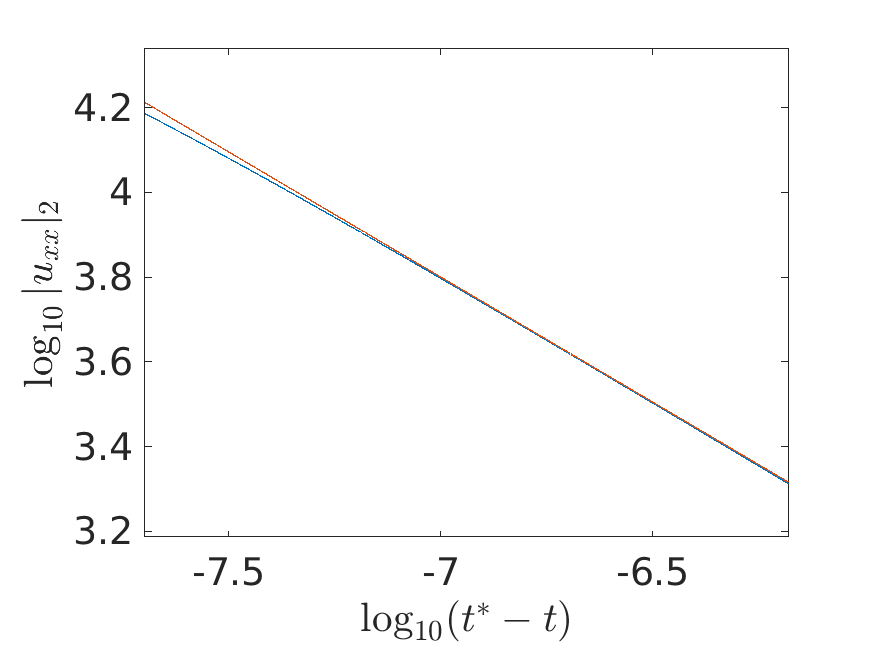}\\
\caption{\footnotesize Supercritical case $\alpha = 10$. Blow up profile and fitting with the rescaled ground state for the solution of \eqref{biNLS} with $a = 1$ and  Gaussian initial data $u(x,0)=1.7 e^{-x^2}$. The solution blows up at $t^* =  0.00375322$.  Snapshot at  $t_m: t^*-t_m = 2.0\times 10^{-8}$.
Clockwise from top left: Profile $|u(t)|$ at time $t$ (blue) fitted to a re-scaled ground state (red, given by the asymptotic solution for $a = 0$); Difference on a log scale between the solution and the fitted ground state;  Blow-up of $\|u_{xx}\|_{L^2}$ experimentally fitted to rate $0.59$; Blow-up of $\|u\|_{L^\infty}$ experimentally fitted to rate of $0.1$. The case with negative $a$ does not provide qualitatively new results.}
\label{F:Gauss_A+1super}
\end{figure}

In particular, in the power case $\alpha=10$, the rate is $1/10$, which we confirm in our fittings in Fig.~\ref{F:Gauss_A+1super}. We note that in this supercritical case, it is easier to check the rate numerically, as it converges to a specific profile very fast (unlike the critical case, so these features are similar to the standard NLS equation).  

As far as the profile is concerned, one can see in the top left plot of Fig. \ref{F:Gauss_A+1super} that the profile is different from the ground state in this case $Q^{(0)}$ (with $\alpha=10$). This is consistent with the standard NLS equation. 

To check further on rates and profiles of the solution, we considered super-Gaussian and sech initial data, and obtained similar results on the rate and profiles.  

The principal difference of the supercritical case to the critical case is that the asymptotic profile and the rescaled ground state $Q^{(0)}$ do not coincide, while the rate is easier to track and fit, and it is consistent with the theoretical prediction as in \eqref{E:infty-norm-super}. 

\section{Conclusion}
In this paper, we have presented a detailed numerical study of 
solutions to the general 4th order (bi-harmonic) NLS equation \eqref{biNLS} in 1d. 
Ground state 
solutions have been constructed numerically, and 
stability of ground states has been investigated. 
In the subcritical cases ($\alpha = 4,6$) we found that there are two branches of ground states, leading to a {\it stable} and an {\it unstable} branches of ground state solutions (which are determined by their energy vs. mass dependence). 
In the critical case, we found that a richer dynamics of solutions than a dichotomy in a standard NLS equation holds: smaller amplitude solutions tend to disperse; solutions which are close to the mass of either of the ground states (a non-scale invariant case ground state $Q^{(a)}$ and the scaling-invariant ground state $Q^{(0)}$) do not disperse but may approach a different branch of ground states, which is a {\it new} phenomenon. 
 
In the critical case the branching also occurs, and besides the previous two behaviors, there is also a blow-up in finite time.  However, a typical threshold for the scattering vs. blow-up as it is in the standard NLS given by a ground state solution, does not work here: in the non-invariant cases, there is a gap where solutions will neither disperse to zero nor blow-up; instead they will approach  a different (stable branch) ground state in an oscillatory manner. This is a new occurrence. 

We conjecture that the blow-up in finite time occurs in a self-similar manner with the profile given by the scale-invariant $Q^{(0)}$ in any case of $a <\sqrt b$. This is new.  
In the supercritical case Schwartz class initial data of 
sufficient mass are shown to blow-up in finite time with a 
self-similar blow-up mechanism.

\section*{Competing interests declaration}
All authors that have contributed to the  submission declare that 
they have no competing interests. 




\bibliographystyle{abbrv}
	
\bibliography{references}

\end{document}